\documentclass[11pt]{article}
\usepackage{amsbsy}
\usepackage{stmaryrd}
\usepackage{bm}
\usepackage{ulem}
\usepackage{tikz}
\usepackage{caption}
\usepackage{subcaption}
\usepackage{amssymb,amsmath,amsfonts,amsthm}
\usepackage{amssymb,amsmath,amsfonts,amsthm}
\usepackage{color}

\newtheorem {Proposition}{Proposition}[section]
\newtheorem {Lemma}{Lemma}[section]

\newtheorem {Theorem}{Theorem}[section]

\usepackage{amssymb,amsmath}
\usepackage{authblk}
\usepackage{epsfig}
\usepackage{enumerate,graphicx}
\usepackage{amsfonts}
\usepackage{graphicx}
\DeclareGraphicsExtensions{.pdf, .jpg, .png , .bmp, .ps}
\def\build#1_#2^#3{\mathrel{\mathop{\kern 0pt#1}\limits_{#2}^{#3}}}
\def\cvl{\build{\ \longrightarrow\ }_{\nti}^{{\mathcal L}}}

\def\nti{n \rightarrow \infty}
\DeclareMathOperator{\R}{\mathbb{R}}
\DeclareMathOperator{\X}{\textbf{X}}
\DeclareMathOperator{\Z}{\textbf{Z}}

\DeclareMathOperator{\Y}{\textbf{Y}}

\DeclareMathOperator{\E}{\textbf{E}}

\def\nti{n \rightarrow \infty}

\def\Var{\mathop{\rm Var}\nolimits}%

\newcommand\ind{{ {{1}}\hspace{-0,8mm}{\mathrm I}}}

\newcommand{\cov}{\text{Cov}}

\begin{document}
\title{A model-free  Screening procedure}

\author[1]{J. Dedecker}
\affil[1]{Université Paris Cité, UMR CNRS  8145, Laboratoire MAP5, 75006, Paris, France.}
\author[2]{M.L. Taupin}
\affil[2]{Université Paris-Saclay, UMR CNRS 8071, Univ Evry, Laboratoire de Mathématiques et Modélisation d'Evry (LaMME), 91037, Evry, France.}
\author[2]{A.S. Tocquet}

\date{}

\maketitle



\begin{abstract}
In this article, we propose a generic screening method for selecting explanatory variables correlated with the response variable $Y$. We make no assumptions about the existence of a model that could link $Y$ with a subset of explanatory variables, nor about the distribution of the variables. Our procedure can therefore be described as ``model-free" and can be applied in a wide range of situations. In order to obtain precise theoretical guarantees (Sure Screening Property and control of the False Positive Rate), we establish a Berry-Esseen type inequality for the studentized statistic of the slope estimator. We illustrate our selection procedure using two simulated examples and a real data set.

\end{abstract}


\smallskip


\tableofcontents

\section{Introduction, objectives and notations}
\setcounter{equation}{0}
Let $(Y,X^{(1)},\cdots,X^{(p)})$ be a random vector of $\mathbb{R}^{p+1}$ where $Y$ is the response variable and $(X^{(1)},\cdots,X^{(p)})$ is the vector of explanatory variables. Assume that all the variables are square integrable, and  let $\mathcal{M}^*$ be the set of indexes $j\in\llbracket 1,p\rrbracket$ for which $X^{(j)}$ is correlated with $Y$, that is 
\begin{align}
\label{M*0}
\mathcal{M}^*=\{j\in\llbracket 1,p\rrbracket \, ; \, \mbox{\cov}(X^{(j)},Y)\not=0\}\, . 
\end{align}

Suppose we have $n$ observations from $n$ independent copies $(Y_i,X_i^{(1)},\ldots,X_i^{(p)})_{1\leq i\leq n}$ of the vector 
$(Y,X^{(1)},\ldots,X^{(p)})$. As in \cite{Fan2008}, we will propose a statistical procedure for selecting from the sample $(Y_i,X_i^{(1)},\ldots,X_i^{(p)})_{1\leq i\leq n}$ a random set $\widehat{\mathcal{M}}$ of indices, by considering only some marginal statistics involving the variables $(Y_i, X^{(j)}_i)_{1 \leq i \leq n}$ for all  $j\in\llbracket 1,p\rrbracket$. This method therefore belongs to the family of procedures known as ``screening", which do not require multivariate models adjustment (sometimes costly from a computational point of view, or not feasible for identifiability reasons). We shall see that the marginal statistics that we propose are particularly easy to compute (not  requiring the implementation of any minimisation procedures), so that our selection  rule can be effectively implemented even when the number of explanatory variables $p$ is very large.  The random set $\widehat{\mathcal{M}}$ is constructed  so as to satisfy the Sure Screening Property (SSP), i.e. 
$$\lim_{\nti}\mathbb{P}(\mathcal{M}^*\subset \widehat{\mathcal{M}})=1 \, .$$ 
More precisely, as in \cite{Fan2008}, we will give a non asymptotic control of the quantity $\mathbb{P}(\mathcal{M}^*\subset \widehat{\mathcal{M}})$. Furthermore, as in \cite{CoxScreening2012}, our statistical procedure will also enable us to control the False Positive Rate (FPR). More specifically, given a rate $q \in (0,1)$ defined prior to the study, our procedure (which depends on $q$) will satisfy the following property 
\begin{equation*}
\lim_{\nti} \left\vert \mathbb{E}\Bigg[\frac{|{\mathcal{M}^*}^c\cap \widehat{\mathcal{M}}|}{|{\mathcal{M}^*}^c|}\Bigg]-q \right\vert = 0 \, .
\end{equation*}
More precisely, as in \cite{CoxScreening2012}, we will give a non-asymptotic control of the quantity
$$
\left\vert \mathbb{E}\Bigg[\frac{|{\mathcal{M}^*}^c\cap \widehat{\mathcal{M}}|}{|{\mathcal{M}^*}^c|}\Bigg]-q \right\vert \, .
$$


To achieve these two objectives (SSP and  control of the FPR), we will not make any assumptions about the existence of an underlying model that would link $Y$ with a subset of the explanatory variables. Nor will we make any assumptions about the distribution of the variables $Y$ or $X^{(j)}$'s. 
These variables may be continuous, binary, discrete, etc. Lastly we will not assume that the random variables indexed in ${\mathcal M}^*$ are globally independent of the variables indexed in  ${\mathcal{M}^*}^c$ as is required in some papers about screening (see for instance Assumption 8 in 
\cite{CoxScreening2012}).
We can therefore say that our procedure is ``model-free" and ``robust". In particular, we do not assume that the variable $Y$  has a conditional distribution given the $X^{(j)}$'s belonging to the generalized linear model, which is a fundamental difference between our work and that proposed in \cite{FSW2009} and \cite{ScreeningGLM2010}, for example. We therefore propose a unique procedure, valid only under reasonable moment conditions, and offering precise theoretical guarantees. 

It should be noted that model-free methods have already been proposed: in particular, we can mention the article \cite{Panetal2019} (see also some of the references given therein) in which the authors propose a procedure that does not require moment assumptions - and which is, in this sense, more robust than ours. However, the approach proposed in this article is very different from ours. First, the authors assume that $\log(p)$ is of the order $o(n)$, which is a rather unnatural assumption for a screening method – an assumption that is unnecessary for our procedure. Second, their selection rule is based on prior information about the Ball correlation coefficients $\rho_i$, which are unknown quantities (see the definition of $\widehat{{\mathcal A}}_n^*$ and the condition (C$_1$) page 930). This last point  can be problematic in practice.

Our selection rule, which we will describe in detail in the following paragraphs, had already been proposed in \cite{DGT2025}, but without theoretical guarantees. In \cite{DGT2025}, we showed through a simulation study that the selection rule worked well when  $Y$ is a binary variable, whose conditional distribution given a subset of the $X^{(j)}$'s is given by a logistic regression model. In this article, we will give two other examples of simulations (continuous outputs and discrete outputs) where $Y$ does not come from a generalized linear model, and for which our procedure continues to perform well. 

Let us now describe our selection rule.
For any $j\in\llbracket 1,p\rrbracket$, assume that $\Var(X^{(j)})>0$ and let
\begin{align}
\label{bj_tauj}
 \tau_j=\frac{\cov(X^{(j)},Y)}{\Var(X^{(j)})}\quad \mbox{ and }\quad b_j=\mathbb{E}(Y)-\tau_j\mathbb{E}(X^{(j)}).
\end{align} 
Then we can write,  
\begin{align}
Y= b_j+\tau_j X^{(j)} +\varepsilon^{(j)} 
\label{model}
\end{align}
where $\mathbb{E}\left(\varepsilon^{(j)}\right)=0$ and $\mathbb{E}\left(\varepsilon^{(j)}X^{(j)}\right)=0$.
From the definition of the $\tau_j$'s, the set $\mathcal{M}^*$ defined in \eqref{M*0} can be written as 
\begin{align}
\label{M*}
\mathcal{M}^*=\{j\in\llbracket 1,p\rrbracket \, ; \, \tau_j\not=0\}.
\end{align}
Consider now $n$ independent copies $(Y_i,X_i^{(1)},\ldots,X_i^{(p)})_{1\leq i\leq n}$ of the vector 
$(Y,X^{(1)},\ldots,X^{(p)})$.
Then for all $i=1,\ldots,n$ and $j\in\llbracket 1,p\rrbracket$, we define the errors $\varepsilon_i^{(j)}$'s via the equations 
\begin{align*}
Y_i&=b_j+\tau_jX_i^{(j)}+\varepsilon_i^{(j)}. 
\end{align*}
Let now 
$$\Y=(Y_1,\ldots,Y_n)^T, \qquad \X^{(j)}=(X_1^{(j)},\ldots,X_n^{(j)})^T
$$
and $$\overline{\Y}=\frac{1}{n}\sum_{i=1}^n{Y_i}, \qquad
\overline{\X}^{(j)}=\frac{1}{n}\sum_{i=1}^n{X_i^{(j)}}.$$
For any $j\in\llbracket 1,p\rrbracket$, we then consider the usual estimators of $\tau_j$ and  $b_j$, 
\begin{align}
\label{tauchapoj_bchapoj}
\hat \tau_j=\frac{\cov_n(\Y,\X^{(j)})}{\Var_n(\X^{(j)})}=\frac{\sum_{i=1}^n(X_i^{(j)}-\overline{\X}^{(j)})(Y_i-\overline{\Y})}{\sum_{i=1}^n\left(X_i^{(j)}-\overline{\X}^{(j)}\right)^2} \quad \mbox{ and }\quad \hat b_j=\overline{\Y}- \hat \tau_j \overline{\X}^{(j)}.
\end{align}
It follows that, for $j\in\llbracket 1,p\rrbracket$,
\begin{align*}
\hat \tau_j-\tau_j=\frac{n^{-1}
\sum_{i=1}^n(X_i^{(j)}-\overline{\X}^{(j)})\varepsilon_i^{(j)}
}{n^{-1}\sum_{i=1}^n{\left(X_i^{(j)}-\overline{\X}^{(j)}\right)^2}}.
\end{align*}
Let now $\hat\varepsilon_i^{(j)}$'s be the residuals of the $j$-th linear regression
\begin{align*}
 \hat\varepsilon_i^{(j)}=Y_i-\hat b_j-\hat\tau_jX_i^{(j)} \, ,
\end{align*} 
and let
\begin{eqnarray*}
\widehat{v}_j &=&  \frac{n^{-1} \sum_{i=1}^n{\left(X_i^{(j)}-\overline{\X}^{(j)} \right)^2 
    \left( \widehat \varepsilon^{(j)}_i \right)^2}}
{
\left( n^{-1}\sum_{i=1}^n {\left(X_i^{(j)}-\overline{\X}^{(j)}\right)^2} \right)^2}.
\end{eqnarray*}
According to \cite{W} (see also \cite{DGT2025}),  assuming 
$\mathbb{E}\big({X^{(j)}}^4\big)<\infty$ and $\mathbb{E}\left({Y}^4\right)<\infty$ for any $j\in\llbracket 1,p\rrbracket$,
we get 
$$\sqrt{n}\frac{\widehat{\tau}_j - \tau_j}{\sqrt{\widehat{v}_j}} \cvl \mathcal{N}(0,1).$$
This implies that a reasonable screening rule is 
\begin{align}
\label{Mchapo}
\widehat{\mathcal M}=\left\{ j\in\llbracket 1,p\rrbracket \, ; \,\frac{\sqrt{n}|\widehat{\tau}_j|}{\sqrt{\widehat{v}_j}} \geq \gamma  \right\}
\end{align}
with $\gamma=\Phi^{-1}(1-q/2)$, where $\Phi$ is the cumulative distribution function of the standard Gaussian random variable and 
$q$ is the false positive rate that we choose. 

Our purpose is to show that this rule allows us to obtain the SSP as well as a satisfactory control of the FPR. To do this, we need to establish a new probabilistic result: a non-asymptotic  Berry-Esseen type inequality for the studentized statistic
$$\sqrt{n}\frac{\widehat{\tau}_j - \tau_j}{\sqrt{\widehat{v}_j}} \, .$$
This inequality is of statistical interest in itself, since it is a robust result valid without the structural assumption that $\varepsilon^{(j)}$ is independent of $X^{(j)}$ (by definition of  $\varepsilon^{(j)}$, the variables $\varepsilon^{(j)}$ and $X^{(j)}$ are only uncorellated).  Its proof is quite technical and requires the combination of many probabilistic tools (that are listed in Subsection \ref{5.1}). The main technical difficulty lies in replacing the residuals $\hat\varepsilon_i^{(j)}$'s in the expression of $\widehat{v}_j$ with the unobserved errors $\varepsilon_i^{(j)}$'s (see the comments after Proposition \ref{propositionattenduesansj}). 

This Berry-Esseen type inequality is given in Section \ref{sec2}, which also contains two other preliminary lemmas.   These results allow us to establish the theoretical properties of our screening method in Section \ref{sec3}. In Subsection \ref{sec3.1} we check the SSP, and in Subsection \ref{sec3.2} we give a non asymptotic control of the FPR. 
In Section \ref{Simulations}, we evaluate the performance of our screening method on different sets of simulations, and we present an application to a real data set. Finally, the proofs of our results are given in Section \ref{Proofs}. 

In the rest of the paper, we denote by $\Phi$  the cumulative distribution function of the standard Gaussian random variable (as already introduced right after \eqref{Mchapo}). 

\section{Berry-Esseen bound for studentized statistic of the slope estimator}\label{sec2}
\setcounter{equation}{0}
Let  $n\in\mathbb{N}^*$ and $(X_1,Y_1),\ldots, (X_n,Y_n)$ be $n$ independent copies of a random  couple $(X,Y)$. 
We use the following notations
\begin{align}
\mu_X =\E(X), \quad\tilde{X}=X-\mu_X, \quad \tilde{X}_i=X_i-\mu_X, \quad
\sigma_{X}^2=\Var(X), \quad \overline{\X}=\frac{1}{n}\sum_{i=1}^n{X_i}, \label{notationsX} \\
\mu_Y =\E(Y), \quad \tilde{Y}=Y-\mu_Y , \quad \tilde{Y}_i=Y_i-\mu_Y, \quad  \sigma_{Y}^2=\Var(Y), \quad \overline{\Y}=\frac{1}{n}\sum_{i=1}^n{Y_i}, \label{notationsY} \\
\tau=\frac{\cov(Y,X)}{\Var(X)}, \quad \varepsilon=\tilde{Y}-\tau\tilde{X}, \quad \varepsilon_i=\tilde{Y}_i-\tau\tilde{X}_i,\quad \sigma_{\varepsilon}^2=\Var(\varepsilon).  
\label{notationsepsilon}
\end{align}
First, we provide two intermediate lemmas that enable us to prove the main result of this section as well as the results of the following section. 
\begin{Lemma}
\label{prelim}
Let $K_1$ and $K_2$ be defined by 
\begin{align}
\label{kappa1_kappa2}
K_1(X,Y)=\frac{\sigma_{\varepsilon}\sigma_X}{
\sqrt{2\pi}\big\Vert\tilde{X}\varepsilon \big\Vert_2} \quad\mbox{ and }\quad
K_2(X,Y)=2+c_{\rm{\scriptscriptstyle BE}}\left[ \frac{ 
\big\Vert\tilde{X}\varepsilon \big\Vert_3^3
  }{ \big\Vert\tilde{X}\varepsilon \big\Vert_2^3  
  }  +\frac{2\left\Vert\varepsilon\right\Vert_3^3
  }{\sigma_{\varepsilon}^3} 
+\frac{2\big\Vert\tilde{X}\big\Vert_3^3}{\sigma_X^{3}}\right]
\end{align}
where $c_{\rm{\scriptscriptstyle BE}}$ is the numerical constant 
involved in the Berry-Esseen inequality (\ref{BE}). 
Then, 
\begin{equation}
\label{inegalite_Lemme}
\sup_{x\in \mathbb{R}}\left\vert\mathbb{P}\left(\frac{\sum_{i=1}^n
(X_i-\overline{\X})\varepsilon_i}
{
\sqrt{n \big\Vert\tilde{X}\varepsilon\big\Vert_2^2}
}\leq x \right)-\Phi(x)\right\vert\leq \frac{\log(n)}{\sqrt{n}}K_1(X,Y)+\frac{K_2(X,Y)}{\sqrt{n}} \, .
\end{equation}
\end{Lemma}
For the next lemma and the next proposition, we need to define $\hat\tau$ and $\hat \varepsilon_i$ by 
\begin{equation}\label{deftauhat} \hat \tau=\frac{\sum_{i=1}^n{(X_i-\overline{\X})(Y_i-\overline{\Y})}}{\sum_{i=1}^n{\left(X_i-\overline{\X}\right)^2}}, \quad \hat b=\overline{\Y}- \hat \tau \overline{\X} \quad\mbox{ and } \quad \hat \varepsilon_i=Y_i- \hat b - \hat\tau X_i\, .
\end{equation}
For readability reasons, the Lyapunov ratios involved in Lemma \ref{lemme2} and Proposition \ref{propositionattenduesansj} are defined in Section \ref{definitionG_ietH_i}. 
\begin{Lemma}
\label{lemme2}
Let $G_1$, $G_2$ and $G_3$ be defined by \eqref{termeG1}, \eqref{termeG2} and \eqref{termeG3}. 
Then
\begin{align}
\mathbb{P}&\left(\left\vert \frac{1}{n}\sum_{i=1}^n
{\left(X_i-\overline{\X}\right)^2\widehat\varepsilon_i^2}-\big\Vert \tilde{X}\varepsilon\big\Vert_2^2 \right\vert  \geq\frac{\big\Vert \tilde{X}\varepsilon\big\Vert_2^2}{2}\right) \leq 
\frac{G_1(X,Y)}{\sqrt{n}}+ \frac{4G_2(X,Y)}{n}+\frac{G_3(X,Y)}{n^{3/2}}.
\label{majoration_terme2_2_reste_prop}
\end{align}
\end{Lemma}
Note that the terms $G_1(X,Y), G_2(X,Y), G_3(X,Y)$ are finite as soon as the variables $X$ and $Y$ have a finite moment of order 6, and $\tilde X$, $\tilde X \varepsilon$, $\tilde X \tilde Y$ are not equal to zero almost surely. 

The following proposition is the main result of this section providing  a Berry-Esseen bound for the studentized statistic of the slope estimator $\hat \tau$. 
\begin{Proposition}
\label{propositionattenduesansj}
Let  $(X_1,Y_1),\ldots, (X_n,Y_n)$ be $n$ independent copies of $(X,Y)$. 
Let $\widehat{v}$ be defined by
\begin{equation}\label{defvhat}
\widehat{v} =  \frac{n^{-1} \sum_{i=1}^n\left (X_i-\overline{\X} \right )^2 
    \widehat \varepsilon_i^2}
{
\left( n^{-1}\sum_{i=1}^n {\left(X_i-\overline{\X}\right)^2} \right)^2}.
\end{equation}
Let $H_1$, $H_2$, $H_3$, $H_4$ and $H_5$ be defined by \eqref{termeH1}-\eqref{termeH5}. The following inequality holds 
\begin{equation*}
\sup_{x\in \mathbb{R}}
\left\vert\mathbb{P}\Bigg(\frac{\sqrt{n}\left(
\hat \tau-\tau\right)}{\sqrt{\widehat{v}}}    \leq x    \Bigg)-
\Phi(x)\right\vert
\leq \frac{\log(n)}{\sqrt{n}}H_1(X,Y)+ \frac{H_2(X,Y)}{\sqrt{n}}
+ \frac{H_3(X,Y)}{\sqrt{n\log(n)}}+ \frac{H_4(X,Y) }{n}+\frac{H_5(X,Y)}{n^{3/2}}.
\end{equation*}

\end{Proposition}

Note that the terms $H_1(X,Y), H_2(X,Y),\ldots, H_5(X,Y)$ are finite as soon as the variables $X$ and $Y$ have a finite moment of order 12, and $\tilde X$, $\tilde Y$, $\varepsilon$, $\tilde X \varepsilon$, $\tilde X \tilde Y$ are not equal to zero almost surely. If the variables $X$ and $Y$ have a moment $q \in [6, 12)$, we can adapt the proof of Proposition \ref{propositionattenduesansj} to get a rate of order $n^{-(q-4)/(q+4)}\log(n)$, which gives a rate of order $n^{-1/5}\log(n)$ for $q=6$ and $n^{-1/3}\log(n)$ for $q=8$.

A natural question is: is the rate $n^{-1/2}\log(n)$  valid if the variables $X$ and $Y$ only have moments of order 6? One might think that this is the case, since a result of Bentkus-G\"otze \cite{BG} implies that 
\begin{equation*}
\sup_{x\in \mathbb{R}}\left\vert\mathbb{P}\left(\frac{\sum_{i=1}^n
(X_i-\mu_X)\varepsilon_i}
{
\sqrt{\sum_{i=1}^n
(X_i-\mu_X)^2\varepsilon_i^2}
}\leq x \right)-\Phi(x)\right\vert\leq \frac{c_{\rm{\scriptscriptstyle BG}}}{\sqrt{n}}\frac{\|\tilde X \varepsilon \|_3^3}{\|\tilde X \varepsilon \|_2^3}\, ,
\end{equation*}
where $c_{\rm{\scriptscriptstyle BG}}$ is a numerical constant.
However, the proof of Bentkus-G\"otze does not seem to adapt easily to our case. The main problem is that the term  $\widehat{v}$ in \eqref{defvhat} is calculated using the residuals $\hat \varepsilon_i$'s and not   the errors  $\varepsilon_i$'s (that are not observed). The main technical difficulty is therefore to replace the residuals  $\hat\varepsilon_i$'s in the expression of $\widehat{v}_j$ with the unobserved errors $\varepsilon_i$'s. 

\section{Theoretical guarantees for robust screening}
\label{sec3}
\setcounter{equation}{0}
Let  $n\in\mathbb{N}^*$ and let $(Y_i,X_i^{(1)},\ldots,X_i^{(p)})_{1\leq i\leq n}$ be $n$ independent copies of
$(Y,X^{(1)},\ldots,X^{(p)})$ where $Y$ is the output  and $X^{(1)},\ldots,X^{(p)}$ are the explanatory  random variables.

Recall that the quantities $\tau_j$, $\widehat{\tau}_j$,
the random variables $\varepsilon^{(j)}$ and the sets $\mathcal{M}^*$ and $\widehat{\mathcal{M}}$ are defined in the introduction (see (\ref{bj_tauj}), (\ref{tauchapoj_bchapoj}),  (\ref{model}), (\ref{M*}), (\ref{Mchapo})).
Moreover, for  $j=1,\ldots,p$, we define 
\begin{align*}
\mu_X^{(j)}=\E(X^{(j)}), \quad\tilde{X}^{(j)}=X^{(j)}-\mu_X^{(j)}, \quad \tilde{X}^{(j)}_i=X^{(j)}_i-\mu_X^{(j)}, 
\sigma_{X^{(j)}}^2=\Var\left(X^{(j)}\right), \quad \sigma_{\varepsilon^{(j)}}^2=\Var\left(\varepsilon^{(j)}\right).
\end{align*}
\subsection{Sure Screening Property} \label{sec3.1}
The following proposition ensures that, under moment assumptions, 
$$\lim_{\nti}\mathbb{P}(\mathcal{M}^*\subset \widehat{\mathcal{M}})=1.$$
\begin{Proposition}
\label{propositionSureScreeningProperty}
 Let $\widehat{\mathcal{M}}$ be defined by  \eqref{Mchapo}, with   $\gamma=\Phi^{-1}(1-q/2)$ where  $q$ is the chosen false positive rate. Let $K_1$ and $K_2$ be defined by \eqref{kappa1_kappa2} and $G_1$, $G_2$ and $G_3$ be defined by \eqref{termeG1}, \eqref{termeG2} and \eqref{termeG3}. 
 Then  
\begin{align*}
\mathbb{P}(\mathcal{M}^*\subset \widehat{\mathcal{M}})&\geq 
1-\sum_{j\in\mathcal{M}^*}{\left[
  \exp\left(-\frac{1}{4}\left(\frac{\sqrt{n}\sigma_{X^{(j)}}^2|\tau_j|}{2\sqrt{2}\big\| \tilde{X}^{(j)}\varepsilon^{(j)}\big\|_2}-\gamma\right)^2\,\right)\right]} - \frac{\log(n)}{\sqrt{n}}\left\vert\mathcal{M}^*\right\vert\max_{j\in\mathcal{M}^*}{\left(2K_1(X^{(j)},Y)\right)}\\
  & \qquad  - \frac{1}{\sqrt{n}}\left\vert\mathcal{M}^*\right\vert\max_{j\in\mathcal{M}^*}{\left(2K_2(X^{(j)},Y)+\frac{2\big\Vert \tilde{X}^{(j)}\big\Vert_4^2}{\sigma_{X^{(j)}}^2} +G_1(X^{(j)},Y)\right)} \\
  & \qquad  - \frac{1}{n}\left\vert\mathcal{M}^*\right\vert\max_{j\in\mathcal{M}^*}{\left(4G_2(X^{(j)},Y)+2\right)}- \frac{1}{n^{3/2}}\left\vert\mathcal{M}^*\right\vert\max_{j\in\mathcal{M}^*}{\left(G_3(X^{(j)},Y)\right)}.
\end{align*}
\end{Proposition}
In Proposition \ref{propositionSureScreeningProperty} we see that the rate of convergence is of order $ n^{-1/2}\log(n) \left\vert\mathcal{M}^*\right\vert$ if all the variables have a finite moment of order 6. Note that our proof combines Lemma \ref{prelim} and Lemma \ref{lemme2}, so it is not possible to get a better rate with this approach. However, if instead of Lemma \ref{prelim} we use the Fuk-Nagaev Inequality \eqref{FN} below (see \cite{FN}), and if we use the same inequality to get an upper bound for the left hand side of  \eqref{majoration_terme2_2_reste_prop} of Lemma \ref{lemme2}, then we could obtain a rate of order $n^{-(q-4)/4} \left\vert\mathcal{M}^*\right\vert$ if all the variables have a finite moment of order $q>6$. For instance, if all the variables have a finite moment of order 12, this gives a rate of order $n^{-2}\left\vert\mathcal{M}^*\right\vert$.

For the sake of completeness, let us recall Fuk-Nagaev's Inequality given in \cite{FN} (in the i.i.d. setting):
Let $Z_1,\ldots, Z_n$,  be $n$ i.i.d. random variables such that $\mathbb{E}(\vert Z_1\vert^q )<\infty$ for $q > 2$ and ${\mathbb E}(Z_1)=0$. Let  
 $\sigma_Z^2=\Var(Z)$.
Then, for any $x>0$,
 \begin{align}
 \label{FN}
 \mathbb{P}\Bigg[  \left\vert  \frac{1}{n}\sum_{i=1}^n Z_i       \right\vert >x         \Bigg] \leq 2\exp\bigg(-\frac{nx^2}{c_{\texttt{FN}}\sigma_Z^2}\bigg)+c(q)\frac{\mathbb{e}\mathbb{E}(\vert Z_1\vert^q)}{x^{q}n^{q-1}} \, ,
 \end{align}
 where $c_{\texttt{FN}}$ is a numerical constant and $c(q)$ is a constant depending on $q$.
 
 Let us continue the discussion about the rates. If we assume now that all the variables  $(X^{(j)})^4$'s   and $Y^2(X^{(j)})^2$'s have an exponential moment, then, applying Bernstein Inequality, we could get a rate of order $\exp(-\kappa n) \left\vert\mathcal{M}^*\right\vert$ for some $\kappa>0$. Regarding the cardinal $\left\vert\mathcal{M}^*\right\vert$, we could also be more precise: for instance, if $Y$ is bounded, and  if $m_1$ is the number of variables indexed in $\mathcal{M}^*$ having a moment of order $q>6$, and $m_2$ is the the number of variables indexed in $\mathcal{M}^*$ that are bounded (with $m_1+m_2= \left\vert\mathcal{M}^*\right\vert$), then we could get a rate of order 
 $n^{-(q-4)/4} m_1
 + \exp(-\kappa n) m_2$ for some $\kappa>0$.
 
\subsection{Control of the False Positive Rate}
\label{sec3.2}

In the following proposition, we give a non asymptotic control of the quantity 
$$
\left\vert \mathbb{E}\Bigg[\frac{|{\mathcal{M}^*}^c\cap \widehat{\mathcal{M}}|}{|{\mathcal{M}^*}^c|}\Bigg]-q \right\vert.
$$

\begin{Proposition}
\label{propositionFPR}
Let $\widehat{\mathcal{M}}$ be defined by  \eqref{Mchapo}, with $\gamma=\Phi^{-1}(1-q/2)$ where  $q$ is the chosen false positive rate.
Let $H_1$, $H_2$, $H_3$, $H_4$ and $H_5$ be defined by \eqref{termeH1}-\eqref{termeH5}. 
Then
\begin{multline*}
\left\vert \mathbb{E}\Bigg[\frac{|{\mathcal{M}^*}^c\cap \widehat{\mathcal{M}}|}{|{\mathcal{M}^*}^c|}\Bigg]-q \right\vert
\leq 
\frac{\log(n)}{\sqrt{n}}\max_{j\in{\mathcal{M}^*}^c}{\left(2H_1(X^{(j)},Y)\right)}+ \frac{1}{\sqrt{n}}\max_{j\in{\mathcal{M}^*}^c}{\left(2H_2(X^{(j)},Y)\right)} \\
+\frac{1}{\sqrt{n\log(n)}}\max_{j\in{\mathcal{M}^*}^c}{\left(2H_3(X^{(j)},Y)\right)}+ \frac{1}{n}\max_{j\in{\mathcal{M}^*}^c}{\left(2H_4(X^{(j)},Y) \right)}+\frac{1}{n^{3/2}}\max_{j\in{\mathcal{M}^*}^c}{\left(2H_5(X^{(j)},Y)\right)}.
\end{multline*}
\end{Proposition}

In Proposition \ref{propositionFPR} we see that the rate of convergence is of order $n^{-1/2}\log(n)$ if all the variables have a finite moment of order 12. Since the proof of Proposition \ref{propositionFPR} is based on Proposition \ref{propositionattenduesansj}, the comments after this proposition apply: if all the variables have a moment of order $q \in [6,12)$ then the rate in Proposition \ref{propositionFPR} is of order $n^{-(q-4)/(q+4)}\log(n)$, which give a rate of order $n^{-1/5}\log(n)$ for $q=6$ and $n^{-1/3}\log(n)$ for $q=8$.

\section{Simulations and real data}
\label{Simulations}
\setcounter{equation}{0}
\subsection{Simulations}

\subsubsection{Description of the explanatory variables}
\label{descriptionvarexplicative} 

We have simulated  $p=2000$ random explanatory variables according to 4 groups of variables. 
\begin{itemize}
\item 
The random variables $X^{(1)}, \ldots, X^{(10)}$ are used to generate the output $Y$.

\begin{itemize}
\item The variables
$X^{(1)},X^{(2)},X^{(3)},X^{(4)},X^{(5)}$ are standard dependent Gaussian random variables. Let  $\rho=0.57$.  The covariance matrix of ($X^{(1)},X^{(2)},X^{(3)},X^{(4)},X^{(5)}$) is given by
$$\Sigma=\begin{pmatrix}
1 & \rho & \rho^2 & \rho^3 & \rho^4 & \\
\rho & 1 &\rho &\rho^2 & \rho^3 \\
\rho^2 &\rho &1 &\rho &\rho^2& \\
\rho^3 &\rho^2 &\rho & 1 &\rho &\\ 
\rho^4 &\rho^3 &\rho^2& \rho& 1  \\
\end{pmatrix}.$$

\item $X^{(6)}\sim \mathcal{B}(0.35)$ (Bernoulli) and is  independent of the other explanatory variables.

\item $X^{(7)}\sim \chi^2(2)$ ($\chi^2$-square) and is independent of the other explanatory variables. 

\item  $X^{(8)}=X^{(1)}\times Z_1$ where $Z_1\sim \mathcal{P}(2)$ (Poisson) is independent of the other explanatory variables. 

\item  $X^{(9)}=X^{(2)}\times Z_2$ where $Z_2 \sim \mathcal{N}(1,1)$ is independent of the other explanatory variables. 

\item  $X^{(10)}\sim \mathcal{S}t(14)$ (Student's distribution with 14 degrees of freedom) and is independent of the other explanatory variables. 
\end{itemize}

\item 
The variables $X^{(11)}$, $X^{(12)}$, $X^{(13)}$, $X^{(14)}$ are noisy versions of $X^{(1)}$, $X^{(3)}$, $X^{(4)}$, $X^{(5)}$. They are correlated with $X^{(1)}$, $X^{(3)}$, $X^{(4)}$ and $X^{(5)}$, and with the response  $Y$,  but they are not used to generate $Y$.

Consider $4$ random variables $Z_3$, $Z_4$, $Z_5$, $Z_6$ independent of $X^{(1)},\ldots,X^{(10)}$ and such that
$$Z_3\sim \mathcal{N}(0,0.4), \qquad Z_4\sim \mathcal{N}(0,0.02), \qquad 
Z_5\sim \mathcal{N}(0,0.35), \qquad Z_6\sim \mathcal{N}(0,0.55).$$
We then define the variables  $X^{(11)}$,\ldots,$X^{(14)}$ as follows:\\
$X^{(11)}=X^{(1)}+Z_3$;
$X^{(12)}=X^{(3)}+Z_4$;
$X^{(13)}=X^{(4)}+Z_5$;
$X^{(14)}=X^{(5)}+Z_6$.

\item The variables $X^{(15)}$, $X^{(16)}$, $X^{(17)}$, $X^{(18)}$ are correlated with $ X^{(11)}$, $X^{(12)}$, $X^{(13)}$, $X^{(14)}$ but are independent of the output $Y$.
More precisely:
$X^{(15)}=Z_3$;
$X^{(16)}=Z_4$;
$X^{(17)}=Z_5$;
$X^{(18)}=Z_6$.

\item The last random variables $X^{(19)}, X^{(20)}, \ldots, X^{(2000)}$ are independent of $(Y, X^{(1)}, \ldots, X^{(18)})$. These variables form a stationary Gaussian sequence with mean 0, variance 1 and $1982\times 1982$  covariance matrix similar to the matrix $\Sigma$ above with $\rho=0.7$. 
\end{itemize}

\subsubsection{Outputs}

We chose to illustrate our simulation study on two models:
one with continuous output, the other with a discrete output, both models  not belonging to the class of GLM models.
Furthermore, for the model with continuous output, we consider an heteroscedastic noise.
In both models, the explanatory variable $X^{(10)}$ only admits  a moment of order 13.

\paragraph{Continous output}
Let $\theta^*$ be defined by
\begin{align}\theta^*= (-1,  3.619350, -3.274923,  2.963273, -2.681280,  2, 4,6,  3,  2,  4)^T\label{thetastar}
\end{align}
and $Y$ such that
$$Y= \left\vert \theta_0^* + \theta_1^* X^{(1)}+ \cdots + \theta^*_{10} X^{(10)}\right\vert^{0.8}+X^{(5)}\varepsilon/3,$$
where
$\theta^*=(\theta_0^*,\ldots,\theta^*_{10})$ is given by (\ref{thetastar}) and $\varepsilon$ is a centered random variable independent of the $X^{(j)}$'s with $\varepsilon\sim\mathcal{X}^2(3)-3$.


\paragraph{Discrete output}
The conditional distribution of $Y$ given $(X^{(1)},\ldots,X^{(10)})$, is a Poisson distribution with parameter 
$$\left\vert \theta_0^* + \theta_1^* X^{(1)}+ \cdots + \theta^*_{10} X^{(10)}\right\vert$$ where $\theta^*=(\theta_0^*,\ldots,\theta^*_{10})$ is given by (\ref{thetastar}).


\subsubsection{Simulations results}

We consider a sample size $n$ equal to $n=500$, $n=1000$, $n=2000$, $n=2500$. We perform a Monte Carlo step over $N=500$ samples to evaluate the frequency of selection of each explanatory random variable, the True Positive Rate, the False Positive Rate and the mean of $\vert\widehat{\mathcal{M}}\vert$. A True Positive is an index between 1 and 14 selected by the procedure. A False Positive is an index between 15 and 2000 selected by the procedure. The results for the continuous output with $q=0.10$ are given in Table 1. Table 2 presents the results for the discrete output with $q=0.15$. 

\begin{table}

\centering
\begin{tabular}[t]{|l|r|r|r|r|}
\hline
Explanatory Variable & $n=500$ & $n=1000$ & $n=2000$ & $n=2500$\\
\hline
$X^{(1)}$ & 1.000&  1.000& 1.000& 1.000\\
\hline
$X^{(2)}$ &0.996&1.000& 1.000 &1.000\\
\hline
$X^{(3)}$ &0.432&0.720& 0.912&0.952\\
\hline
$X^{(4)}$ &0.638&0.838& 0.988&0.992\\
\hline
$X^{(5)}$ &0.730&0.956& 1.000&1.000\\
\hline
$X^{(6)}$ &0.388&0.570& 0.820&0.906\\
\hline
$X^{(7)}$ &1.000&1.000& 1.000&1.000\\
\hline
$X^{(8)}$ &0.980&1.000& 1.000&1.000\\
\hline
$X^{(9)}$ &0.990&1.000& 1.000&1.000\\
\hline
$X^{(10)}$ &0.994&1.000& 1.000&1.000\\
\hline
$X^{(11)}$ &0.998&1.000& 1.000&1.000\\
\hline
$X^{(12)}$ &0.282&0.498& 0.716&0.748\\
\hline
$X^{(13)}$ &0.558&0.800& 0.978& 0.976\\
\hline
$X^{(14)}$ &0.650&0.884& 0.994&0.998\\
\hline
$X^{(15)}$ &0.118&0.740& 0.114&0.104\\
\hline
$X^{(16)}$ &0.108&0.086& 0.118&0.112\\
\hline
$X^{(17)}$ &0.120&0.092& 0.092&0.122\\
\hline
$X^{(18)}$&0.116&0.120& 0.116&0.108\\
\hline
True Positive Rate & 0.760&0.876&0.958& 0.969\\
\hline
False Positive Rate & 0.104&0.102&0.101&0.100\\
\hline
mean $|\widehat{\mathcal{M}}|$ & 217.124&214.422&214.040&212.682\\
\hline
\end{tabular}
\caption{$N=500$, $q=0.10$, continuous output}
\end{table}

All variables correlated with the output $Y$ are correctly selected by the procedure and eight of them are systematically selected when $n=2500$. We note that the frequency of selection of the variable $X^{(12)}$, which is weakly correlated with $Y$, is further from 1 than that of the other explanatory variables. It reaches, when $n=2500$, $74.8 \%$  for the continuous output and $79.8 \%$ for the discrete output. We also note that the selection rate of each variable correlated with $Y$ increases as $n$ increases, as predicted by the theoretical results.
Moreover the True Positive Rate is close to 1 when the sample size is large enough. For the continuous output, it increases from $76 \%$ when $n=500$ to $96.9\%$ when $n=2500$ and for the discrete  output, it increases from $79.4 \%$ when $n=500$ to $97.8\%$ when $n=2500$. As expected, the False Positive Rate is close to $10\%$ for the continuous output (with $q=0.10$) and close to $15\%$ for the discrete output (with $q=0.15$). 
Finally, the mean cardinality of $\widehat{\mathcal{M}}$ is close to its expected value (212.6 when $q=0.10$ and 311.9 when $q=0.15$).

As indicated in the introduction, our procedure does not  assume that the random variables indexed in ${\mathcal M}^*$ are globally independent of the variables indexed in  ${\mathcal{M}^*}^c$. Here the random variables $X^{(15)}$, $X^{(16)}$, $X^{(17)}$, $X^{(18)}$ belong to
${\mathcal{M}^*}^c$ but are correlated to some variables in ${\mathcal{M}^*}$. For these four variables, as predicted by the theoretical results, the selection rate is close to the false positive rate $q$.

\begin{table}
\centering
\begin{tabular}[t]{|l|r|r|r|r|}
\hline
Explanatory Variable & $n=500$ & $n=1000$ & $n=2000$ & $n=2500$ \\
\hline
$X^{(1)}$ & 1.000& 1.000 &1.000 & 1.000\\
\hline
$X^{(2)}$ & 0.996 & 1.000 & 1.000 & 1.000\\
\hline
$X^{(3)}$ & 0.554 & 0.696& 0.934 & 0.982\\
\hline
$X^{(4)}$ & 0.654& 0.872& 0.988& 0.996\\
\hline
$X^{(5)}$ & 0.814& 0.952& 0.996 & 1.000\\
\hline
$X^{(6)}$ & 0.476& 0.600 & 0.864 & 0.926\\
\hline
$X^{(7)}$ & 1.000& 1.000& 1.000 & 1.000\\
\hline
$X^{(8)}$ & 0.986& 1.000& 1.000& 1.000\\
\hline
$X^{(9)}$ & 0.992& 1.000& 1.000& 1.000\\
\hline
$X^{(10)}$ & 0.996& 1.000& 1.000& 1.000\\
\hline
$X^{(11)}$ & 0.996& 1.000& 1.000& 1.000\\
\hline
$X^{(12)}$ & 0.334 & 0.520& 0.756& 0.798\\
\hline
$X^{(13)}$ & 0.638& 0.814& 0.982& 0.990\\
\hline
$X^{(14)}$ & 0.686& 0.912& 0.986 & 0.996\\
\hline
$X^{(15)}$ & 0.174& 0.144 & 0.160 & 0.162\\
\hline
$X^{(16)}$ & 0.144& 0.154 & 0.122 & 0.138\\
\hline
$X^{(17)}$ & 0.138& 0.130& 0.160&  0.158\\
\hline
$X^{(18)}$ & 0.160& 0.172& 0.126 & 0.152\\
\hline
True Positive Rate & 0.794& 0.883& 0.965& 0.978\\
\hline
False Positive Rate & 0.154 & 0.152& 0.151& 0.150\\
\hline
mean $|\widehat{\mathcal{M}}|$ & 316.936& 313.432 & 312.310 & 312.782\\
\hline
\end{tabular}
\caption{$N=500$, $q=0.15$, discrete output}
\end{table}

\subsection{Real data set}
\label{dataset}

In this section, we will apply our method to a genomic dataset, the single-cell RNA-seq from Shalek \textit{et al.} \cite{Shalek2014}. These data have already been examined using statistical methods in two articles (\cite{Cai} and \cite{Hou2024}), and our goal is to see if our selection rule is compatible with the results obtained in these articles, but also to highlight some differences. 

Let us begin by reviewing the description of these data, as in \cite{Hou2024}: the data comprise gene expression profiles of 27723 genes from 1861 primary dendritic cells derived from mouse bone marrow under several experimental conditions.

As in \cite{Hou2024}, we focus on a subset of data consisting of 96 cells stimulated by the PIC pathogen component (viral-like double-stranded RNA) and 96 unstimulated control cells, whose gene expression was measured six hours after stimulation. The goal is therefore to select the gene expression variables that are correlated with the variable $Y$, which is a binary variable coded 0 for ``unstimulated'' and 1 for ``stimulated''. 

In principle, we do not need to perform any pre-processing to apply our method (the calculation times for 27723 genes being quite reasonable). However, in order to start from a basis comparable to \cite{Cai} and \cite{Hou2024}, we will still filter out genes that are not expressed in more than 80\% of cells: this gives us a total of 6966 genes. Note that, unlike \cite{Cai} and \cite{Hou2024}, we do not perform the other pre-processing step of eliminating the 90\% of genes with the lowest variance among the 6966, as this second pre-processing step seems somewhat arbitrary compared to our selection rule. 

\begin{figure}[h!]
    \centering
\includegraphics[height=7cm, width=13cm]{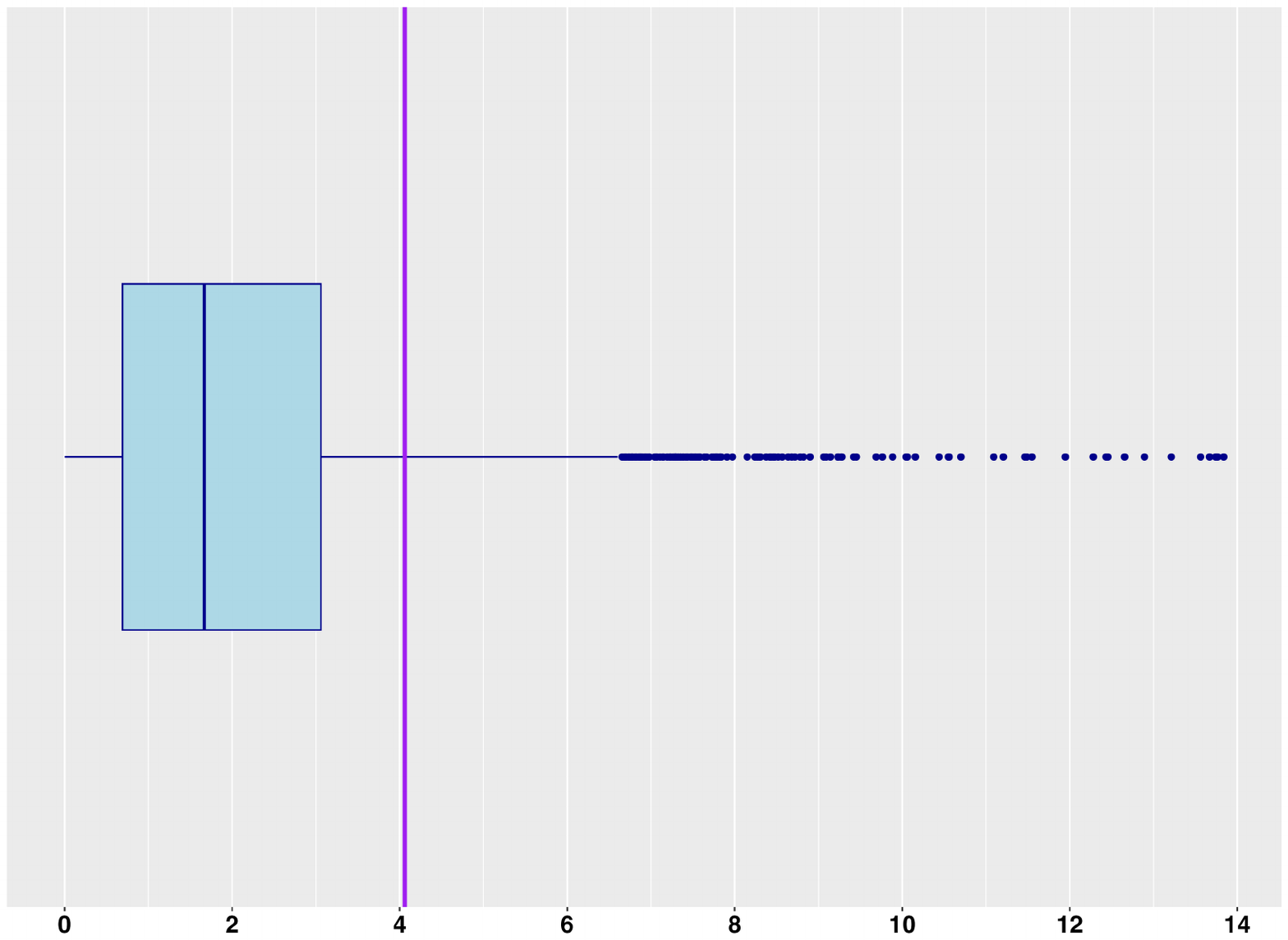}
    \caption{Boxplot of the absolute values of the $\hat \gamma_j$'s. The vertical line represents
    the threshlold 4.056}
    \label{fig:boxplot}
    \end{figure}

For the $j$-th gene among the 6966,  we compute its {\it score} $\hat \gamma_j$, that is the observed value of 
$
\sqrt{n} {\widehat{v}_j}^{-1/2}\widehat{\tau}_j
$
According to our selection rule \eqref{Mchapo}, we will keep all the genes whose score are greater (in absolute value) than a given threshold. To have an idea of the distribution of the scores, the boxplot of the $| \hat \gamma_j|$ is given in Figure \ref{fig:boxplot}.
From this boxplot, we can deduce that a large number of genes appear to be expressed differently depending on whether the cells are stimulated or not. For example, the median is 1.67 and the 90\% quantile is 4.4. This means that, for this data set, we are a priori very far from the situation known as sparsity.

We therefore propose to apply our method with a threshold of 4.056, which corresponds to a false positive rate of $q=5\times 10^{-5}$. Doing so, it is likely that we will not select all the genes whose expression is associated with the PIC cell stimulation, but only a subset of those genes for which the association appears to be particularly strong. 

With this threshold of 4.056, we select 918 genes from the initial 6966, representing a selection rate of $13.2\%$. We can now list the genes in our selection that are common to those found in \cite{Hou2024}, based on the adjustment of multivariate models. In \cite{Hou2024}, 9 genes were identified  as belonging to a relevant multivariate model (the largest models containing at most 6 of these genes). Among these 9 genes, 6 are also selected by our method. They are listed in Table 3 below, where we also give their $\hat \gamma$ scores. In particular, we find the RSAD2 gene (with a very high score of 8.64), which had already been identified in \cite{Cai}. 
In Table 4, we present the 3 genes that are not selected by our method, among the 9 identified in \cite{Hou2024}.

\begin{table}[h!]
\begin{subtable}{0.5\textwidth}
\centering
\begin{tabular}{|c|c|}\hline
Gene & $\hat \gamma$\\
\hline
RSAD2 & 8.64
\\
\hline
ACTB & 9.76
\\
\hline
BC044745 & -4.24
\\
\hline
HMGN2 & -6.87
\\
\hline
IFI47 & 8.31
\\
\hline
ZFP488 & 4.93
\\
\hline
\end{tabular}
 \caption{Table 3}
\end{subtable}
 \label{Table_3}
\begin{subtable}{0.5\textwidth}
\centering
\begin{tabular}{|c|c|}
\hline
Gene & $\hat \gamma$\\
\hline
IFT80 & 2.08\\
\hline
IFIT1 & 2\\
\hline
AK217941 & -2.88\\
\hline
\end{tabular}
 \caption{Table 4}
\end{subtable}
    \label{Table_4}
\end{table}

 We see that the scores of the 3 genes IFT80, IFIT1, AK217941 are well below our threshold. The main difference with \cite{Hou2024} concerns the gene AK217941, which appears in several multivariate models proposed in \cite{Hou2024}. In \cite{Hou2024}, the authors write ``AK217941, though not studied in the literature, deserves further attention''. Our analysis shows that the association between the expression of the gene AK217941 and PIC cell stimulation is not particularly significant compared to that of many other genes. 




We can also extract from our analysis the two selected genes with the highest $\hat \gamma$ scores. These are genes PRDX6-RS1-PS and PIGR, with respective scores of -13.84 and -13.77. The PRDX6 gene in mice plays an important role in  defence against oxidative stress (see \cite{Simeone2005}). The PIGR gene is known to respond to pathogenic stimulation. The most prominent modulators of PIGR regulation consist of TLR4 and TLR3, which recognize bacterial lipopolysaccharide and viral dsRNA respectively (see \cite{Kaetzel2005}). Given this information, it therefore seems perfectly natural to find these two genes among those ultimately selected by our procedure. On the contrary, it is surprising that these two genes – whose main effect seems indisputable – do not appear in the multivariate models proposed in \cite{Hou2024} (which can be explained by the second pre-processing step performed in this article). 


        

\section{Proofs} \label{Proofs}
\setcounter{equation}{0} 

\subsection{Known results and preliminary lemmas} \label{5.1}

In this section, we first recall some known results used all along the proofs. 

According to Esseen \cite{Esseen1942}, we have the following theorem.
\begin{Theorem}\textbf{[Berry-Esseen inequality]}
 Consider $Z_1,\ldots, Z_n$, $n$ independent and identically distributed centered random variables such that $\mathbb{E}(\vert Z_1 \vert ^3)< \infty$ and let $\overline{\Z}=n^{-1}\sum_{i=1}^n{Z_i}$.
 Set $\rho_Z=\mathbb{E}(\vert Z_1 \vert ^3)$ and $\sigma^2_Z=\Var(Z_1).$
 There exists a numerical constant $c_{\rm{\scriptscriptstyle BE}}$ such that for any $x\in \mathbb{R}$ and for any positive integer $n$,
 \begin{align}
 \label{BE}
\left\vert \mathbb{P}\Bigg[ \frac{ \sqrt{n}\,\overline{\Z}}{\sigma_Z} \leq x\Bigg]-\Phi(x)\right\vert \leq c_{\rm{\scriptscriptstyle BE}}\frac{   \rho_Z    }{ \sigma_Z^3\sqrt{n}}.
 \end{align}
 \end{Theorem}
The following lemma is due to Cuny \textit{et al.} \cite{CDMP2022}.
\begin{Lemma} 
\label{CRAS}
Let $(T_n)_{n\in \mathbb{N}}$ and $(R_n)_{n\in \mathbb{N}}$ be two sequences of random variables. Assume that there exist three sequences of positives numbers $(a_n)_{n\in \mathbb{N}}$, $(b_n)_{n\in \mathbb{N}}$ and $(c_n)_{n\in \mathbb{N}}$, going to infinity as $n\rightarrow\infty$, and a positive constant $s$ such that, for any positive integer $n$,
$$\sup_{t\in \mathbb{R}}\left\vert\mathbb{P}( T_n \leq t\sqrt{n})-\Phi(t/s)\right\vert \leq \frac{1}{a_n}, \mbox{ and } \mathbb{P}(\vert R_n\vert \geq
s\sqrt{n} / b_n)\leq \frac{1}{c_n}.$$
Then,
$$\sup_{t\in \mathbb{R}}\left\vert \mathbb{P}(T_n+R_n\leq t\sqrt{n})-\Phi(t/s)\right\vert\leq \frac{1}{a_n}+\frac{1}{\sqrt{2\pi}b_n}+\frac{1}{c_n}.$$
\end{Lemma}

According to Rosenthal \cite{Rosenthal1970}, we have the following theorem. The constants can be found in Pinelis \cite{Pinelis1994}.
\begin{Theorem}\textbf{[Rosenthal Inequality]}
 Consider $Z_1,\ldots, Z_n$, $n$ independent and identically distributed centered random variables and let $q\geq 2$. Then, there exist two numerical constants $C_1(q)=c_{1,\rm{\scriptscriptstyle R}}\sqrt{q}$ and  $C_2(q)=c_{2,\rm{\scriptscriptstyle R}} q$  such that
 \begin{align}
 \label{rosenthal}
 \left \| \frac{1}{n}\sum_{i=1}^n Z_i \right \|_q \leq C_1(q)\sqrt{
 \frac{\Var(Z_1)}{n}}+C_2(q)\frac{1}{n^{(q-1)/q}}\left\| Z_1\right\|_q.
\end{align}
\end{Theorem}
According to von Bahr and Esseen \cite{BE1965}, we have the following theorem. The constant $2^{2-q}$  can be found in Pinelis \cite{Pinelis2015}.
\begin{Theorem}\textbf{[von Bahr - Esseen Inequality]}
Consider $Z_1,\ldots, Z_n$, $n$ independent and identically distributed centered random variables and let $q\in (1,2)$. Then
\begin{align}
 \label{Von_Bahr_Esseen}
\left\Vert \frac{1}{n}\sum_{i=1}^n Z_i \right \Vert_q^q \leq \frac{2^{2-q}}{n^{q-1}}  \Vert Z_1\Vert_q^q \, .
\end{align}
\end{Theorem}
According to Pollak \cite{Pollak1956} we have the following lemma.
\begin{Lemma}\textbf{[Pollak Inequality]}
For any $x\in\R_+$, 
\begin{align}
1-\Phi(x)\leq \frac{1}{2}\exp\left(-\frac{x^2}{2}\right)
\label{Pollak}
\end{align}
\end{Lemma}
We also use several times the following elementary lemmas. The proofs are left to the reader except for Lemma \ref{lemmemajoration_var_X} for which we give a short proof.
\begin{Lemma}
\label{lemmeproduit}
Consider $Z_1, \ldots Z_k$, $k$ positive random variables and let $a_2,\ldots,a_k$  and $c$ be some positive constants.
Then
$$\mathbb{P} \Big(  Z_1 Z_2\ldots Z_k >c       \Big)\leq \mathbb{P}(Z_1a_2a_3\ldots a_k>c)+\sum_{j=2}^k\mathbb{P}(Z_j>a_j).$$
\end{Lemma}

 
\begin{Lemma}
\label{tauchapo}
Let $(X_1,Y_1),\ldots,(X_n,Y_n)$ be $n$ independent copies of a random couple $(X,Y)$ and $\widehat{\tau}$ be defined by (\ref{deftauhat}). Then
\begin{align}
|\widehat{\tau}|\leq \frac{ \sqrt{\frac{1}{n}\sum_{i=1}^n(Y_i-\overline{\Y})^2       }}{\sqrt{\frac{1}{n}\sum_{i=1}^n(X_i-\overline{\X})^2}}.
\label{majorationtauchapeau}
\end{align}
\end{Lemma}

\begin{Lemma}
\label{lemmemajoration_var_Y}
Let $Y_1,\ldots,Y_n$ be $n$ independent copies of a random variable $Y$ and consider the notations (\ref{notationsY}). Then  
\begin{eqnarray}
\label{majoration_var_Y}
\mathbb{P}\Bigg[   \frac{1}{n}\sum_{i=1}^n(Y_i-\overline{\Y})^2 \geq 2\sigma_Y^2   \Bigg]\leq \frac{ \Vert\tilde{Y}\Vert_4^2}{\sigma_Y^2\sqrt{n}}.
\end{eqnarray}
\end{Lemma}
\begin{Lemma}
\label{lemmemajoration_var_X}
Let $X_1,\ldots,X_n$ be $n$ independent copies of a random variable $X$ and consider the notations (\ref{notationsX}). Then  
\begin{eqnarray}
\label{majoration_var_X}
\mathbb{P}\Bigg[   \frac{1}{n}\sum_{i=1}^n(X_i-
\overline{\X})^2\leq\frac{\sigma_X^2}{2}\Bigg]
\leq \frac{2}{\sqrt{n}\sigma_X^2} \Vert\tilde{X}\Vert_4^2
+\frac{2}{n}. 
\end{eqnarray}
\end{Lemma}

\subsubsection*{Proof of Lemma \ref{lemmemajoration_var_X}.}
\begin{align*}
\mathbb{P}\Bigg[   \frac{1}{n}\sum_{i=1}^n(X_i-
\overline{\X})^2\leq\frac{\sigma_X^2}{2}\Bigg]
&\leq \mathbb{P}\Bigg[  \left\vert \frac{1}{n}\sum_{i=1}^n(X_i-
\overline{\X})^2-\sigma_X^2\right\vert\geq\frac{\sigma_X^2}{2}\Bigg]
\leq 
\frac{2}{\sigma_X^2}\left\Vert \frac{1}{n}\sum_{i=1}^n{\tilde{X}_i^2}-(\overline{\X}-\mu_X)^2-\sigma_X^2\right\Vert_1\nonumber\\
&\leq \frac{2}{\sigma_X^2}\left\Vert \frac{1}{n}\sum_{i}^n\tilde{X}_i^2-\sigma_X^2\right\Vert_1+\frac{2}{\sigma_X^2}\left\Vert (\overline{\X}-\mu_X)^2\right\Vert_1
\leq \frac{2}{\sqrt{n}\sigma_X^2} \Vert\tilde{X}\Vert_4^2
+\frac{2}{n}. 
\end{align*}

\subsection{Definition of the Lyapunov ratios involved in Lemma \ref{lemme2} and Proposition \ref{propositionattenduesansj}}
\label{definitionG_ietH_i}
The Lyapunov ratios $G_1$, $G_2$ and $G_3$ involved in Lemma \ref{lemme2} are defined by 
\begin{align}
G_1(X,Y)&=\frac{64\big\Vert
\tilde{Y}^2 \tilde{X}\big\Vert_1\sigma_X}
{\big\Vert \tilde{X}\varepsilon\big\Vert_2^2}
+\frac{1120 \big\Vert \tilde{X}\big\Vert_3^3\sigma_Y^2}{\big\Vert \tilde{X}\varepsilon\big\Vert_2^2\sigma_X}
+\frac{64\sqrt{2}\big\Vert \tilde{X}\varepsilon\big\Vert_3^3}{\big\Vert \tilde{X}\varepsilon\big\Vert_2^3}+\frac{288\,\sigma_Y\big\Vert\tilde{X}^2 \tilde{Y}\big\Vert_1}{
\big\Vert \tilde{X}\varepsilon\big\Vert_2^2}
+
\frac{432 \vert\tau \vert\big\Vert \tilde{X}\big\Vert_3^3\sigma_Y}{\big\Vert \tilde{X}\varepsilon\big\Vert_2^2}\nonumber\\
&+
\frac{1152\big\Vert \tilde{X}^3\tilde{Y}\big\Vert_1}{\sigma_X^2\big\Vert \tilde{X}\varepsilon\big\Vert_2}+
\frac{1152 \big\Vert \tilde{X}\big\Vert_4^4\sigma_Y}{\sigma_X^3\big\Vert \tilde{X}\varepsilon \big\Vert_2} 
+\frac{576\vert\tau\vert \big\Vert \tilde{X}\big\Vert_4^4}{\sigma_X^2\big\Vert \tilde{X}\epsilon\big\Vert_2} 
+
\frac{576 \big\Vert \tilde{X}^2\tilde{Y}\big\Vert_1\sigma_Y}{\big\Vert \tilde{X}\varepsilon\big\Vert_2^2\sigma_X}+
\frac{2\big\Vert\tilde{X}^2 \tilde{Y}\big\Vert_2}{\big\Vert \tilde{X}^2\tilde{Y}\big\Vert_1} \nonumber\\
&+
\frac{1}{2\big\Vert \tilde{Y}^2 \tilde{X}\big\Vert_1}\left(2\big\Vert\tilde{Y}^2\tilde{X}\big\Vert_2+\sigma_X\big\Vert \tilde{Y}\big\Vert_4^2+
4\sigma_Y
\big\Vert\tilde{X}\tilde{Y}\big\Vert_2\right) 
+ \frac{\sigma_Y\big\Vert \tilde{X}\big\Vert_4^2}{\big\Vert \tilde{X}^2\tilde{Y}\big\Vert_1}+ \frac{13\sigma_X \big\Vert \tilde{X}\big\Vert_4^2}{4\big\Vert \tilde{X}\big\Vert_3^3}\nonumber\\
& + \frac{9\big\Vert \tilde{X}\big\Vert_6^3}{2\big\Vert \tilde{X}\big\Vert_3^3} 
+ \frac{6\big\Vert \tilde{Y}\big\Vert_4^2}{\sigma_Y^2}
+ \frac{16\big\Vert \tilde{X}\big\Vert_4^2}{\sigma_X^2}+ \frac{2^{7/2}\big\Vert\tilde{X}^3\tilde{Y}\big\Vert_{3/2}^{3/2}}{\big\Vert \tilde{X}^3 \tilde{Y}\big\Vert_1^{3/2}}+ \frac{\sigma_Y\big\Vert \tilde{X}\big\Vert_6^3}{2\big\Vert \tilde{X}^3 \tilde{Y}\big\Vert_1}
+\frac{8\big\Vert \tilde{X}\big\Vert_6^6}{\big\Vert \tilde{X}\big\Vert_4^6}
+ \frac{\sigma_X\big\Vert \tilde{X}\big\Vert_6^3}{2\big\Vert \tilde{X}\big\Vert_4^4}, \label{termeG1}
\end{align}
\begin{align}
G_2(X,Y) &= 4+\frac{288 \big\Vert \tilde{X}^3 \tilde{Y}\big\Vert_1
\sigma_{\varepsilon}}{
\sigma_X\big\Vert \tilde{X}\varepsilon\big\Vert_2^2
}+\frac{288 \,\sigma_Y \sigma_{\varepsilon}\big\Vert \tilde{X}\big\Vert_4^4}{\big\Vert \tilde{X}\varepsilon\big\Vert_2^2\sigma_X^2}
+\frac{144  \vert\tau\vert\sigma_{\varepsilon}\big\Vert \tilde{X}\big\Vert_4^4 }{\sigma_X\big\Vert \tilde{X}\varepsilon\big\Vert_2^2}
+\frac{c_{1,\rm{\scriptscriptstyle R}}^2\sigma_X^3}{\big\Vert \tilde{X}\big\Vert_3^3}+\frac{2c_{1,\rm{\scriptscriptstyle R}}^2\sigma_X\sigma_Y^2}{\big\Vert \tilde{Y}^2 \tilde{X}\big\Vert_1}, \label{termeG2}
\end{align}
\begin{align}
G_3(X,Y)&=16c_{2,\rm{\scriptscriptstyle R}}^2\left[\frac{\sigma_X\big\Vert \tilde{X}\big\Vert_4^2}{\big\Vert \tilde{X}\big\Vert_3^3 } +
\frac{2\sigma_X\big\Vert \tilde{Y}\big\Vert_4^2}{\big\Vert \tilde{Y}^2 \tilde{X}\big\Vert_1}\right] \label{termeG3}.
\end{align}
Let $K$ be defined by 
\begin{multline}
    K= \frac{\big\Vert \tilde{X}\varepsilon\big\Vert_2}{4}\min \left(\frac{\big\Vert \tilde{X}\varepsilon\big\Vert_2}{36\big\Vert \tilde{X}^2\tilde{Y}\big\Vert_1\sigma_Y}, \frac{\big\Vert \tilde{X}\varepsilon\big\Vert_2\sigma_X}{72 \big\Vert \tilde{X}^2\tilde{Y}\big\Vert_1\sigma_Y}
,\frac{\big\Vert \tilde{X}\varepsilon\big\Vert_2\sigma_X}{72 \big\Vert \tilde{X}\big\Vert_3^3\sigma_Y^2}
,\frac{\big\Vert \tilde{X}\varepsilon\big\Vert_2}{36  \vert\tau\vert\big\Vert \tilde{X}\big\Vert_3^3\sigma_Y}, 
\right. \\
\left.\frac{\sigma_X^2}{288
    \big\Vert \tilde{X}^3\tilde{Y}\big\Vert_1}, 
\frac{\sigma_X^3}{288\big\Vert \tilde{X}\big\Vert_4^4
\sigma_Y},
\frac{\sigma_X^2}{72\vert\tau\vert\big\Vert \tilde{X}\big\Vert_4^4}, 
\frac{\big\Vert \tilde{X}\varepsilon\big\Vert_2}{\big\Vert \tilde{X}\varepsilon\big\Vert_4^2},\frac{
\big\Vert \tilde{X}\varepsilon\big\Vert_2}{8 \big\Vert
\tilde{Y}^2 \tilde{X}\big\Vert_1\sigma_X}\right).
\label{K_final}
\end{multline}
The Lyapunov ratios $H_1,\ldots,H_5$ involved in Proposition \ref{propositionattenduesansj} are defined by
\begin{align}
H_1(X,Y)& =\frac{3\sigma_{\varepsilon}\sigma_X}{
\sqrt{2\pi}\big\Vert\tilde{X}\varepsilon \big\Vert_2}+\frac{1}{\sqrt{2\pi}K}, \label{termeH1}
\end{align}
\begin{align}
H_2(X,Y)&=19+\frac{1}{\big\Vert \tilde{Y}^2 \tilde{X}\big\Vert_1}\left(2\big\Vert\tilde{Y}^2\tilde{X}\big\Vert_2+\sigma_X\big\Vert \tilde{Y}\big\Vert_4^2+
4\sigma_Y
\big\Vert\tilde{X}\tilde{Y}\big\Vert_2\right) \nonumber \\
&+\frac{9c_{\rm{\scriptscriptstyle BE}}
\big\Vert\tilde{X} \varepsilon\big\Vert_3^3}{\big\Vert \tilde{X}\varepsilon\big\Vert_2^3
}+\frac{6c_{\rm{\scriptscriptstyle BE}}\big\Vert \varepsilon\big\Vert_3^3
  }{\sigma_{\varepsilon}^3} 
+\frac{16c_{\rm{\scriptscriptstyle BE}}\big\Vert \tilde{X}\big\Vert_3^3}{\sigma_X^{3}}
+\frac{2c_{\rm{\scriptscriptstyle BE}}
\big\Vert \tilde{X}\varepsilon\big\Vert_6^6}{\left[\Var\left(\tilde{X}^2\varepsilon^2\right)\right]^{3/2}}+\frac{6c_{\rm{\scriptscriptstyle BE}}\big\Vert \tilde{Y}\big\Vert_3^3}{\sigma_Y^3} \nonumber\\
&\left.+
4\frac{\big\Vert\tilde{X}^2 \tilde{Y}\big\Vert_2}{\big\Vert \tilde{X}^2\tilde{Y}\big\Vert_1}
+ 2\frac{\sigma_Y\big\Vert \tilde{X}\big\Vert_4^2}{\big\Vert \tilde{X}^2\tilde{Y}\big\Vert_1}
+ \frac{13}{2}\frac{\sigma_X \big\Vert \tilde{X}\big\Vert_4^2}{\big\Vert \tilde{X}\big\Vert_3^3}
+ 9\frac{\big\Vert \tilde{X}\big\Vert_6^3}{\big\Vert \tilde{X}\big\Vert_3^3} 
+ 12\frac{\big\Vert \tilde{Y}\big\Vert_4^2}{\sigma_Y^2}
+ 32\frac{\big\Vert \tilde{X}\big\Vert_4^2}{\sigma_X^2}
\right. \nonumber\\
&+ 2^{11/2}\frac{\big\Vert\tilde{X}^3\tilde{Y}\big\Vert_{3/2}^{3/2}}{\big\Vert \tilde{X}^3 \tilde{Y}\big\Vert_1^{3/2}}
+ \frac{3}{4}\frac{\sigma_Y\big\Vert \tilde{X}\big\Vert_6^3}{\big\Vert \tilde{X}^3 \tilde{Y}\big\Vert_1} +16\frac{\big\Vert \tilde{X}\big\Vert_6^6}{\big\Vert \tilde{X}\big\Vert_4^6}
+ \frac{\sigma_X\big\Vert \tilde{X}\big\Vert_6^3}{\big\Vert \tilde{X}\big\Vert_4^4} 
+
\frac{64\big\Vert
\tilde{Y}^2 \tilde{X}\big\Vert_1\sigma_X}
{\big\Vert \tilde{X}\varepsilon\big\Vert_2^2}
\nonumber\\
&\left.+
 \frac{1120 \big\Vert \tilde{X}\big\Vert_3^3\sigma_Y^2}{\big\Vert \tilde{X}\varepsilon\big\Vert_2^2\sigma_X}+
64\sqrt{2}\frac{\big\Vert \tilde{X}\varepsilon\big\Vert_3^3}{\big\Vert \tilde{X}\varepsilon\big\Vert_2^3}+\frac{288\,\sigma_Y\big\Vert\tilde{X}^2 \tilde{Y}\big\Vert_1}{
\big\Vert \tilde{X}\varepsilon\big\Vert_2^2}+
\frac{432 \vert\tau \vert\big\Vert \tilde{X}\big\Vert_3^3\sigma_Y}{\big\Vert \tilde{X}\varepsilon\big\Vert_2^2}\right.\nonumber\\
&+
\frac{1152\big\Vert \tilde{X}^3\tilde{Y}\big\Vert_1}{\sigma_X^2\big\Vert \tilde{X}\varepsilon\big\Vert_2}+
\frac{1152 \big\Vert \tilde{X}\big\Vert_4^4\sigma_Y}{\sigma_X^3\big\Vert \tilde{X}\varepsilon \big\Vert_2} 
+\frac{576\vert\tau\vert \big\Vert \tilde{X}\big\Vert_4^4}{\sigma_X^2\big\Vert \tilde{X}\epsilon\big\Vert_2} 
+
\frac{576 \big\Vert \tilde{X}^2\tilde{Y}\big\Vert_1\sigma_Y}{\big\Vert \tilde{X}\varepsilon\big\Vert_2^2\sigma_X}\label{termeH2},
\end{align}
\begin{align}
H_3(X,Y)&=\frac{288K\sigma_{\varepsilon}}{\big\Vert \tilde{X}\varepsilon\big\Vert_2^2
\sigma_X}\left[4\big\Vert \tilde{X}^3 \tilde{Y}\big\Vert_1 +4\frac{\sigma_Y \big\Vert \tilde{X}\big\Vert_4^4}{\sigma_X}
+ \vert\tau\vert\big\Vert \tilde{X}\big\Vert_4^4
\right], \label{termeH3}
\end{align}
\begin{align}
H_4(X,Y)& = 8\left[4+\frac{144 \big\Vert \tilde{X}^3 \tilde{Y}\big\Vert_1
\sigma_{\varepsilon}}{
\sigma_X\big\Vert \tilde{X}\varepsilon\big\Vert_2^2
}+\frac{144 \,\sigma_Y \sigma_{\varepsilon}\big\Vert \tilde{X}\big\Vert_4^4}{\big\Vert \tilde{X}\varepsilon\big\Vert_2^2\sigma_X^2}
+\frac{72  \vert\tau\vert\sigma_{\varepsilon}\big\Vert \tilde{X}\big\Vert_4^4 }{\sigma_X\big\Vert \tilde{X}\varepsilon\big\Vert_2^2}
+\frac{c_{1,\rm{\scriptscriptstyle R}}^2\sigma_X^3}{\big\Vert \tilde{X}\big\Vert_3^3}+\frac{2c_{1,\rm{\scriptscriptstyle R}}^2\sigma_X\sigma_Y^2}{\big\Vert \tilde{Y}^2 \tilde{X}\big\Vert_1} \right],\label{termeH4}
\end{align}
\begin{align}
H_5(X,Y) &= 32c_{2,\rm{\scriptscriptstyle R}}^2\left[\frac{\sigma_X\big\Vert \tilde{X}\big\Vert_4^2}{\big\Vert \tilde{X}\big\Vert_3^3 } +
\frac{2\sigma_X\big\Vert \tilde{Y}\big\Vert_4^2}{\big\Vert \tilde{Y}^2 \tilde{X}\big\Vert_1}\right].\label{termeH5}
\end{align}
\subsection{Proof of  Proposition \ref{propositionSureScreeningProperty}}
Let $\gamma$ be defined by $\gamma=\Phi^{-1}(1-q/2)$. We have
\begin{align*}
\mathbb{P}(\mathcal{M}^*\subset \widehat{\mathcal{M}})&=
\mathbb{P}\left(\bigcap_{j\in\mathcal{M}^*}{\frac{\sqrt{n}|\widehat{\tau}_j|}{\sqrt{\widehat{v}_j}}\geq \gamma } \right)=1-\mathbb{P}\left(\bigcup_{j\in\mathcal{M}^*}{\frac{\sqrt{n}|\widehat{\tau}_j|}{\sqrt{\widehat{v}_j}}< \gamma } \right)\geq 1-\sum_{j\in\mathcal{M}^*}{\mathbb{P}\left(\frac{\sqrt{n}|\widehat{\tau}_j|}{\sqrt{\widehat{v}_j}}< \gamma  \right)}.
\end{align*}
For $j\in\llbracket 1,p\rrbracket$, consider the following quantities:
\begin{eqnarray*}
A_{n,j}&=&\frac{1}{n}\sum_{i=1}^n
{\left(X_i^{(j)}-\overline{\X}^{(j)}\right)\varepsilon_i^{(j)}},\qquad
B_{n,j}=\frac{1}{n}\sum_{i=1}^n
{\left(X_i^{(j)}-\overline{\X}^{(j)}\right)^2\left(\widehat\varepsilon_i^{(j)}\right)^2},\\
B^{(j)}&=&\big\| \tilde{X}^{(j)}\varepsilon^{(j)}\big\|_2^2, \qquad \mbox{ and }\Var_{n}\left(\X^{(j)}\right)=\frac{1}{n}\sum_{i=1}^n
{\left(X_i^{(j)}-\overline{\X}^{(j)}\right)^2}.\end{eqnarray*}
Then
\begin{eqnarray*}
   \frac{\sqrt{n}|\widehat{\tau}_j|}{\sqrt{\widehat{v}_j}}\geq
   \frac{\sqrt{n}|\tau_j|}{\sqrt{\widehat{v}_j}}-\frac{\sqrt{n}|\widehat{\tau}_j-\tau_j|}{\sqrt{\widehat{v}_j}} \geq \frac{\sqrt{n}\Var_{n}\left(\X^{(j)}\right)|\tau_j|}{\sqrt{B_{n,j}}}-\frac{\sqrt{n}|A_{n,j}|}{\sqrt{B_{n,j}}}.
\end{eqnarray*}
Consequently, the probability
$$\mathbb{P}\left(\frac{\sqrt{n}|\widehat{\tau}_j|}{\sqrt{\widehat{v}_j}}< \gamma  \right)$$ is bounded by

\begin{eqnarray*}    
  \mathbb{P}\left( \left( \frac{\sqrt{n}\Var_{n}\left(\X^{(j)}\right)|\tau_j|}{\sqrt{B_{n,j}}}-\frac{\sqrt{n}|A_{n,j}|}{\sqrt{B_{n,j}}}    < \gamma\right) \bigcap \left(\frac{B^{(j)}}{2}\leq B_{n,j}\leq 2B^{(j)}
  \right)\bigcap \left(\Var_{n}\left(\X^{(j)}\right)\geq \frac{\sigma_{X^{(j)}}^2}{2}
  \right)\right)\nonumber\\ \!\!\! \!\!\!
  +\mathbb{P}\left(B_{n,j}< \frac{B^{(j)}}{2}
  \right)+\mathbb{P}\left(B_{n,j}> 2B^{(j)} \right)+\mathbb{P}\left(\Var_{n}\left(\X^{(j)}\right)< \frac{\sigma_{X^{(j)}}^2}{2}
  \right)
  \end{eqnarray*}
  which is less than
  \begin{eqnarray}
    \mathbb{P}\left( \frac{\sqrt{n}|A_{n,j}|}{\sqrt{B^{(j)}}}    > \frac{\sqrt{n}\sigma_{X^{(j)}}^2|\tau_j|}{4\sqrt{B^{(j)}}}-\frac{\gamma}{\sqrt{2}}\right) 
  +\mathbb{P}\left(|B_{n,j}-B^{(j)}|> \frac{B^{(j)}}{2}
  \right)+\mathbb{P}\left(\Var_{n}\left(\X^{(j)}\right)< \frac{\sigma_{X^{(j)}}^2}{2}\label{borne_tau}
  \right).\label{borne_tau_chapo}
\end{eqnarray}
The first probability in (\ref{borne_tau_chapo}) is bounded by applying Lemma \ref{prelim}. For $K_1$ and $K_2$ as in (\ref{kappa1_kappa2}) we have
\begin{multline*}
  \mathbb{P}\left( \frac{\sqrt{n}|A_{n,j}|}{\sqrt{B^{(j)}}}    > \frac{\sqrt{n}\sigma_{X^{(j)}}^2|\tau_j|}{4\sqrt{B^{(j)}}}-\frac{\gamma}{\sqrt{2}}\right) \\\leq 
  2\mathbb{P}\left(\mathcal{N}(0,1)\geq \frac{\sqrt{n}\sigma_{X^{(j)}}^2|\tau_j|}{4\sqrt{B^{(j)}}}-\frac{\gamma}{\sqrt{2}}\right)
  +\frac{2K_1(X^{(j)},Y)\log(n)}{\sqrt{n}}+\frac{2K_2(X^{(j)},Y)}{\sqrt{n}}.
\end{multline*}
Then, according to  Pollak Inequality (\ref{Pollak}), we write
\begin{multline*}
  \mathbb{P}\left( \frac{\sqrt{n}|A_{n,j}|}{\sqrt{B^{(j)}}}    > \frac{\sqrt{n}\sigma_{X^{(j)}}^2|\tau_j|}{4\sqrt{B^{(j)}}}-\frac{\gamma}{\sqrt{2}}\right) \\\leq  
  \exp\left(-\frac{1}{2}\left(\frac{\sqrt{n}\sigma_{X^{(j)}}^2|\tau_j|}{4\sqrt{B^{(j)}}}-\frac{\gamma}{\sqrt{2}}\right)^2\,\right)
  +\frac{2K_1(X^{(j)},Y)\log(n)}{\sqrt{n}}+\frac{2K_2(X^{(j)},Y)}{\sqrt{n}}.
\end{multline*}
The second probability in (\ref{borne_tau_chapo}) is bounded by applying Lemma \ref{lemme2}. For $G_1$, $G_2$ and $G_3$ as in (\ref{termeG1}), (\ref{termeG2}), (\ref{termeG3}) we have
\begin{align*}
\mathbb{P}\left(|B_{n,j}-B^{(j)}|> \frac{B^{(j)}}{2}
  \right)\leq 
\frac{G_1(X^{(j)},Y)}{\sqrt{n}}+ \frac{4G_2(X^{(j)},Y)}{n}+\frac{G_3(X^{(j)},Y)}{n^{3/2}}. 
\end{align*}
Finally, the last probability in (\ref{borne_tau_chapo}) is bounded by applying Lemma \ref{lemmemajoration_var_X}.
Consequently,
\begin{multline*}
\mathbb{P}(\mathcal{M}^*\subset \widehat{\mathcal{M}})\geq
1-\sum_{j\in\mathcal{M}^*}{\left[
  \exp\left(-\frac{1}{4}\left(\frac{\sqrt{n}\sigma_{X^{(j)}}^2|\tau_j|}{2\sqrt{2}\big\| \tilde{X}^{(j)}\varepsilon^{(j)}\big\|_2}-\gamma\right)^2\,\right)\right]} \\
    - \left\vert\mathcal{M}^*\right\vert\max_{j\in\mathcal{M}^*}
  \Bigg[
  \frac{2K_1(X^{(j)},Y)\log(n)}{\sqrt{n}}+\frac{1}{\sqrt{n}}\Bigg(2K_2(X^{(j)},Y)+\frac{2\big\Vert \tilde{X}^{(j)}\big\Vert_4^2}{\sigma_{X^{(j)}}^2} +G_1(X^{(j)},Y)\Bigg)
  \\\qquad + 
  \frac{4G_2(X^{(j)},Y)+2}{n}+\frac{G_3(X^{(j)},Y)}{n^{3/2}}\Bigg],
  \end{multline*}
  and Proposition \ref{propositionSureScreeningProperty} is proved.
\subsection{Proof of  Proposition \ref{propositionFPR}}
We write
\begin{align*}
\mathbb{E}\Bigg[\frac{|{\mathcal{M}^*}^c\cap \widehat{\mathcal{M}}|}{|{\mathcal{M}^*}^c|}\Bigg]
&=\frac{1}{{|{\mathcal{M}^*}^c|}}\mathbb{E}\Bigg[
\sum_{ j\in {\mathcal{M}^*}^c}\ind_{\frac{\sqrt{n}|\hat \tau_j|}{\sqrt{\widehat{v}_j}}\geq \gamma}
\Bigg]=\frac{1}{{|{\mathcal{M}^*}^c|}}\sum_{ j\in {\mathcal{M}^*}^c} \mathbb{P}\Bigg[   \frac{\sqrt{n}|\hat \tau_j|}{\sqrt{\widehat{v}_j}}\geq \Phi^{-1}(1-q/2)    \Bigg]\\
&=
\frac{1}{{|{\mathcal{M}^*}^c|}}\sum_{ j\in {\mathcal{M}^*}^c} \Bigg[\mathbb{P}\Bigg[   \frac{\sqrt{n}|\hat \tau_j|}{\sqrt{\widehat{v}_j}}\geq \Phi^{-1}(1-q/2)     \Bigg]-\mathbb{P}(|\mathcal{N}(0,1)|\geq \Phi^{-1}(1-q/2))\Bigg]\\& \qquad+\mathbb{P}\left(|\mathcal{N}(0,1)|\geq \Phi^{-1}(1-q/2)\right)\\
&=\frac{1}{{|{\mathcal{M}^*}^c|}}\sum_{ j\in {\mathcal{M}^*}^c} \Bigg[\mathbb{P}\Bigg[   \frac{\sqrt{n}|\hat \tau_j|}{\sqrt{\widehat{v}_j}}\geq \Phi^{-1}(1-q/2)      \Bigg]-\mathbb{P}(|\mathcal{N}(0,1)|\geq \Phi^{-1}(1-q/2))\Bigg]+q.
\end{align*}
The proof of Proposition \ref{propositionFPR} is complete by applying Proposition \ref{propositionattenduesansj}.

\subsection{Proof of Lemma \ref{prelim}}

Note first that, applying Berry-Esseen Inequality (\ref{BE}), 
\begin{equation}
\label{etape1BE}
\sup_{x\in \mathbb{R}}\left\vert\mathbb{P}\left(\frac{\sum_{i=1}^n \tilde{X}_i\varepsilon_i}
{
\sqrt{n}
\Vert\tilde{X}\varepsilon\Vert_2
}\leq x \right)-\Phi(x)\right\vert\leq c_{\rm{\scriptscriptstyle BE}}\frac{    \Vert\tilde{X}\varepsilon\Vert_3^3  }{ 
\Vert\tilde{X}\varepsilon\Vert_2^{3}}\frac{1}{\sqrt{n}}.
\end{equation}
Next we write that
\begin{equation*}
\frac{\sum_{i=1}^n
(X_i-\overline{\X})\varepsilon_i}{\sqrt{n} \Vert\tilde{X}\varepsilon\Vert_2}=
\frac{\sum_{i=1}^n \tilde{X_i}
\varepsilon_i}{\sqrt{n} \Vert\tilde{X}\varepsilon\Vert_2
}
+\frac{(\mu_X-\overline{\X})\sum_{i=1}^n \varepsilon_i}{\sqrt{n}
\Vert\tilde{X}\varepsilon\Vert_2}.
\end{equation*}
Then, we apply Lemma \ref{CRAS} with  $$R_n=
\frac{(\mu_X-\overline{\X})\sum_{i=1}^n \varepsilon_i}{
\Vert\tilde{X}\varepsilon\Vert_2} \mbox{ and }T_n=\frac{\sum_{i=1}^n
\tilde{X}_i\varepsilon_i}{\Vert\tilde{X}\varepsilon\Vert_2}.$$
According to (\ref{etape1BE})
\begin{equation*}
\sup_{x\in \mathbb{R}}\left\vert\mathbb{P}\left(\frac{T_n}
{
\sqrt{n}
}\leq x \right)-\Phi(x)\right\vert\leq c_{\rm{\scriptscriptstyle BE}}\frac{ \Vert\tilde{X}\varepsilon\Vert_3^3 
  }{ \Vert\tilde{X}\varepsilon\Vert_2^3} \frac{1}{\sqrt{n}}.
\end{equation*}
Let  $(b_n)_{n\in \mathbb{N}}$ and $(d_n)_{n\in \mathbb{N}}$ be two real and positive sequences  tending to infinity. Applying Lemma \ref{lemmeproduit}, we get 
\begin{equation*}
\mathbb{P}\left(\frac{|R_n|}{\sqrt{n}}\geq \frac{1}{b_n}\right)\leq 
\mathbb{P}\left(\frac{\left |\sum_{i=1}^n \varepsilon_i\right |}{\sqrt{n}\sigma_{\varepsilon}}\geq d_n\right)
+ \mathbb{P}\left(d_n\frac{\left |\overline{\X}-\mu_X\right |\sigma_{\varepsilon}}{\Vert\tilde{X}\varepsilon\Vert_2
}       \geq \frac{1}{b_n}\right).
\end{equation*}
Applying Berry-Esseen  Inequality (\ref{BE}) and Pollak Inequality (\ref{Pollak}),
\begin{align*}
\mathbb{P}\left(\frac{\left |\sum_{i=1}^n \varepsilon_i\right |}{\sqrt{n}\sigma_{\varepsilon}}\geq d_n\right)
&\leq \mathbb{P}\left(\left |\mathcal{N}(0,1)\right |\geq d_n\right) + 2c_{\rm{\scriptscriptstyle BE}}\frac{ \left\Vert\varepsilon \right\Vert_3^3}
{\sigma_{\varepsilon}^3}\frac{1}{\sqrt{n}} 
 \leq \exp{\left(-\frac{d_n^2}{2} \right)}
+ 2c_{\rm{\scriptscriptstyle BE}}\frac{\left\Vert\varepsilon\right\Vert_3^3
}{\sigma_{\varepsilon}^3}\frac{1}{\sqrt{n}}.
\end{align*}
Thus we take $d_n=\sqrt{\log(n)}$
and we obtain 
$$
\mathbb{P}\left(\frac{\left |\sum_{i=1}^n \varepsilon_i\right |}{\sqrt{n}\sigma_{\varepsilon}}\geq d_n\right)
\leq \frac{1}{\sqrt{n}} \left(1+2c_{\rm{\scriptscriptstyle BE}}\frac{
\left\Vert\varepsilon \right\Vert_3^3
}{\sigma_{\varepsilon}^3}\right).
$$
Applying   Berry-Esseen Inequality (\ref{BE}) and  Pollak Inequality (\ref{Pollak}), we get
\begin{align*}
\mathbb{P}\left(d_n\frac{\left |\overline{\X}-\mu_X\right |\sigma_{\varepsilon}}{ \Vert\tilde{X}\varepsilon\Vert_2
}       \geq \frac{1}{b_n}\right)
 & = \mathbb{P}\left(\frac{\sqrt{n}\left |\overline{\X}-\mu_X\right |}{\sigma_X}\geq \frac{\sqrt{n} \Vert\tilde{X}\varepsilon\Vert_2
 }{b_n\sqrt{\log(n)}\sigma_{\varepsilon} \sigma_X }\right) \\
 &\leq \mathbb{P}\left(\left |\mathcal{N}(0,1)\right |\geq \frac{\sqrt{n}
 \Vert\tilde{X}\varepsilon\Vert_2
 }{b_n\sqrt{\log(n)}\sigma_{\varepsilon} \sigma_X }
 \right) + 2c_{\rm{\scriptscriptstyle BE}}\frac{\Vert\tilde{X}\Vert_3^3
 }{\sigma_X^{3}}\frac{1}{\sqrt{n}} \\
 & \leq \exp{\left(-  \frac{n \Vert\tilde{X}\varepsilon\Vert_2^2
 }{2b_n^2\log(n)\sigma_{\varepsilon}^2 \sigma_X^2 }
\right)}
+2c_{\rm{\scriptscriptstyle BE}}\frac{ \Vert\tilde{X}\Vert_3^3
}{\sigma_X^{3}}\frac{1}{\sqrt{n}} .
 \end{align*}
Hence, taking
 \begin{align*}
 b_n 
 & =  \frac{\sqrt{n}}{\log(n)}
  \frac{ \Vert\tilde{X}\varepsilon\Vert_2
  }{\sigma_{\varepsilon}\sigma_X},
   \end{align*}
we obtain
\begin{equation*}
\mathbb{P}\left(\frac{|R_n|}{\sqrt{n}}\geq \frac{1}{b_n}\right)
 \leq \frac{2}{\sqrt{n}}\left[1+
c_{\rm{\scriptscriptstyle BE}}\left(\frac{\Vert\varepsilon\Vert_3^3
}{\sigma_{\varepsilon}^3}+\frac{\Vert\tilde{X}\Vert_3^3}
{\sigma_X^{3}}\right)\right].
 \end{equation*}
Applying Lemma \ref{CRAS}, we get that
\begin{multline*}
\sup_{x\in \mathbb{R}} \Big\vert\mathbb{P}  \left(\frac{\sum_{i=1}^n
(X_i-\overline{\X})\varepsilon_i}
{
\sqrt{n} \Vert\tilde{X}\varepsilon\Vert_2
}\leq x \right) - \Phi(x)\Big\vert \\
\leq
c_{\rm{\scriptscriptstyle BE}}\frac{  \Vert\tilde{X}\varepsilon\Vert_3^3
  }{ \Vert\tilde{X}\varepsilon\Vert_2^3}\frac{1}{\sqrt{n}}+
\frac{2}{\sqrt{n}}\left[1+
c_{\rm{\scriptscriptstyle BE}}\left(\frac{ \left\Vert\varepsilon\right\Vert_3^3
}{\sigma_{\varepsilon}^3}+\frac{\Vert\tilde{X}\Vert_3^3
}{\sigma_X^{3}}\right)\right] +\frac{\log(n)}{\sqrt{n}}\frac{\sigma_{\varepsilon}\sigma_X}{\sqrt{2\pi}\Vert\tilde{X}\varepsilon\Vert_2}, 
\end{multline*}
and Lemma \ref{prelim} is proved with $K_1$ and $K_2$ given by (\ref{kappa1_kappa2}).

\subsection{Proof of  Proposition \ref{propositionattenduesansj}}
Let $A_n$, $B_n$ and $B$ be defined by
\begin{eqnarray}
A_n=\frac{1}{n}\sum_{i=1}^n
{\left(X_i-\overline{\X}\right)\varepsilon_i},\quad  
B_n=\frac{1}{n}\sum_{i=1}^n
{\left(X_i-\overline{\X}\right)^2\left(\widehat\varepsilon_i\right)^2}, \quad
B=\Vert\tilde{X}\varepsilon\Vert_2^2.
\label{notations_An_Bn_B}
\end{eqnarray}
Applying Lemma \ref{prelim}, we get 
\begin{equation}
\label{inegalite2_Lemme}
\sup_{x\in \mathbb{R}}\left\vert\mathbb{P}\left(\sqrt{n}\frac{A_n}
{
\sqrt{B}
}\leq x \right)-\Phi(x)\right\vert\leq \frac{K_1(X,Y)\log(n)}{\sqrt{n}} + \frac{K_2(X,Y)}{\sqrt{n}}.
\end{equation}
Our aim is to prove an inequality similar to  (\ref{inegalite2_Lemme}) for
\begin{equation*}
\sup_{x\in \mathbb{R}}\left\vert\mathbb{P}\left(\sqrt{n}\frac{A_n}
{
\sqrt{B_n}
}\leq x \right)-\Phi(x)\right\vert.
\end{equation*}
We start by writing
\begin{align*}
\sqrt{n}\frac{A_n}
{
\sqrt{B_n}}&=
\sqrt{n}\frac{A_n}
{
\sqrt{B}}
+
\sqrt{n}\frac{A_n}
{
\sqrt{B}}\left(\frac{\sqrt{B}}
{
\sqrt{B_n}}-1
\right) = \sqrt{n}\frac{A_n}
{
\sqrt{B}}
+\frac{R_n}{\sqrt{n}}.
\end{align*}
where $$R_n=\frac{nA_n}
{
\sqrt{B}}\left(\frac{\sqrt{B}}
{
\sqrt{B_n}}-1
\right),$$ 
To apply Lemma \ref{CRAS}, according to (\ref{inegalite2_Lemme}), we take
\begin{eqnarray}
\displaystyle a_n=\frac{\sqrt{n}}{K_1(X,Y)\log(n)+ K_2(X,Y)}.
\label{a_n}
\end{eqnarray}
We aim at finding $b_n$ and $c_n$ such that
$$
\mathbb{P}\left(\frac{\vert R_n\vert}
{
\sqrt{n}}\geq\frac{1}{b_n}\right)\leq \frac{1}{c_n}.
$$
For $d_n\rightarrow \infty$, applying Lemma \ref{lemmeproduit}, we get
\begin{align}
\mathbb{P}\left(\frac{\vert R_n\vert}
{
\sqrt{n}}\geq\frac{1}{b_n}\right) &= \mathbb{P}\left(\frac{\sqrt{n}\vert A_n\vert}
{
\sqrt{B}}\left\vert\frac{\sqrt{B}}
{
\sqrt{B_n}}-1 \right\vert  \geq\frac{1}{b_n}\right) \nonumber\\
& \leq \mathbb{P}\left(\frac{\sqrt{n}\vert A_n\vert}
{
\sqrt{B}}\geq d_n\right) + \mathbb{P}\left(d_n\left\vert\frac{\sqrt{B}}
{
\sqrt{B_n}}-1 \right\vert  \geq\frac{1}{b_n}\right).
\label{majoration_reste_prop}
\end{align}
According to (\ref{inegalite2_Lemme}) and Pollak Inequality (\ref{Pollak}), for the first term in (\ref{majoration_reste_prop}) we obtain 
\begin{align*}
\mathbb{P}\left(\frac{\sqrt{n}\vert A_n\vert}
{
\sqrt{B}}\geq d_n\right) & \leq \mathbb{P}\left(\vert \mathcal{N}(0,1)\vert \geq d_n\right)+
 \frac{2K_1(X,Y)\log(n)}{\sqrt{n}} + \frac{2K_2(X,Y)}{\sqrt{n}} \\
& \leq  \exp\left(-\frac{d_n^2}{2}\right)+ \frac{2K_1(X,Y)\log(n)}{\sqrt{n}} + \frac{2K_2(X,Y)}{\sqrt{n}} .
\end{align*}
Thus we choose $d_n = \sqrt{\log(n)}$ and we get
\begin{align}
\mathbb{P}\left(\frac{\sqrt{n}\vert A_n\vert}
{
\sqrt{B}}\geq d_n\right) 
& \leq  \frac{2K_1(X,Y)\log(n)}{\sqrt{n}} + \frac{2K_2(X,Y)+1}{\sqrt{n}}. 
\label{majoration_terme1_reste_prop}
\end{align}
For the second term in (\ref{majoration_reste_prop}) we write 
\begin{align*}
\mathbb{P}\left(\left\vert\frac{\sqrt{B}}
{
\sqrt{B_n}}-1 \right\vert  \geq\frac{1}{b_n\sqrt{\log(n)}}\right)&=
\mathbb{P}\left(\left\vert\sqrt{B}-\sqrt{B_n} \right\vert  \geq\frac{\sqrt{B_n}}{b_n\sqrt{\log(n)}}\right) \nonumber \\
&\leq 
\mathbb{P}\left[\left(\left\vert\sqrt{B}-\sqrt{B_n} \right\vert  \geq\frac{\sqrt{B/2}}{b_n\sqrt{\log(n)}}\right)\bigcap\left(B_n  > \frac{B}{2} \right)\right]+\mathbb{P}\left(B_n\leq \frac{B}{2}\right).
\end{align*}
First we use the following inequality
$$\left\vert \sqrt{B}-\sqrt{B_n} \right\vert \leq \frac{1}{2\sqrt{B/2}}\left\vert B-B_n \right\vert $$
valid for $B_n> \frac{B}{2}$. 
Then we get
\begin{align*}
\mathbb{P}\left[\left(\left\vert\sqrt{B}-\sqrt{B_n} \right\vert  \geq\frac{\sqrt{B/2}}{b_n\sqrt{\log(n)}}\right)\bigcap\left(B_n  > \frac{B}{2} \right)\right]
&\leq 
\mathbb{P}\left(\left\vert B_n-B \right\vert  \geq\frac{B}{b_n\sqrt{\log(n)}}\right). \\
\end{align*}
Besides 
\begin{align*}
\mathbb{P}\left(B_n\leq \frac{B}{2}\right)
\leq \mathbb{P}\left(\left\vert B_n-B\right\vert \geq \frac{B}{2}\right).
\end{align*}
Hence we get
\begin{align}
\label{majoration_terme2_reste_prop}
\mathbb{P}\left(\left\vert\frac{\sqrt{B}}
{
\sqrt{B_n}}-1 \right\vert  \geq\frac{1}{b_n\sqrt{\log(n)}}\right)
\leq P_1+P_2
\end{align}
where
\begin{align}
\label{P1-P2}
P_1=
\mathbb{P}\left(\left\vert B_n-B \right\vert  \geq\frac{B}{b_n\sqrt{\log(n)}}\right) \mbox{ and } P_2=\mathbb{P}\left(\left\vert B_n-B\right\vert \geq \frac{B}{2}\right).
\end{align}
\noindent
We now choose 
\begin{eqnarray}
b_n=\frac{K\sqrt{n}}{\log(n)}
\label{b_n}
\end{eqnarray} 
where $K$ is given by  \eqref{K_final}. We have to bound up the two 
probabilities $P_1$ and $P_2$ in (\ref{P1-P2}). 
The probability $P_2$ is controlled by applying Lemma \ref{lemme2} (whose proof is given in Section \ref{section5.8}), and the probability $P_1$ is controlled in Lemma \ref{lemme_majoration_P1} below (whose proof is given in Section \ref{preuve_lemme_majorationP1}).
\begin{Lemma}
\label{lemme_majoration_P1}
Let $K$ be defined by \eqref{K_final}. The term $P_1$ given by \eqref{P1-P2} is controlled as follows.
\begin{align}
P_1&\leq \frac{1}{\sqrt{n}}\left[12+\frac{1}{2\big\Vert \tilde{Y}^2 \tilde{X}\big\Vert_1}\left(2\big\Vert\tilde{Y}^2\tilde{X}\big\Vert_2+\sigma_X\big\Vert \tilde{Y}\big\Vert_4^2+
4\sigma_Y
\big\Vert\tilde{X}\tilde{Y}\big\Vert_2\right) \right.\nonumber \\ 
&\left.+2c_{\rm{\scriptscriptstyle BE}}\left(\frac{
\big\Vert \tilde{X}\varepsilon\big\Vert_6^6}{\left[\Var\left(\tilde{X}^2\varepsilon^2\right)\right]^{3/2}}+3\frac{\big\Vert \tilde{Y}\big\Vert_3^3}{\sigma_Y^3}+
3\frac{
\big\Vert\tilde{X} \varepsilon\big\Vert_3^3}{\big\Vert \tilde{X}\varepsilon\big\Vert_2^3
}+5\frac{\big\Vert \tilde{X}\big\Vert_3^3}{\sigma_X^3}\right)+
2\frac{\big\Vert\tilde{X}^2 \tilde{Y}\big\Vert_2}{\big\Vert \tilde{X}^2\tilde{Y}\big\Vert_1}
+ \frac{\sigma_Y\big\Vert \tilde{X}\big\Vert_4^2}{\big\Vert \tilde{X}^2\tilde{Y}\big\Vert_1}+ \frac{13}{4}\frac{\sigma_X \big\Vert \tilde{X}\big\Vert_4^2}{\big\Vert \tilde{X}\big\Vert_3^3}
\right. \nonumber\\
& \qquad + \frac{9}{2}\frac{\big\Vert \tilde{X}\big\Vert_6^3}{\big\Vert \tilde{X}\big\Vert_3^3} 
+ \left.6\frac{\big\Vert \tilde{Y}\big\Vert_4^2}{\sigma_Y^2}
+ 16\frac{\big\Vert \tilde{X}\big\Vert_4^2}{\sigma_X^2}+ 4\frac{\big\Vert\tilde{X}^3\tilde{Y}\big\Vert_{3/2}^{3/2}}{\big\Vert \tilde{X}^3 \tilde{Y}\big\Vert_1^{3/2}}+ \frac{\sigma_Y\big\Vert \tilde{X}\big\Vert_6^3}{4\big\Vert \tilde{X}^3 \tilde{Y}\big\Vert_1}
+8\frac{\big\Vert \tilde{X}\big\Vert_6^6}{\big\Vert \tilde{X}\big\Vert_4^6}
+ \frac{\sigma_X\big\Vert \tilde{X}\big\Vert_6^3}{2\big\Vert \tilde{X}\big\Vert_4^4} \right]\nonumber\\
& \qquad+ \frac{1}{\sqrt{n\log(n)}}\frac{288K\sigma_{\varepsilon}}{\big\Vert \tilde{X}\varepsilon\big\Vert_2^2
\sigma_X}\left[4\big\Vert \tilde{X}^3 \tilde{Y}\big\Vert_1 +4\frac{\sigma_Y \big\Vert \tilde{X}\big\Vert_4^4}{\sigma_X}
+ \vert\tau\vert\big\Vert \tilde{X}\big\Vert_4^4
\right] \nonumber \\
&\qquad +\frac{2}{n}\left(8+\frac{2c_{1,\rm{\scriptscriptstyle R}}^2\sigma_X^3}{\big\Vert\tilde{X}\big\Vert_3^3}+\frac{4c_{1,\rm{\scriptscriptstyle R}}^2\sigma_X\sigma_Y^2}{\big\Vert\tilde{Y}^2 \tilde{X}\big\Vert_1}\right) +\frac{16c_{2,\rm{\scriptscriptstyle R}}^2}{n^{3/2}}\left(\frac{\sigma_X\big\Vert\tilde{X}\big\Vert_4^2}{\big\Vert\tilde{X}\big\Vert_3^3 } +
\frac{2\sigma_X\big\Vert\tilde{Y}\big\Vert_4^2}{\big\Vert\tilde{Y}^2 \tilde{X}\big\Vert_1}\right). 
\label{majoration_terme2_1_reste_prop}
\end{align}
\end{Lemma}
We come back to the study of the term (\ref{majoration_reste_prop}), by using the inequalities (\ref{majoration_terme1_reste_prop}),  (\ref{majoration_terme2_reste_prop}) and (\ref{majoration_terme2_1_reste_prop}) as well as Inequality  (\ref{majoration_terme2_2_reste_prop}) in Lemma \ref{lemme2}. Let $K_1$, $K_2$ be given by (\ref{kappa1_kappa2}),  $H_3$, $H_4$, $H_5$ be defined in Section \ref{definitionG_ietH_i} by (\ref{termeH3}), (\ref{termeH4}) and (\ref{termeH5}) and $\widetilde{H}_2$ be defined by 
\begin{align*}
\widetilde{H}_2(X,Y)&=12+\frac{1}{\big\Vert \tilde{Y}^2 \tilde{X}\big\Vert_1}\left(2\big\Vert\tilde{Y}^2\tilde{X}\big\Vert_2+\sigma_X\big\Vert \tilde{Y}\big\Vert_4^2+
4\sigma_Y
\big\Vert\tilde{X}\tilde{Y}\big\Vert_2\right) \nonumber \\ 
&\left.+2c_{\rm{\scriptscriptstyle BE}}\left(\frac{
\big\Vert \tilde{X}\varepsilon\big\Vert_6^6}{\left[\Var\left(\tilde{X}^2\varepsilon^2\right)\right]^{3/2}}+3\frac{\big\Vert \tilde{Y}\big\Vert_3^3}{\sigma_Y^3}+
3\frac{
\big\Vert\tilde{X} \varepsilon\big\Vert_3^3}{\big\Vert \tilde{X}\varepsilon\big\Vert_2^3
}+5\frac{\big\Vert \tilde{X}\big\Vert_3^3}{\sigma_X^3}\right)+
4\frac{\big\Vert\tilde{X}^2 \tilde{Y}\big\Vert_2}{\big\Vert \tilde{X}^2\tilde{Y}\big\Vert_1}
+ 2\frac{\sigma_Y\big\Vert \tilde{X}\big\Vert_4^2}{\big\Vert \tilde{X}^2\tilde{Y}\big\Vert_1}
\right. \nonumber\\
&\qquad + \frac{13}{2}\frac{\sigma_X \big\Vert \tilde{X}\big\Vert_4^2}{\big\Vert \tilde{X}\big\Vert_3^3}
+ 9\frac{\big\Vert \tilde{X}\big\Vert_6^3}{\big\Vert \tilde{X}\big\Vert_3^3} 
+ 12\frac{\big\Vert \tilde{Y}\big\Vert_4^2}{\sigma_Y^2}
+ 32\frac{\big\Vert \tilde{X}\big\Vert_4^2}{\sigma_X^2}+ 2^{11/2}\frac{\big\Vert\tilde{X}^3\tilde{Y}\big\Vert_{3/2}^{3/2}}{\big\Vert \tilde{X}^3 \tilde{Y}\big\Vert_1^{3/2}}+ \frac{3}{4}\frac{\sigma_Y\big\Vert \tilde{X}\big\Vert_6^3}{\big\Vert \tilde{X}^3 \tilde{Y}\big\Vert_1} \nonumber\\
&\left.+16\frac{\big\Vert \tilde{X}\big\Vert_6^6}{\big\Vert \tilde{X}\big\Vert_4^6}
+ \frac{\sigma_X\big\Vert \tilde{X}\big\Vert_6^3}{\big\Vert \tilde{X}\big\Vert_4^4} 
+
\frac{64\big\Vert
\tilde{Y}^2 \tilde{X}\big\Vert_1\sigma_X}
{\big\Vert \tilde{X}\varepsilon\big\Vert_2^2}+
 \frac{1120 \big\Vert \tilde{X}\big\Vert_3^3\sigma_Y^2}{\big\Vert \tilde{X}\varepsilon\big\Vert_2^2\sigma_X}+
64\sqrt{2}\frac{\big\Vert \tilde{X}\varepsilon\big\Vert_3^3}{\big\Vert \tilde{X}\varepsilon\big\Vert_2^3}+\frac{288\,\sigma_Y\big\Vert\tilde{X}^2 \tilde{Y}\big\Vert_1}{
\big\Vert \tilde{X}\varepsilon\big\Vert_2^2}\right.\nonumber\\
&+
\frac{432 \vert\tau \vert\big\Vert \tilde{X}\big\Vert_3^3\sigma_Y}{\big\Vert \tilde{X}\varepsilon\big\Vert_2^2}+
\frac{1152\big\Vert \tilde{X}^3\tilde{Y}\big\Vert_1}{\sigma_X^2\big\Vert \tilde{X}\varepsilon\big\Vert_2}+
\frac{1152 \big\Vert \tilde{X}\big\Vert_4^4\sigma_Y}{\sigma_X^3\big\Vert \tilde{X}\varepsilon \big\Vert_2} 
+\frac{576\vert\tau\vert \big\Vert \tilde{X}\big\Vert_4^4}{\sigma_X^2\big\Vert \tilde{X}\epsilon\big\Vert_2} 
+
\frac{576 \big\Vert \tilde{X}^2\tilde{Y}\big\Vert_1\sigma_Y}{\big\Vert \tilde{X}\varepsilon\big\Vert_2^2\sigma_X}. \nonumber
\end{align*}
We get 
\begin{align*}
\mathbb{P}\left(\frac{\vert R_n\vert}
{
\sqrt{n}}\geq\frac{1}{b_n}\right)\leq \frac{1}{c_n}
\end{align*}
where 
\begin{align}
\frac{1}{c_n}= \frac{2K_1(X,Y)\log(n)}{\sqrt{n}} + \frac{2K_2(X,Y)+\widetilde{H}_2(X,Y)+1}{\sqrt{n}}+ \frac{H_3(X,Y)}{\sqrt{n\log(n)}}+ \frac{H_4(X,Y) }{n}+\frac{H_5(X,Y)}{n^{3/2}}.
\label{c_n}
\end{align}
Then by applying Lemma \ref{CRAS}, with $a_n$, $b_n$ and $c_n$ given by (\ref{a_n}), (\ref{b_n}) and (\ref{c_n}) respectively, we obtain that  
\begin{align*}
\sup_{x\in \mathbb{R}}\left\vert\mathbb{P}\left(\sqrt{n}\frac{A_n}
{
\sqrt{B_n}
}\leq x \right)-\Phi(x)\right\vert
&\leq \frac{\log(n)H_1(X,Y)}{\sqrt{n}}+ \frac{H_2(X,Y)}{\sqrt{n}}
+ \frac{H_3(X,Y)}{\sqrt{n\log(n)}}+ \frac{H_4(X,Y) }{n}+\frac{H_5(X,Y)}{n^{3/2}},
\end{align*}
which completes the proof of Proposition \ref{propositionattenduesansj}.
\subsection{Proof of Lemma \ref{lemme_majoration_P1}}
\label{preuve_lemme_majorationP1}
Let $P_1$ be defined by \eqref{P1-P2}. Recall that $b_n$ is defined by (\ref{b_n}) with $K$ defined in (\ref{K_final}). 
We start by writing
\begin{align}
\label{Majoration_P1}
P_1
&\leq
D_1+D_2
\end{align}
with
\begin{align*}
D_1=\mathbb{P}\Bigg[ \left\vert    \frac{1}{n}\sum_{i=1}^n\big[(X_i-\overline{\X})^2-\tilde{X}_i^2\Big]{\widehat{\varepsilon_i}}^2
  \right\vert    \geq \frac{B}{2b_n\sqrt{\log(n)}}        \Bigg],
  \,
D_2=\mathbb{P} &\Bigg[   \left\vert   \frac{1}{n}\sum_{i=1}^n\tilde{X}_i^2{\widehat{\varepsilon_i}}^2 -\Vert\tilde{X}\varepsilon\Vert_2^2\right\vert \geq\frac{B}{2b_n\sqrt{\log(n)}}  \Bigg]\! .
\end{align*}
\subsubsection*{Study of the term $D_1$}
We have
\begin{align*}
D_1\leq  \mathbb{P}\Bigg[ \left\vert  \overline{\X}-\mu_{X}\right\vert \frac{1}{n}\sum_{i=1}^n   {\widehat{\varepsilon}_i}^2  \left\vert  2X_i-\overline{\X}-\mu_{X}\right\vert \geq\frac{B\sqrt{\log(n)}}{2K\sqrt{n}}  \Bigg].
\end{align*}
Using the following bound for ${\widehat{\varepsilon_i}}^2$
\begin{align}
\label{varepsilonchapo2}
{\widehat{\varepsilon_i}}^2=
[(Y_i-\overline{\Y})-\widehat{\tau}(X_i-\overline{\X})]^2\leq 2(Y_i-\overline{\Y})^2+2\widehat{\tau}^2(X_i-\overline{\X})^2,
\end{align}
we get 
\begin{align*}
D_1\leq D_{1,1}+D_{1,2}
\end{align*}
with
\begin{align*}
D_{1,1}&=\mathbb{P} \Bigg[  \left\vert  \overline{\X}-\mu_{X}\right\vert    \frac{1}{n}\sum_{i=1}^n(Y_i-
\overline{\Y})^2\left\vert 2X_i-\overline{\X}-\mu_{X}\right\vert \geq
\frac{B\sqrt{\log(n)}}{8K\sqrt{n}}     \Bigg]
\end{align*}
and 
\begin{align*}
D_{1,2}&=\mathbb{P} \Bigg[  \left\vert  \overline{\X}-\mu_{X}\right\vert  \widehat{\tau}^2  \frac{1}{n}\sum_{i=1}^n(X_i-
\overline{\X})^2\left\vert 2X_i-\overline{\X}-\mu_{X}\right\vert \geq
\frac{B\sqrt{\log(n)}}{8K\sqrt{n}}     \Bigg].
\end{align*}

\noindent
\underline{Study of the term $D_{1,1}$.} 
We apply Lemma \ref{lemmeproduit} with $a=4\Vert\tilde{Y}^2 \tilde{X}\Vert_1$ and obtain that
\begin{align}
\label{D1_1}
D_{1,1}&\leq \mathbb{P}\Bigg[a\vert \overline{\X}-\mu_{X}\vert \geq \frac{B\sqrt{\log(n)}}{8K\sqrt{n}}\Bigg]+
\mathbb{P}\Bigg[   \frac{1}{n}\sum_{i=1}^n(Y_i-
\overline{\Y})^2\left\vert 2X_i-\overline{\X}-\mu_{X}\right\vert \geq a   \Bigg].
\end{align}
The first term in (\ref{D1_1}) is bounded by applying Berry-Esseen Inequality (\ref{BE}) and Pollak Inequality (\ref{Pollak}). Since $$\displaystyle 0<K\leq \frac{B}{8a\sigma_X}=\frac{\big\Vert\tilde{X}\varepsilon\Vert_2^2}{32\big\Vert\tilde{Y}^2\tilde{X}\Vert_1\sigma_X},$$ 
then
\begin{align}
\label{D1_1_termeBE}
\mathbb{P}\Bigg[\sqrt{n}\frac{\vert \overline{\X}-\mu_{X}\vert}{\sigma_X} \geq \frac{B\sqrt{\log(n)}}{8Ka\sigma_X}\Bigg]&\leq 
\mathbb{P}\Bigg(\vert \mathcal{N}(0,1)\vert \geq \frac{B\sqrt{\log(n)}}{8Ka\sigma_X}\Bigg)+2c_{\rm{\scriptscriptstyle BE}}\frac{\Vert\tilde{X}\Vert_3^3}{\sigma_X^3\sqrt{n}} \nonumber \\
&\leq \exp\left(-\frac{1}{2}\frac{B^2\log(n)}{64K^2a^2\sigma_X^2}\right)+2c_{\rm{\scriptscriptstyle BE}}\frac{\Vert\tilde{X}\Vert_3^3}{\sigma_X^3\sqrt{n}} \nonumber \\
&\leq \frac{1}{\sqrt{n}}\left(1+2c_{\rm{\scriptscriptstyle BE}}\frac{\Vert\tilde{X}\Vert_3^3}{\sigma_X^3}\right).
\end{align}
For the second term in (\ref{D1_1}), we have to bound up
\begin{multline*}\mathbb{P}\Bigg[   \frac{1}{n}\sum_{i=1}^n(Y_i-
\overline{\Y})^2\left\vert 2X_i-\overline{\X}-\mu_{X}\right\vert \geq a   \Bigg]= \mathbb{P}\Bigg[   \frac{1}{n}\sum_{i=1}^n(Y_i-
\overline{\Y})^2\left\vert 2X_i-\overline{\X}-\mu_{X}\right\vert -\frac{a}{2} \geq \frac{a}{2}  \Bigg]\\
 \leq \frac{2}{a}\left\Vert \frac{1}{n}\sum_{i=1}^n(Y_i-
\overline{\Y})^2\left\vert 2X_i-\overline{\X}-\mu_{X}\right\vert -\frac{a}{2}\right\Vert_1
\\ \qquad   \leq \frac{2}{a} \left\Vert \frac{1}{n}\sum_{i=1}^n\tilde{Y}_i^2\left\vert 2\tilde{X}_i\right\vert -\frac{a}{2}\right\Vert_1+\frac{2}{a}\left\Vert  R_n^{(1,1)}\right\Vert_1,
\end{multline*}
where 
\begin{align*}
\vert R_n^{(1,1)}\vert&=\left\vert  \frac{1}{n}\sum_{i=1}^n(Y_i-
\overline{\Y})^2\left\vert 2X_i-\overline{\X}-\mu_{X}\right\vert- \frac{1}{n}\sum_{i=1}^n\tilde{Y}_i^2\left\vert 2\tilde{X}_i\right\vert       \right\vert\\
& \leq \vert \overline{\X}-\mu_X\vert \frac{1}{n}\sum_{i=1}^n(Y_i-\overline{\Y})^2 + 2\vert \overline{\Y}-\mu_Y\vert^2\frac{1}{n}\sum_{i=1}^n \vert \tilde{X}_i\vert+ 4\vert \overline{\Y}-\mu_Y\vert\frac{1}{n}\sum_{i=1}^n\vert \tilde{X}_i\vert\vert \tilde{Y}_i\vert.
\end{align*}
By applying Cauchy-Schwarz Inequality and Rosenthal Inequality (\ref{rosenthal}) to the middle term we get that  
\begin{multline*}
\left\Vert R_n^{(1,1)}\right\Vert_1\leq \left\Vert\overline{\X}-\mu_X\right\Vert_2\Vert\tilde{Y}\Vert_4^2
+ \frac{16c_{1,\rm{\scriptscriptstyle R}}^2\sigma_X\sigma_Y^2}{n}+\frac{64 c_{2,\rm{\scriptscriptstyle R}}^2\sigma_X}{n^{3/2}}\Vert\tilde{Y}\Vert_4^2 +4\frac{\sigma_Y}{\sqrt{n}}\Vert\tilde{X}\tilde{Y}\Vert_2\nonumber\\
\leq \frac{\sigma_X}{\sqrt{n}}\Vert \tilde{Y}\Vert_4^2
+
\frac{16c_{1,\rm{\scriptscriptstyle R}}^2\sigma_X\sigma_Y^2}{n}+\frac{64 c_{2,\rm{\scriptscriptstyle R}}^2\sigma_X}{n^{3/2}}\Vert\tilde{Y}\Vert_4^2+4\frac{\sigma_Y}{\sqrt{n}}\Vert\tilde{X}\tilde{Y}\Vert_2.
\end{multline*}
Hence
\begin{align}
\label{D11partie2}
\mathbb{P}&\Bigg[   \frac{1}{n}\sum_{i=1}^n(Y_i-
\overline{\Y})^2\left\vert 2X_i-\overline{\X}-\mu_{X}\right\vert  \geq a  \Bigg] \leq \frac{2}{a} \left\Vert \frac{1}{n}\sum_{i=1}^n\tilde{Y}_i^2\left\vert 2\tilde{X}_i\right\vert -\frac{a}{2}\right\Vert_1 
+\frac{2}{a}\left\Vert R_n^{(1,1)}\right\Vert_1\nonumber\\
&\qquad \leq \frac{4}{a\sqrt{n}}\sqrt{\Var[\tilde{Y}^2\vert \tilde{X}\vert]}+\frac{2}{a}\frac{\sigma_X}{\sqrt{n}}\Vert \tilde{Y}\Vert_4^2+
\frac{8\sigma_Y}{a\sqrt{n}}
\Vert\tilde{X}\tilde{Y}\Vert_2   +
\frac{32c_{1,\rm{\scriptscriptstyle R}}^2\sigma_X\sigma_Y^2}{an}+\frac{128 c_{2,\rm{\scriptscriptstyle R}}^2\sigma_X}{an^{3/2}}\Vert\tilde{Y}\Vert_4^2.
\end{align}
Combining (\ref{D1_1_termeBE}) and (\ref{D11partie2}), we get
\begin{align}
\label{majorationD11}
D_{1,1}&\leq
\frac{1}{\sqrt{n}}\left(1+2c_{\rm{\scriptscriptstyle BE}}\frac{\Vert\tilde{X} \Vert_3^3}{\sigma_X^3}\right) +\frac{1}{\sqrt{n}}\frac{1}{2\Vert\tilde{Y}^2 \tilde{X}\Vert_1}\left(2\Vert
\tilde{Y}^2\tilde{X}\Vert_2
+\sigma_X\Vert \tilde{Y}\Vert_4^2+
4\sigma_Y
\Vert \tilde{X} \tilde{Y}\Vert_2\right)
\nonumber \\ & \qquad +
\frac{8c_{1,\rm{\scriptscriptstyle R}}^2}{n}\frac{\sigma_X\sigma_Y^2}{\Vert\tilde{Y}^2 \tilde{X}\Vert_1}+\frac{32c_{2,\rm{\scriptscriptstyle R}}^2}{n^{3/2}}\frac{\sigma_X\Vert\tilde{Y}\Vert_4^2}{\Vert\tilde{Y}^2 \tilde{X}\Vert_1}.
\end{align}
\underline{Study of the term $D_{1,2}.$}
We apply Lemma \ref{lemmeproduit} with
$a=2\sigma_Y^2$, $b=4 \Vert\tilde{X}\Vert_3^3$, $c=2/\sigma_X^2$.
 Then we write
\begin{eqnarray}
\label{D1_2}
D_{1,2}&\leq& \mathbb{P}\Bigg[abc\vert \overline{\X}-\mu_{X}\vert \geq \frac{B\sqrt{\log(n)}}{8K\sqrt{n}}\Bigg]+
\mathbb{P}\Bigg[   \frac{1}{n}\sum_{i=1}^n(Y_i-
\overline{\Y})^2 \geq a   \Bigg] \nonumber \\
& & + \,\mathbb{P}\Bigg[   \frac{1}{n}\sum_{i=1}^n(X_i-
\overline{\X})^2\left\vert 2X_i-\overline{\X}-\mu_{X}\right\vert \geq b   \Bigg]
+\mathbb{P}\Bigg[   \frac{1}{n}\sum_{i=1}^n(X_i-
\overline{\X})^2\leq\frac{1}{c}\Bigg].
\end{eqnarray}
For the first term in (\ref{D1_2}) we apply Berry-Esseen Inequality (\ref{BE}) and Pollak Inequality (\ref{Pollak}). Since
$$\displaystyle 0<K\leq \frac{B}{[8abc\sigma_X]}=
\frac{\Vert\tilde{X}\varepsilon\Vert_2^2
\sigma_X}{128\Vert\tilde{X}\Vert_3^3
\sigma_Y^2},$$ 
then
\begin{eqnarray}
\label{D1_2_termeBE}
\mathbb{P}\Bigg[\sqrt{n}\frac{\vert \overline{\X}-\mu_{X}\vert}{\sigma_X} \geq \frac{B\sqrt{\log(n)}}{8Kabc\sigma_X}\Bigg]&\leq& 
\mathbb{P}\Bigg(\vert \mathcal{N}(0,1)\vert \geq \frac{B\sqrt{\log(n)}}{8Kabc\sigma_X}\Bigg)+\frac{2c_{\rm{\scriptscriptstyle BE}}}{\sqrt{n}}\frac{\Vert\tilde{X} \Vert_3^3}{\sigma_X^3} \nonumber\\
&\leq& \exp\left(-\frac{1}{2}\frac{B^2\log(n)}{64K^2(abc)^2\sigma_X^2}\right)+2c_{\rm{\scriptscriptstyle BE}}\frac{\Vert\tilde{X} \Vert_3^3}{\sigma_X^3\sqrt{n}} \nonumber\\
&\leq& \frac{1}{\sqrt{n}}\left(1+2c_{\rm{\scriptscriptstyle BE}}\frac{\Vert\tilde{X} \Vert_3^3}{\sigma_X^3}\right).
\end{eqnarray}
The second probability and the last one in (\ref{D1_2})  are bounded using Lemma \ref{lemmemajoration_var_Y} and Lemma \ref{lemmemajoration_var_X}. 
Then we bound up the third probability in (\ref{D1_2}). Note first that 
\begin{align*}
\mathbb{P}\Bigg[   \frac{1}{n}\sum_{i=1}^n(X_i-
\overline{\X})^2\left\vert 2X_i-\overline{\X}-\mu_{X}\right\vert \geq b   \Bigg]\leq
\mathbb{P}\Bigg[   \frac{1}{n}\sum_{i=1}^n\tilde{X}_i^2\left\vert 2\tilde{X}_i\right\vert +\vert R_n^{(1,2)}\vert \geq b   \Bigg],
\end{align*}
where 
\begin{align*}\left\vert R_n^{(1,2)}\right\vert=
\left\vert \frac{1}{n}\sum_{i=1}^n(X_i-
\overline{\X})^2\left\vert 2X_i-\overline{\X}-\mu_{X}\right\vert
- \frac{1}{n}\sum_{i=1}^n\tilde{X}_i^2\left\vert 2\tilde{X}_i\right\vert  \right\vert, 
\end{align*}
with
\begin{align*}
\!\!\! \left\vert R_n^{(1,2)}\right\vert
& \leq \vert \overline{\X}-\mu_X\vert \frac{1}{n}\sum_{i=1}^n(X_i-\overline{\X})^2 + 2\vert \overline{\X}-\mu_X\vert^2\frac{1}{n}\sum_{i=1}^n \vert \tilde{X}_i\vert + 4\vert \overline{\X}-\mu_X\vert\frac{1}{n}\sum_{i=1}^n \tilde{X}_i^2.
\end{align*}
Now, we apply Rosenthal Inequality (\ref{rosenthal}) to the second term of the above inequality and we conclude that 
\begin{align*}
\left\Vert R_n^{(1,2)}\right\Vert_1&\leq \left\Vert\overline{\X}-\mu_X\right\Vert_2\Vert\tilde{X}^2\Vert_2+
\frac{16c_{1,\rm{\scriptscriptstyle R}}^2\sigma_X^3}{n}+\frac{64c_{2,\rm{\scriptscriptstyle R}}^2\sigma_X}{n^{3/2}}\Vert\tilde{X}\Vert_4^2 + \frac{4\sigma_X}{\sqrt{n}}\Vert\tilde{X}\Vert_4^2 \nonumber\\
&\leq \frac{5\sigma_X}{\sqrt{n}}\Vert\tilde{X}\Vert_4^2
+\frac{16c_{1,\rm{\scriptscriptstyle R}}^2\sigma_X^3}{n}+\frac{64c_{2,\rm{\scriptscriptstyle R}}^2\sigma_X}{n^{3/2}}\Vert\tilde{X}\Vert_4^2.
\end{align*}
It follows that 
\begin{multline}
\label{D12partie2}
\mathbb{P}\Bigg[   \frac{1}{n}\sum_{i=1}^n(X_i-
\overline{\X})^2\left\vert 2X_i-\overline{\X}-\mu_{X}\right\vert  \geq b  \Bigg]\leq \frac{1}{b}\left\Vert \frac{1}{n}\sum_{i=1}^n\tilde{X}_i^2\left\vert 2\tilde{X}_i\right\vert \right\Vert_1\!\!+\frac{1}{b}\left\Vert R_n^{(1,2)}\right\Vert_1 \\
\leq \frac{2}{b\sqrt{n}}
\Vert\tilde{X}\Vert_6^3
+
\frac{5\sigma_X}{b\sqrt{n}}\Vert\tilde{X}\Vert_4^2
+\frac{16c_{1,\rm{\scriptscriptstyle R}}^2\sigma_X^3}{bn}+\frac{64c_{2,\rm{\scriptscriptstyle R}}^2\sigma_X}{bn^{3/2}}\Vert\tilde{X}\Vert_4^2.
\end{multline}
By combining (\ref{D1_2_termeBE}), (\ref{majoration_var_Y}), (\ref{majoration_var_X}) and 
(\ref{D12partie2}), 
we conclude that
 \begin{multline}\label{majorationD12}
D_{1,2}\leq \frac{1}{\sqrt{n}}\left(1+2c_{\rm{\scriptscriptstyle BE}}\frac{\Vert\tilde{X}\Vert_3^3}{\sigma_X^3} +\frac{2
\Vert\tilde{X}\Vert_4^2
}{\sigma_X^2}+\frac{
\Vert\tilde{Y}\Vert_4^2
}{\sigma_Y^2}
+\frac{\Vert\tilde{X}\Vert_6^3 
}{2\Vert\tilde{X}\Vert_3^3} + \frac{5\sigma_X\Vert\tilde{X}\Vert_4^2}{4\Vert\tilde{X}\Vert_3^3} \right)\\
+\frac{2}{n}\left(1+\frac{2c_{1,\rm{\scriptscriptstyle R}}^2\sigma_X^3}{\Vert\tilde{X}\Vert_3^3}\right)+\frac{16c_{2,\rm{\scriptscriptstyle R}}^2}{n^{3/2}}\frac{\sigma_X\Vert\tilde{X}\Vert_4^2}{\Vert\tilde{X}\Vert_3^3 }.
\end{multline}
\underline{End of the control of the term $D_1$.}
Combining (\ref{majorationD11}) and (\ref{majorationD12}) we obtain 
\begin{align}
D_{1} \leq &
\frac{1}{\sqrt{n}}\frac{1}{2\Vert\tilde{Y}^2 \tilde{X}\Vert_1}\left(2\Vert\tilde{Y}^2\tilde{X}\Vert_2+\sigma_X\Vert \tilde{Y}\Vert_4^2+
4\sigma_Y
\big\Vert\tilde{X}\tilde{Y}\big\Vert_2\right) \nonumber \\ 
&+\frac{2}{\sqrt{n}}\left(1+2c_{\rm{\scriptscriptstyle BE}}\frac{\Vert\tilde{X}\Vert_3^3}{\sigma_X^3} +\frac{\Vert\tilde{X}\Vert_4^2}{\sigma_X^2}+\frac{\Vert\tilde{Y}\Vert_4^2}{2\sigma_Y^2}
+\frac{\Vert\tilde{X}\Vert_6^3}{4\Vert\tilde{X}\Vert_3^3}+\frac{5\sigma_X\Vert\tilde{X}\Vert_4^2}{8\Vert\tilde{X}\Vert_3^3} \right)\nonumber\\
&+\frac{2}{n}\left(1+\frac{2c_{1,\rm{\scriptscriptstyle R}}^2\sigma_X^3}{\Vert\tilde{X}\Vert_3^3}+\frac{4c_{1,\rm{\scriptscriptstyle R}}^2\sigma_X\sigma_Y^2}{\Vert\tilde{Y}^2 \tilde{X}\Vert_1}\right) +\frac{16c_{2,\rm{\scriptscriptstyle R}}^2}{n^{3/2}}\left(\frac{\sigma_X\Vert\tilde{X}\Vert_4^2}{\Vert\tilde{X}\Vert_3^3 } +
\frac{2\sigma_X\Vert\tilde{Y}\Vert_4^2}{\Vert\tilde{Y}^2 \tilde{X}\Vert_1}\right).
\label{majorationD1}
\end{align}
\subsubsection*{Study of the term $D_2$} 
We write
 \begin{align*}
D_2&=\mathbb{P} \Bigg[   \left\vert   \frac{1}{n}\sum_{i=1}^n\tilde{X}_i^2{\widehat{\varepsilon_i}}^2 -\Vert\tilde{X}\varepsilon\Vert_2^2\right\vert \geq   \frac{B\sqrt{\log(n)}}{2K\sqrt{n}}    \Bigg]\leq D_{2,1}+D_{2,2},
\end{align*}
with
\begin{align*}
D_{2,1}=\mathbb{P} \Bigg[   \left\vert   \frac{1}{n}\sum_{i=1}^n\tilde{X}_i^2{\varepsilon_i}^2 -\Vert\tilde{X}\varepsilon\Vert_2^2\right\vert \geq\frac{B\sqrt{\log(n)}}{4K\sqrt{n}}
    \Bigg]
\mbox{ and } 
D_{2,2}=\mathbb{P} \Bigg[   \left\vert   \frac{1}{n}\sum_{i=1}^n\tilde{X}_i^2\left({\widehat{\varepsilon_i}}^2-{\varepsilon_i}^2\right) \right\vert \geq
\frac{B\sqrt{\log(n)}}{4K\sqrt{n}}      \Bigg].
\end{align*}
\underline{Study of the term $D_{2,1}.$}
We first write 
\begin{align*}
D_{2,1}\leq 
\mathbb{P} \left[ \frac{\displaystyle\left\vert \frac{1}{\sqrt{n}}\sum_{i=1}^n\left[\tilde{X}_i^2{\varepsilon_i}^2 -\Vert\tilde{X}\varepsilon\Vert_2^2
\right]\right\vert}{\sqrt{\Var\left(\tilde{X}^2\varepsilon^2\right)}} \geq\frac{B\sqrt{\log(n)}}{4K\sqrt{\Var\left(\tilde{X}^2\varepsilon^2\right)}}\right].
\end{align*}
We now apply Berry-Esseen Inequality (\ref{BE}) and Pollak Inequality (\ref{Pollak}). Since  
$$\displaystyle 0<K\leq \frac{\Vert\tilde{X}\varepsilon\Vert_2^2}{4\Vert\tilde{X}\varepsilon\Vert_4^2}
,$$ 
then
\begin{align}
D_{2,1}
&\leq \mathbb{P}\Bigg(\vert \mathcal{N}(0,1)\vert \geq \frac{B\sqrt{\log(n)}}{4K
\sqrt{\Var\left(\tilde{X}^2\varepsilon^2\right)}
}\Bigg)
+\frac{2c_{\rm{\scriptscriptstyle BE}}}{\sqrt{n}}\frac{\mathbb{E}\left[\left(\tilde{X}^2  \varepsilon^2\right)^3\right]}{\left[\Var\left(\tilde{X}^2\varepsilon^2\right)\right]^{3/2}
} \nonumber\\
&\leq \exp\left(-\frac{1}{2}\frac{B^2\log(n)}{16K^2\big\Vert\tilde{X}\varepsilon\big\Vert_4^4}\right)
+\frac{2c_{\rm{\scriptscriptstyle BE}}}{\sqrt{n}}\frac{\big\Vert\tilde{X}\varepsilon\big\Vert_6^6
}{\left[\Var\left(\tilde{X}^2\varepsilon^2\right)\right]^{3/2}} \leq \frac{1}{\sqrt{n}}\left(1+2c_{\rm{\scriptscriptstyle BE}}\frac{
\big\Vert\tilde{X}\varepsilon\big\Vert_6^6}{\left[\Var\left(\tilde{X}^2\varepsilon^2\right)\right]^{3/2}}\right).
\label{majorationD21}
\end{align}

\noindent
\underline{Study of the term $D_{2,2}.$}
We start by writing ${\widehat{\varepsilon_i}}^2-{\varepsilon_i}^2$ as 
\begin{align*}
{\widehat{\varepsilon_i}}^2-{\varepsilon_i}^2=(\widehat{\varepsilon_i}-\varepsilon_i)(\widehat{\varepsilon_i}+\varepsilon_i)
= \! \Big(\mu_Y-\overline{\Y}+\tilde{X}_i(\tau-\widehat{\tau})+\widehat{\tau}(\overline{\X}-\mu_{X}  )           \Big)
\Big(  2Y_i-\overline{\Y}-\mu_Y-\widehat{\tau}(X_i-\overline{\X})-\tau    \tilde{X}_i          \Big). 
\end{align*}
The term $D_{2,2}$ is then bounded by the sum of the $9$ following  terms:
\begin{align}
\mathbb{P}&\Bigg[ \frac{1}{n}\sum_{i=1}^n\tilde{X}_i^2\vert\mu_Y-\overline{\Y}\vert\vert 2Y_i-\overline{\Y}-\mu_Y\vert\geq\frac{B\sqrt{\log(n)}}{36K\sqrt{n}} \Bigg],\label{I_I}\\
 \mathbb{P}&\Bigg[ \frac{1}{n}\sum_{i=1}^n\tilde{X}_i^2\vert \mu_Y-\overline{\Y}\vert \vert\widehat{\tau}\vert\vert X_i-\overline{\X}\vert  \geq\frac{B\sqrt{\log(n)}}{36K\sqrt{n}}\Bigg],\label{I_II}
\\
 \mathbb{P}&\Bigg[ \frac{1}{n}\sum_{i=1}^n\tilde{X}_i^2\vert\mu_Y-\overline{\Y}\vert \vert \tau  \vert   \vert \tilde{X}_i\vert \geq\frac{B\sqrt{\log(n)}}{36K\sqrt{n}}\Bigg],\label{I_III}
 \\
 \mathbb{P}&\Bigg[  \frac{1}{n}\sum_{i=1}^n\vert\tilde{X}_i\vert^3 \vert\widehat{\tau}-\tau \vert \vert2Y_i-\overline{\Y}-\mu_Y\vert 
\geq\frac{B\sqrt{\log(n)}}{36K\sqrt{n}}\Bigg],\label{II_I} \\
 \mathbb{P}&\Bigg[  \frac{1}{n}\sum_{i=1}^n\vert\tilde{X}_i\vert^3\vert \widehat{\tau}-\tau \vert \vert \widehat{\tau}\vert\vert X_i-\overline{\X}\vert
  \geq\frac{B\sqrt{\log(n)}}{36K\sqrt{n}}\Bigg],\label{II_II}\\
  \mathbb{P}&\Bigg[  \frac{1}{n}\sum_{i=1}^n\vert\tilde{X}_i\vert^3\vert \widehat{\tau}-\tau \vert  \vert \tau \vert\vert   \tilde{X}_i\vert\geq\frac{B\sqrt{\log(n)}}{36K\sqrt{n}}\Bigg],\label{II_III}\\
 \mathbb{P}&\Bigg[  \frac{1}{n}\sum_{i=1}^n\tilde{X}_i^2\vert \widehat{\tau}\vert \vert\overline{\X}-\mu_{X}\vert\vert2Y_i-\overline{\Y}-\mu_Y\vert\geq\frac{B\sqrt{\log(n)}}{36K\sqrt{n}}\Bigg],\label{III_I}\\
   \mathbb{P}&\Bigg[  \frac{1}{n}\sum_{i=1}^n\tilde{X}_i^2 \widehat{\tau}^2 \vert\overline{\X}-\mu_{X}\vert \vert X_i-\overline{\X}\vert  \geq\frac{B\sqrt{\log(n)}}{36K\sqrt{n}}\Bigg],\label{III_II}\\
   \mathbb{P}&\Bigg[  \frac{1}{n}\sum_{i=1}^n\tilde{X}_i^2\vert \widehat{\tau}\vert \vert\overline{\X}-\mu_{X}\vert   \vert \tau  \vert\vert  \tilde{X}_i\vert
     \geq\frac{B\sqrt{\log(n)}}{36K\sqrt{n}}\Bigg]. \label{III_III}
\end{align}
We then have to bound each of these terms. \\

\noindent
\underline{Control of the term (\ref{I_I}) of the upper bound of $D_{2,2}$.}
Let $a=4\big\Vert\tilde{X}^2\tilde{Y}\big\Vert_1$. Applying Lemma \ref{lemmeproduit}, we get 
\begin{multline}
\mathbb{P}\Bigg[ \vert\overline{\Y}-\mu_Y\vert\frac{1}{n}\sum_{i=1}^n\tilde{X}_i^2\vert 2Y_i-\overline{\Y}-\mu_Y\vert \geq\frac{B\sqrt{\log(n)}}{36K\sqrt{n}} \Bigg]\leq  \mathbb{P}\Bigg[ a\vert\overline{\Y}-\mu_Y\vert \geq\frac{B\sqrt{\log(n)}}{36K\sqrt{n}}\Bigg]\\
+\mathbb{P}\Bigg[ \frac{1}{n}\sum_{i=1}^n\tilde{X}_i^2\vert 2Y_i-\overline{\Y}-\mu_Y\vert \geq a \Bigg]. \label{term_I_I}
\end{multline}
The first probability in (\ref{term_I_I}) is bounded using Berry-Esseen Inequality (\ref{BE}) and Pollak Inequality (\ref{Pollak}). Since $$\displaystyle 0<K\leq \frac{B}{36a\sigma_Y}
=\frac{
\big\Vert\tilde{X}\varepsilon\big\Vert_2^2}{144\big\Vert\tilde{X}^2 \tilde{Y}\big\Vert_1\sigma_Y},$$  then
\begin{align*}
\mathbb{P}\Bigg[ \sqrt{n}\frac{\vert\overline{\Y}-\mu_Y\vert}{\sigma_Y}\geq\frac{B\sqrt{\log(n)}}{36Ka\sigma_Y} \Bigg]
&\leq 
\mathbb{P}\Bigg(\vert \mathcal{N}(0,1)\vert \geq \frac{B\sqrt{\log(n)}}{36Ka\sigma_Y}\Bigg)+2c_{\rm{\scriptscriptstyle BE}}\frac{\mathbb{E}\left( \vert \tilde{Y}\vert^3\right)}{\sigma_Y^3\sqrt{n}} \\
&\leq \exp\left(-\frac{1}{2}\frac{B^2\log(n)}{36^2K^2a^2\sigma_Y^2}\right)+2c_{\rm{\scriptscriptstyle BE}}\frac{\big\Vert\tilde{Y}\big\Vert_3^3}{\sigma_Y^3\sqrt{n}} \leq \frac{1}{\sqrt{n}}\left(1+2c_{\rm{\scriptscriptstyle BE}}\frac{\big\Vert\tilde{Y}\big\Vert_3^3
}{\sigma_Y^3}\right).
\end{align*}

\noindent
The second probability in (\ref{term_I_I}) is controlled by applying Markov Inequality  
\begin{align}
\mathbb{P}\Bigg[ \frac{1}{n}\sum_{i=1}^n\tilde{X}_i^2\vert 2Y_i-\overline{\Y}-\mu_Y\vert \geq a \Bigg] &\leq 
\mathbb{P}\Bigg[ \frac{1}{n}\sum_{i=1}^n\tilde{X}_i^2\vert 2Y_i-\overline{\Y}-\mu_Y\vert -\frac{a}{2}\geq \frac{a}{2} \Bigg] \nonumber\\
&\leq \frac{2}{a}\left\Vert\frac{1}{n}\sum_{i=1}^n\tilde{X}_i^2\vert 2Y_i-\overline{\Y}-\mu_Y\vert -\frac{a}{2}\right\Vert_1 \nonumber\\
&\leq \frac{2}{a}\left\Vert\frac{1}{n}\sum_{i=1}^n\tilde{X}_i^22\vert \tilde{Y}_i\vert -\frac{a}{2}\right\Vert_1
+ \frac{2}{a}\left\Vert\vert \overline{\Y}-\mu_Y\vert\frac{1}{n}\sum_{i=1}^n\tilde{X}_i^2\right\Vert_1 \nonumber\\
&\leq \frac{4}{a\sqrt{n}}\big\Vert\tilde{X}^2\tilde{Y}\big\Vert_2
+ \frac{2\sigma_Y}{a\sqrt{n}}\big\Vert \tilde{X}\big\Vert_4^2
\label{majoration_terme2_de_II}.
\end{align}
We then obtain the following bound for the term (\ref{I_I}):  
\begin{multline}
\mathbb{P}\Bigg[ \vert\overline{\Y}-\mu_Y\vert\frac{1}{n}\sum_{i=1}^n\tilde{X}_i^2\vert 2Y_i-\overline{\Y}-\mu_Y\vert \geq\frac{B\sqrt{\log(n)}}{36K\sqrt{n}} \Bigg] \\
\qquad  \leq 
\frac{1}{\sqrt{n}}\left(1+2c_{\rm{\scriptscriptstyle BE}}\frac{\big\Vert \tilde{Y}\big\Vert_3^3}{\sigma_Y^3}+
\frac{\big\Vert\tilde{X}^2\tilde{Y}\big\Vert_2}{\big\Vert\tilde{X}^2 \tilde{Y}\big\Vert_1}
+ \frac{\sigma_Y\big\Vert \tilde{X}\big\Vert_4^2}{2\big\Vert\tilde{X}^2 \tilde{Y}\big\Vert_1}\right).\label{borne_I_I}
\end{multline}
\underline{Control of the term (\ref{I_II}) of the upper bound of $D_{2,2}$.
}
We start by bounding the term (\ref{I_II}) using the bound on $ \widehat{\tau}$ given in (\ref{majorationtauchapeau}). Then, applying Lemma \ref{lemmeproduit} with
$a=2\big\Vert\tilde{X}\big\Vert_3^3$, $b=\sqrt{2}\sigma_Y$ and $c=\sqrt{2}/\sigma_X$, we get 
\begin{align}
\mathbb{P}&\Bigg[ \vert\overline{\Y}-\mu_Y\vert\frac{1}{n}\sum_{i=1}^n{\tilde{X}_i^2\vert X_i-\overline{\X}\vert} \frac{\sqrt{\frac{1}{n}\sum_{i=1}^n{(Y_i-\overline{\Y})^2}}}{\sqrt{\frac{1}{n}\sum_{i=1}^n{(X_i-\overline{\X})^2}}} \geq\frac{B\sqrt{\log(n)}}{36K\sqrt{n}}\Bigg] \nonumber\\
&\leq \mathbb{P}\Bigg[ abc\vert\overline{\Y}-\mu_Y\vert \geq\frac{B\sqrt{\log(n)}}{36K\sqrt{n}}\Bigg] 
+ \mathbb{P}\Bigg[ \frac{1}{n}\sum_{i=1}^n{\tilde{X}_i^2\vert X_i-\overline{\X}\vert}  \geq a\Bigg]
+ \mathbb{P}\Bigg[ \frac{1}{n}\sum_{i=1}^n{(Y_i-\overline{\Y})^2} \geq b^2\Bigg]\nonumber\\
& \qquad + \mathbb{P}\Bigg[ \sqrt{\frac{1}{n}\sum_{i=1}^n{(X_i-\overline{\X})^2}} \leq 1/c\Bigg]. \label{term_I_II}
\end{align}
The first probability in (\ref{term_I_II}) is bounded by applying Berry-Esseen Inequality (\ref{BE}) and Pollak Inequality (\ref{Pollak}). Since  $$\displaystyle 0<K\leq \frac{B}{36abc\sigma_Y}
=
\frac{\big\Vert\tilde{X}\varepsilon\big\Vert_2^2\sigma_X}{144 \big\Vert\tilde{X}\big\Vert_3^3\sigma_Y^2}
,$$
then
\begin{align*}
\mathbb{P}\Bigg[ \sqrt{n}\frac{\vert\overline{\Y}-\mu_Y\vert}{\sigma_Y}\geq \frac{B\sqrt{\log(n)}}{36Kabc\sigma_Y} \Bigg]
&\leq 
\mathbb{P}\Bigg(\vert \mathcal{N}(0,1)\vert \geq \frac{B\sqrt{\log(n)}}{36Kabc\sigma_Y}\Bigg)+2c_{\rm{\scriptscriptstyle BE}}\frac{\mathbb{E}\left( \vert \tilde{Y}\vert^3\right)}{\sigma_Y^3\sqrt{n}} \\
&\leq \exp\left(-\frac{1}{2}\frac{B^2\log(n)}{36^2K^2(abc)^2\sigma_Y^2}\right)+2c_{\rm{\scriptscriptstyle BE}}\frac{
\big\Vert\tilde{Y}\big\Vert_3^3}{\sigma_Y^3\sqrt{n}} \leq \frac{1}{\sqrt{n}}\left(1+2c_{\rm{\scriptscriptstyle BE}}\frac{\big\Vert\tilde{Y}\big\Vert_3^3}{\sigma_Y^3}\right).
\end{align*}

\noindent
The second probability in (\ref{term_I_II}) is controlled by applying Markov Inequality: 
\begin{align}
\mathbb{P}\Bigg[ \frac{1}{n}\sum_{i=1}^n\tilde{X}_i^2\vert X_i-\overline{\X}\vert \geq a \Bigg] &\leq 
\mathbb{P}\Bigg[ \frac{1}{n}\sum_{i=1}^n\tilde{X}_i^2\vert X_i-\overline{\X}\vert -\frac{a}{2}\geq \frac{a}{2} \Bigg]  \leq \frac{2}{a}\left\Vert\frac{1}{n}\sum_{i=1}^n\tilde{X}_i^2\vert X_i-\overline{\X}\vert -\frac{a}{2}\right\Vert_1 \nonumber\\
& \leq \frac{2}{a}\left\Vert\frac{1}{n}\sum_{i=1}^n\tilde{X}_i^2\vert \tilde{X}_i\vert -\frac{a}{2}\right\Vert_1
+ \frac{2}{a}\left\Vert\vert \overline{\X}-\mu_X\vert\frac{1}{n}\sum_{i=1}^n\tilde{X}_i^2\right\Vert_1 \nonumber\\
&\leq \frac{2}{a\sqrt{n}}\big\Vert\tilde{X}^3\big\Vert_2
+ \frac{2\sigma_X}{a\sqrt{n}}\big\Vert \tilde{X}\big\Vert_4^2
\label{majorationterme1-2_Markov}.
\end{align}
The two last probabilities in (\ref{term_I_II}) are bounded by (\ref{majoration_var_Y}) and (\ref{majoration_var_X}). 
Hence, we obtain the following bound for the term (\ref{I_II}) 
\begin{multline}
\mathbb{P}\Bigg[\vert \overline{\Y}-\mu_Y\vert \vert\widehat{\tau}\vert\frac{1}{n}\sum_{i=1}^n\tilde{X}_i^2 \vert X_i-\overline{\X}\vert \geq\frac{B\sqrt{\log(n)}}{36K\sqrt{n}}\Bigg] \\
\leq \frac{1}{\sqrt{n}}\left(1+2c_{\rm{\scriptscriptstyle BE}}\frac{\big\Vert \tilde{Y}\big\Vert_3^3}{\sigma_Y^3}+ \frac{\big\Vert \tilde{X}\big\Vert_6^3}{\big\Vert \tilde{X}\big\Vert_3^3}
+ \frac{\sigma_X \big\Vert \tilde{X}\big\Vert_4^2}{\big\Vert \tilde{X}\big\Vert_3^3}
+ \frac{\big\Vert\tilde{Y}\big\Vert_4^2}{\sigma_Y^2}
+ \frac{2\big\Vert\tilde{X}\big\Vert_4^2}{\sigma_X^2}\right)+\frac{2}{n}.\label{borne_I_II}
\end{multline}
\underline{Control of the term (\ref{I_III}) of the upper bound of $D_{2,2}$.
}
Applying Lemma \ref{lemmeproduit} with $a=2\vert\tau \vert
\Vert \tilde{X}\Vert^3_3$, we get 
\begin{multline} \label{term_I_III}
\mathbb{P}\Bigg[ \vert\overline{\Y}-\mu_Y\vert\vert\tau \vert\frac{1}{n}\sum_{i=1}^n{\vert \tilde{X}_i\vert^3} \geq\frac{B\sqrt{\log(n)}}{36K\sqrt{n}} \Bigg]
 \leq \mathbb{P}\Bigg[ a\vert\overline{\Y}-\mu_Y\vert \geq\frac{B\sqrt{\log(n)}}{36K\sqrt{n}}\Bigg]+\mathbb{P}\Bigg[ \vert\tau \vert\frac{1}{n}\sum_{i=1}^n{\vert \tilde{X}_i\vert^3} \geq a \Bigg].
\end{multline}
The first probability in (\ref{term_I_III}) is bounded by applying Berry-Esseen Inequality (\ref{BE}) and Pollak inequality (\ref{Pollak}). Since $$\displaystyle 0<K\leq \frac{B}{36a\sigma_Y}
=
\frac{\big\Vert\tilde{X}\varepsilon\big\Vert_2^2}{72 \vert\tau \vert\big\Vert\tilde{X}\big\Vert_3^3\sigma_Y},
$$
then
\begin{align*}
\mathbb{P}\Bigg[ \sqrt{n}\frac{\vert\overline{\Y}-\mu_Y\vert}{\sigma_Y}\geq\frac{B\sqrt{\log(n)}}{36Ka\sigma_Y} \Bigg]
&\leq 
\mathbb{P}\Bigg(\vert \mathcal{N}(0,1)\vert \geq \frac{B\sqrt{\log(n)}}{36Ka\sigma_Y}\Bigg)+2c_{\rm{\scriptscriptstyle BE}}\frac{\mathbb{E}\left( \vert \tilde{Y}\vert^3\right)}{\sigma_Y^3\sqrt{n}} \\
&\leq \exp\left(-\frac{1}{2}\frac{B^2\log(n)}{36^2K^2a^2\sigma_Y^2}\right)+2c_{\rm{\scriptscriptstyle BE}}\frac{\big\Vert\tilde{Y}\big\Vert_3^3}{\sigma_Y^3\sqrt{n}} \leq \frac{1}{\sqrt{n}}\left(1+2c_{\rm{\scriptscriptstyle BE}}\frac{\big\Vert\tilde{Y}\big\Vert_3^3}{\sigma_Y^3}\right).
\end{align*}

\noindent
The second probability in (\ref{term_I_III}) is bounded by applying Markov Inequality:
\begin{align}
\mathbb{P}\Bigg[  \vert\tau \vert \frac{1}{n}\sum_{i=1}^n\vert \tilde{X_i}\vert^3 \geq a   \Bigg]=
\mathbb{P}\Bigg[   \frac{1}{n}\sum_{i=1}^n\vert \tilde{X_i}\vert^3-\frac{a}{2\vert\tau \vert}\geq \frac{a}{2\vert\tau \vert}   \Bigg]\leq 
\frac{2\vert\tau \vert}{a}\left\Vert \frac{1}{n}\sum_{i=1}^n
\vert \tilde{X}_i\vert^3-\frac{a}{2\vert\tau \vert}\right\Vert_1 \leq \frac{2\vert\tau \vert 
\big\Vert \tilde{X}  \big\Vert_6^3
}{a\sqrt{n}}
\label{majorationterme1-3_Markov}.
\end{align}
Consequently, we get the following bound for the term (\ref{I_III}): 
\begin{align}
\mathbb{P}\Bigg[ \vert\overline{\Y}-\mu_Y\vert\vert\tau \vert\frac{1}{n}\sum_{i=1}^n{\vert \tilde{X}_i\vert^3} \geq\frac{B\sqrt{\log(n)}}{36K\sqrt{n}} \Bigg]&\leq 
\frac{1}{\sqrt{n}}\left(1+2c_{\rm{\scriptscriptstyle BE}}\frac{\big\Vert \tilde{Y}   \big\Vert_3^3}{\sigma_Y^3}+ \frac{ 
\big\Vert \tilde{X}   \big\Vert_6^3
}{\big\Vert \tilde{X}   \big\Vert_3^3}\right) .\label{borne_I_III}
\end{align}
\underline{Control of the term (\ref{II_I}) of the upper bound of $D_{2,2}$
}
Applying Lemma \ref{lemmeproduit} with $a=8\Vert \tilde{X}^3 \tilde{Y} \Vert_1
$ and $b=2/{\sigma_X^2}$, we get
\begin{align}
\mathbb{P}\Bigg[ \vert\widehat{\tau}-\tau \vert\frac{1}{n}\sum_{i=1}^n\vert \tilde{X}_i\vert^3&\vert 2Y_i-\overline{\Y}-\mu_Y\vert \geq \frac{B\sqrt{\log(n)}}{36K\sqrt{n}} \Bigg] \nonumber\\
\leq &\, \mathbb{P}\Bigg[ ab\left\vert\frac{1}{n}\sum_{i=1}^n{(X_i-
\overline{\X})\varepsilon_i}\right\vert  \geq\frac{B\sqrt{\log(n)}}{36K\sqrt{n}}\Bigg] \nonumber
\\
&+\mathbb{P}\Bigg[ \frac{1}{n}\sum_{i=1}^n{\vert \tilde{X}_i\vert^3\vert 2Y_i-\overline{\Y}-\mu_Y\vert} \geq a \Bigg] 
+\mathbb{P}\Bigg[ \frac{1}{n}\sum_{i=1}^n{(X_i-\overline{\X})^2}\leq 1/b \Bigg] \nonumber\\
&\quad\leq \mathbb{P}\Bigg[ ab\left\vert\frac{1}{n}\sum_{i=1}^n{\tilde{X_i}\varepsilon_i}\right\vert  \geq\frac{B\sqrt{\log(n)}}{72K\sqrt{n}}\Bigg]
+\mathbb{P}\Bigg[ ab\vert\overline{\X}-\mu_X\vert\left\vert\frac{1}{n}\sum_{i=1}^n{\varepsilon_i}\right\vert  \geq\frac{B\sqrt{\log(n)}}{72K\sqrt{n}}\Bigg] \nonumber\\
&\qquad+\mathbb{P}\Bigg[ \frac{1}{n}\sum_{i=1}^n{\vert \tilde{X}_i\vert^3\vert 2Y_i-\overline{\Y}-\mu_Y\vert} \geq a \Bigg] 
+\mathbb{P}\Bigg[ \frac{1}{n}\sum_{i=1}^n{(X_i-\overline{\X})^2}\leq 1/b \Bigg]. \label{term_II_I}
\end{align}
The first probability in (\ref{term_II_I}) is bounded by applying Berry-Esseen Inequality (\ref{BE}) and Pollak inequality (\ref{Pollak}). Since 
$$\displaystyle 0<K\leq 
\frac{B}{72ab\big\Vert\tilde{X}\varepsilon\big\Vert_2}=
\frac{\sigma_X^2\big\Vert\tilde{X}\varepsilon\big\Vert_2}{1152\big\Vert \tilde{X}^3\tilde{Y}\big\Vert_1},
$$
then
\begin{align*}
\mathbb{P}\Bigg[ \sqrt{n}\frac{\vert 
n^{-1}\sum_{i=1}^n{\tilde{X}_i\varepsilon_i}\vert}{\big\Vert\tilde{X}\varepsilon\big\Vert_2
}\geq\frac{B\sqrt{\log(n)}}{72Kab 
\big\Vert\tilde{X}\varepsilon\big\Vert_2} \Bigg]
&\leq 
\mathbb{P}\Bigg(\vert \mathcal{N}(0,1)\vert \geq \frac{B\sqrt{\log(n)}}{72Kab
\big\Vert\tilde{X}\varepsilon\big\Vert_2
}\Bigg)+2c_{\rm{\scriptscriptstyle BE}}\frac{
\big\Vert \tilde{X}\varepsilon \big\Vert_3^3
}{\big
\Vert \tilde{X}\varepsilon\big\Vert_2^3\sqrt{n}} \\
&\leq \exp\left(-\frac{1}{2}\frac{B^2\log(n)}{72^2(Kab)^2
\big\Vert\tilde{X}\varepsilon\big\Vert_2^2}\right)+2c_{\rm{\scriptscriptstyle BE}}\frac{
\big\Vert \tilde{X}\varepsilon \big\Vert_3^3
}{\big\Vert \tilde{X}\varepsilon\big\Vert_2^3\sqrt{n}} \\
&\leq \frac{1}{\sqrt{n}}\left(1+2c_{\rm{\scriptscriptstyle BE}}
\frac{\big\Vert \tilde{X}\varepsilon \big\Vert_3^3}{
\big\Vert \tilde{X}\varepsilon\big\Vert_2^{3}}\right).
\end{align*}

\noindent
The second probability in (\ref{term_II_I}) is bounded by applying Markov Inequality: 
\begin{align}
\label{majoration_avec_epsilonbar}
\mathbb{P}\Bigg[ ab\left\vert
\overline{\X}-\mu_X\right\vert\left\vert\frac{1}{n}\sum_{i=1}^n{\varepsilon_i}\right\vert  \geq\frac{B\sqrt{\log(n)}}{72 K\sqrt{n}}\Bigg]
&\leq \frac{72 Kab\sqrt{n}}{B\sqrt{\log(n)}}\left\Vert
\overline{\X}-\mu_X\right\Vert_2\left\Vert\frac{1}{n}\sum_{i=1}^n{\varepsilon_i}\right\Vert_2 \leq \frac{72 K a b \sigma_X\sigma_{\varepsilon}}{B\sqrt{n\log(n)}}.
\end{align}
For the third probability in (\ref{term_II_I}), we write   
$$
\mathbb{P}\Bigg[ 
\frac{1}{n}\sum_{i=1}^n\vert \tilde{X}_i\vert^3\vert 2Y_i-\overline{\Y}-\mu_Y\vert \geq a \Bigg] \leq 
\mathbb{P}\Bigg[ 
\frac{2}{n}\sum_{i=1}^n \vert \tilde{X}_i\vert^3 \vert\tilde{Y}_i\vert  \geq \frac a 2 \Bigg]
+
\mathbb{P}\Bigg[ \vert \overline{\Y}-\mu_Y\vert
\frac{1}{n}\sum_{i=1}^n\vert \tilde{X}_i\vert^3 \geq \frac a 2 \Bigg].
$$
The first term on right hand side is bounded using von Bahr-Esseen Inequality (\ref{Von_Bahr_Esseen}):
\begin{align*}
\mathbb{P}\Bigg[ 
\frac{2}{n}\sum_{i=1}^n \vert \tilde{X}_i\vert^3 \vert\tilde{Y}_i\vert  \geq \frac a 2 \Bigg] &\leq 
\mathbb{P}\Bigg[ \frac{1}{n}\sum_{i=1}^n \left ( \vert\tilde{X}_i\vert^3 \vert\tilde{Y}_i\vert -{\mathbb E}(\vert\tilde{X}\vert^3 \vert\tilde{Y}\vert)\right)\geq \frac{a}{8} \Bigg]\\&\leq 
\left (\frac{8}{a}\right )^{3/2}\left\Vert \frac{1}{n}\sum_{i=1}^n \left ( \vert\tilde{X}_i\vert^3 \vert\tilde{Y}_i\vert -{\mathbb E}(\vert\tilde{X}\vert^3 \vert\tilde{Y}\vert) \right )\right\Vert_{3/2}^{3/2}\nonumber\\
&\leq \frac{\sqrt 2}{\sqrt n} \left (\frac{8}{a}\right )^{3/2}  \left\Vert  \vert\tilde{X}\vert^3 \vert\tilde{Y}\vert -{\mathbb E}(\vert\tilde{X}\vert^3 \vert\tilde{Y}\vert) \right\Vert_{3/2}^{3/2} \leq  \frac{\sqrt{2}\,2^{6}}{a^{3/2}\sqrt n}  \left\Vert \vert\tilde{X}\vert^3 \vert\tilde{Y}\vert\right\Vert_{3/2}^{3/2}.
\end{align*}
The second term on right hand side is bounded using Markov Inequality: 
$$
\mathbb{P}\Bigg[ \vert \overline{\Y}-\mu_Y\vert
\frac{1}{n}\sum_{i=1}^n\vert \tilde{X}_i\vert^3 \geq \frac a 2 \Bigg]\leq \frac{2}{a}\left\Vert\vert \overline{\Y}-\mu_Y\vert\frac{1}{n}\sum_{i=1}^n\vert \tilde{X}_i\vert^3\right\Vert_1 \leq \frac{2}{a}\left\Vert \overline{\Y}-\mu_Y\right\Vert_2
\big\Vert \tilde{X}^3\big\Vert_2 \leq \frac{2\sigma_Y\big\Vert \tilde{X}\big\Vert_6^3}{a\sqrt{n}}.
$$
Then we obtain the following bound for the third probability in (\ref{term_II_I}): 
\begin{eqnarray}
\label{majoration_Markov_proba3_de_termeII_I}
\mathbb{P}\Bigg[ 
\frac{1}{n}\sum_{i=1}^n\vert \tilde{X}_i\vert^3\vert 2Y_i-\overline{\Y}-\mu_Y\vert \geq a \Bigg] \leq \frac{\sqrt{2}\,2^{6}}{a^{3/2}\sqrt n}  \left\Vert \vert\tilde{X}\vert^3 \vert\tilde{Y}\vert\right\Vert_{3/2}^{3/2}+\frac{2\sigma_Y\big\Vert \tilde{X}\big\Vert_6^3}{a\sqrt{n}}.
\end{eqnarray}
The last probability in (\ref{term_II_I}) is bounded by (\ref{majoration_var_X}). Then we get the following bound for the term (\ref{II_I}):
\begin{align}
\mathbb{P}&\Bigg[ \vert\widehat{\tau}-\tau \vert\frac{1}{n}\sum_{i=1}^n\vert \tilde{X}_i\vert^3\vert 2Y_i-\overline{\Y}-\mu_Y\vert \geq\frac{B\sqrt{\log(n)}}{36K\sqrt{n}} \Bigg] \nonumber\\
&\leq \frac{1}{\sqrt{n}}\left(1+2c_{\rm{\scriptscriptstyle BE}}
\frac{
\big\Vert\tilde{X} \varepsilon\big\Vert_3^3}{\big\Vert \tilde{X}\varepsilon\big\Vert_2^3
}+ \frac{4\big\Vert\tilde{X}^3\tilde{Y}\big\Vert_{3/2}^{3/2}}{\big\Vert \tilde{X}^3 \tilde{Y}\big\Vert_1^{3/2}}+ \frac{\sigma_Y\big\Vert \tilde{X}\big\Vert_6^3}{4\big\Vert \tilde{X}^3 \tilde{Y}\big\Vert_1}+\frac{2\big\Vert \tilde{X}\big\Vert_4^2}{\sigma_X^2}\right)
+ \frac{1152 K \big\Vert\tilde{X}^3 \tilde{Y}\big\Vert_1
\sigma_{\varepsilon}}{
\big\Vert\tilde{X}\varepsilon\big\Vert_2^2
\sigma_X\sqrt{n\log(n)}}+\frac{2}{n}. \label{borne_II_I}
\end{align}
\underline{Control of the term (\ref{II_II}) of the upper bound of $D_{2,2}$.}
We start by bounding the term (\ref{II_II}) using the bound on $ \widehat{\tau}$ given in (\ref{majorationtauchapeau}). Then, applying Lemma \ref{lemmeproduit} with $a=(\sqrt{2}/\sigma_X)^3$, $b=\sqrt{2}\sigma_Y$ and $c=4
\big\Vert \tilde{X}\big\Vert_4^4$, we get 
\begin{align}
\mathbb{P}&\Bigg[ \left\vert
\frac{1}{n}\sum_{i=1}^n{(X_i-
\overline{\X})\varepsilon_i}
\right\vert
\frac{\sqrt{\frac{1}{n}\sum_{i=1}^n{(Y_i-
\overline{\Y})^2}}}{\left[\frac{1}{n}\sum_{i=1}^n{(X_i-
\overline{\X})^2}\right]^{3/2}}
\frac{1}{n}\sum_{i=1}^n\vert \tilde{X}_i\vert^3\vert X_i-\overline{\X}\vert \geq\frac{B\sqrt{\log(n)}}{36K\sqrt{n}} \Bigg] \nonumber\\
&\leq \mathbb{P}\Bigg[ abc\left\vert\frac{1}{n}\sum_{i=1}^n{\tilde{X}_i\varepsilon_i}\right\vert  \geq\frac{B\sqrt{\log(n)}}{72K\sqrt{n}}\Bigg] + \mathbb{P}\Bigg[ abc\vert\overline{\X}-\mu_X\vert\left\vert\frac{1}{n}\sum_{i=1}^n{\varepsilon_i}\right\vert  \geq\frac{B\sqrt{\log(n)}}{72K\sqrt{n}}\Bigg] \nonumber \\
&+\mathbb{P}\Bigg[ \frac{1}{n}\sum_{i=1}^n{(X_i-
\overline{\X})^2}\leq \frac{\sigma_X^2}{2} \Bigg]+\mathbb{P}\Bigg[ \frac{1}{n}\sum_{i=1}^n{(Y_i-
\overline{\Y})^2}\geq 2\sigma_Y^2 \Bigg]+ 
\mathbb{P}\Bigg[ 
\frac{1}{n}\sum_{i=1}^n\vert \tilde{X}_i\vert^3\vert X_i-\overline{\X}\vert \geq c \Bigg]. \label{term_II_II}
\end{align}

\noindent
The first probability in (\ref{term_II_II}) is bounded by applying Berry-Esseen Inequality (\ref{BE}) and Pollak inequality (\ref{Pollak}). Since  
$$\displaystyle 0<K\leq \frac{B}{72abc\big\Vert \tilde{X}\varepsilon \big\Vert_2}=\frac{\sigma_X^3\big\Vert \tilde{X}\varepsilon \big\Vert_2}{1152 \big\Vert \tilde{X}\big\Vert_4^4\sigma_Y},$$ we get 
\begin{align*}
\mathbb{P}\left[ \sqrt{n}\frac{\displaystyle\left\vert 
\frac{1}{n}\sum_{i=1}^n{\tilde{X}_i\varepsilon_i}\right\vert}{\big\Vert \tilde{X}\varepsilon \big\Vert_2}\geq\frac{B\sqrt{\log(n)}}{72Kabc\big\Vert \tilde{X}\varepsilon \big\Vert_2} \right]
&\leq \frac{1}{\sqrt{n}}\left(1+2c_{\rm{\scriptscriptstyle BE}}
\frac{\big\Vert \tilde{X}\varepsilon \big\Vert_3^3}{ \big\Vert \tilde{X}\varepsilon \big\Vert_2^3}\right).
\end{align*}
We bound the second probability in (\ref{term_II_II}) as in   (\ref{majoration_avec_epsilonbar}): 
\begin{align*}
\mathbb{P}\Bigg[ abc\left\vert
\overline{\X}-\mu_X\right\vert\left\vert \frac{1}{n}\sum_{i=1}^n{\varepsilon_i}\right\vert  \geq\frac{B\sqrt{\log(n)}}{72 K\sqrt{n}}\Bigg]
&\leq \frac{72 Kabc\sqrt{n}}{B\sqrt{\log(n)}}\left\Vert
\overline{\X}-\mu_X\right\Vert_2\left\Vert\frac{1}{n}\sum_{i=1}^n{\varepsilon_i}\right\Vert_2 
\leq \frac{72 K a b c\sigma_X\sigma_{\varepsilon}}{B\sqrt{n\log(n)}}.
\end{align*}
The third and fourth probabilities in (\ref{term_II_II}) are bounded by (\ref{majoration_var_X}) and (\ref{majoration_var_Y}) and the last one is bounded as follows: 
$$
\mathbb{P}\Bigg[ 
\frac{1}{n}\sum_{i=1}^n\vert \tilde{X}_i\vert^3\vert X_i-\overline{\X}\vert \geq c \Bigg] \leq 
\mathbb{P}\Bigg[ 
\frac{1}{n}\sum_{i=1}^n \tilde{X}_i^4 \geq \frac c 2 \Bigg]
+
\mathbb{P}\Bigg[ \vert \overline{\X}-\mu_X\vert
\frac{1}{n}\sum_{i=1}^n\vert \tilde{X}_i\vert^3 \geq \frac c 2 \Bigg].
$$
The first term on right hand is bounded by von Bahr-Esseen Inequality (\ref{Von_Bahr_Esseen}):
\begin{align*}
\mathbb{P}\Bigg[ \frac{1}{n}\sum_{i=1}^n \tilde{X}_i^4 \geq \frac c 2 \Bigg] &\leq 
\mathbb{P}\Bigg[ \frac{1}{n}\sum_{i=1}^n \left ( \tilde{X}_i^4 -{\mathbb E}(\tilde{X}^4)\right)\geq \frac{c}{4} \Bigg] \leq 
\left (\frac{4}{c}\right )^{3/2}\left\Vert \frac{1}{n}\sum_{i=1}^n \left (\tilde{X}_i^4 -{\mathbb E}(\tilde{X}^4) \right )\right\Vert_{3/2}^{3/2}\nonumber\\
&\leq \frac{\sqrt 2}{\sqrt n} \left (\frac{4}{c}\right )^{3/2}  \left\Vert  \tilde{X}^4 -{\mathbb E}(\tilde{X}^4) \right\Vert_{3/2}^{3/2} \leq  \frac{2^{5}}{c^{3/2}\sqrt n}  \big\Vert \tilde{X}\big\Vert_6^6.
\end{align*}
The second term on right hand is bounded by Markov Inequality: 
$$
\!\!\!\mathbb{P}\Bigg[ \vert \overline{\X}-\mu_X\vert
\frac{1}{n}\sum_{i=1}^n\vert \tilde{X}_i\vert^3 \geq \frac c 2 \Bigg]\leq \frac{2}{c}\left\Vert\vert \overline{\X}-\mu_X\vert\frac{1}{n}\sum_{i=1}^n\vert \tilde{X}_i\vert^3\right\Vert_1 \!\!\!\leq \frac{2}{c}\left\Vert \overline{\X}-\mu_X\right\Vert_2
\big\Vert \tilde{X}^3\big\Vert_2 \leq \frac{2\sigma_X\big\Vert \tilde{X}\big\Vert_6^3}{c\sqrt{n}}.
$$
Then we have  
\begin{align}
\mathbb{P}\Bigg[ 
\frac{1}{n}\sum_{i=1}^n\vert \tilde{X}_i\vert^3\vert X_i-\overline{\X}\vert \geq c \Bigg] 
&\leq \frac{1}{\sqrt n} \left(\frac{4\big\Vert \tilde{X}\big\Vert_6^6}{\big\Vert \tilde{X}\big\Vert_4^6}+\frac{\sigma_X\big\Vert \tilde{X}\big\Vert_6^3}{2\big\Vert \tilde{X}\big\Vert_4^4}\right).
\label{majoration_terme5deII_II}
\end{align}
Hence we get the following bound for the term (\ref{II_II}):
\begin{multline}
\mathbb{P}\Bigg[ \left\vert \widehat{\tau}(\widehat{\tau}-\tau ) \right\vert\frac{1}{n}\sum_{i=1}^n\vert \tilde{X}_i\vert ^3\vert X_i-\overline{\X}\vert 
 \geq\frac{B\sqrt{\log(n)}}{36K\sqrt{n}}\Bigg] 
 \\
\leq \frac{1}{\sqrt{n}}\left(1+2c_{\rm{\scriptscriptstyle BE}}
\frac{
\big\Vert \tilde{X}\varepsilon\big\Vert_3^3}{
\big\Vert \tilde{X}\varepsilon\big\Vert_2^3}
+\frac{\big\Vert \tilde{Y}\big\Vert_4^2}{\sigma_Y^2}+\frac{2\big\Vert \tilde{X}\big\Vert_4^2}{\sigma_X^2} +\frac{4\big\Vert \tilde{X}\big\Vert_6^6}{\big\Vert \tilde{X}\big\Vert_4^6}
+ \frac{\sigma_X\big\Vert \tilde{X}\big\Vert_6^3}{2\big\Vert \tilde{X}\big\Vert_4^4} \right) 
\\ +\frac{1152 K \sigma_Y \sigma_{\varepsilon}\big\Vert \tilde{X}\big\Vert_4^4}{\big\Vert \tilde{X}\varepsilon\big\Vert_2^2\sigma_X^2\sqrt{n\log(n)}}+\frac{2}{n}. \label{borne_II_II}
\end{multline}
\underline{Control of the term (\ref{II_III}) of the upper bound of $D_{2,2}$.}
Applying Lemma \ref{lemmeproduit} with $a=2\vert\tau\vert \big\Vert \tilde{X} \big\Vert_4^4$ and $b=2/{\sigma_X^2}$, we get 
\begin{multline}
\mathbb{P}\Bigg[ \vert\widehat{\tau}-\tau\vert \vert\tau\vert \frac{1}{n}\sum_{i=1}^n\tilde{X}_i^4\geq\frac{B\sqrt{\log(n)}}{36K\sqrt{n}}\Bigg] \\
\leq \mathbb{P}\Bigg[ ab\left\vert\frac{1}{n}\sum_{i=1}^n{(X_i-
\overline{\X})\varepsilon_i}\right\vert  \geq\frac{B\sqrt{\log(n)}}{36K\sqrt{n}}\Bigg]
+\mathbb{P}\Bigg[ \vert\tau\vert \frac{1}{n}\sum_{i=1}^n{\tilde{X}_i^4} \geq a \Bigg] +\mathbb{P}\Bigg[ \frac{1}{n}\sum_{i=1}^n{(X_i-\overline{\X})^2}\leq 1/b \Bigg] \\
\leq \mathbb{P}\Bigg[ ab\left\vert\frac{1}{n}\sum_{i=1}^n{\tilde{X}_i\varepsilon_i}\right\vert  \geq\frac{B\sqrt{\log(n)}}{72K\sqrt{n}}\Bigg] 
+\mathbb{P}\Bigg[ ab\vert\overline{\X}-\mu_X\vert\left\vert\frac{1}{n}\sum_{i=1}^n{\varepsilon_i}\right\vert  \geq\frac{B\sqrt{\log(n)}}{72K\sqrt{n}}\Bigg] \\
+\mathbb{P}\Bigg[ \vert\tau\vert \frac{1}{n}\sum_{i=1}^n{\tilde{X}_i^4} \geq a \Bigg] +\mathbb{P}\Bigg[ \frac{1}{n}\sum_{i=1}^n{(X_i-\overline{\X})^2}\leq 1/b \Bigg]. \label{term_II_III}
\end{multline}
The first probability in (\ref{term_II_III}) is bounded by applying Berry-Esseen Inequality (\ref{BE}) and Pollak Inequality (\ref{Pollak}). Since 
$$\displaystyle 0<K\leq \frac{B}{72ab\big\Vert \tilde{X}\varepsilon\big\Vert_2}=\frac{\sigma_X^2\big\Vert \tilde{X}\varepsilon\big\Vert_2}{288\vert\tau\vert \big\Vert \tilde{X}\big\Vert_4^4}
,$$
we get 
\begin{multline*}
\mathbb{P}\Bigg[ \sqrt{n}\frac{\vert 
n^{-1}\sum_{i=1}^n{\tilde{X}_i\varepsilon_i}\vert}{
\big\Vert\tilde{X}\varepsilon \big\Vert_2}\geq\frac{B\sqrt{\log(n)}}{72Kab
\big\Vert\tilde{X}\varepsilon \big\Vert_2
} \Bigg]\leq \exp\left(-\frac{1}{2}\frac{B^2\log(n)}{72^2(Kab)^2
\big\Vert \tilde{X}\varepsilon \big\Vert_2^2}\right) +2c_{\rm{\scriptscriptstyle BE}}\frac{
\big\Vert \tilde{X}\varepsilon\big\Vert_3^3
}{
\big\Vert \tilde{X}\varepsilon\big\Vert_2^{3}\sqrt{n}}\\
\leq \frac{1}{\sqrt{n}}\left(1+2c_{\rm{\scriptscriptstyle BE}}
\frac{\big\Vert\tilde{X}\varepsilon\big\Vert_3^3}{
\big\Vert \tilde{X}\varepsilon\big\Vert_2^{3}}\right). 
\end{multline*}
The second probability in (\ref{term_II_III}) is bounded by applying Markov Inequality: 
\begin{align*}
\mathbb{P}\Bigg[ ab\left\vert
\overline{\X}-\mu_X\right\vert\left\vert\frac{1}{n}\sum_{i=1}^n{\varepsilon_i}\right\vert  \geq\frac{B\sqrt{\log(n)}}{72 K\sqrt{n}}\Bigg]
&\leq \frac{72 Kab\sqrt{n}}{B\sqrt{\log(n)}}\left\Vert
\overline{\X}-\mu_X\right\Vert_2\left\Vert\frac{1}{n}\sum_{i=1}^n{\varepsilon_i}\right\Vert_2 \leq \frac{72 K a b \sigma_X\sigma_{\varepsilon}}{B\sqrt{n\log(n)}}.
\end{align*}
The third probability in (\ref{term_II_III}) is bounded by applying Markov Inequality and von Bahr-Esseen Inequality (\ref{Von_Bahr_Esseen}):
\begin{multline}
\mathbb{P}\Bigg[ \vert\tau\vert\frac{1}{n}\sum_{i=1}^n \tilde{X}_i^4 \geq a \Bigg] \leq 
\mathbb{P}\Bigg[ \vert\tau\vert\frac{1}{n}\sum_{i=1}^n \left ( \tilde{X}_i^4 -{\mathbb E}(\tilde{X}^4)\right)\geq \frac{a}{2} \Bigg] \\
\leq 
\left (\frac{2\vert\tau\vert}{a}\right )^{3/2}\left\Vert \frac{1}{n}\sum_{i=1}^n \left (\tilde{X}_i^4 -{\mathbb E}(\tilde{X}^4) \right )\right\Vert_{3/2}^{3/2}\\
\leq \frac{4}{\sqrt n} \left (\frac{\vert\tau\vert}{a}\right )^{3/2} \left\Vert  \tilde{X}^4 -{\mathbb E}(\tilde{X}^4) \right\Vert_{3/2}^{3/2} \leq \frac{8\sqrt 2}{\sqrt n} \left (\frac{\vert\tau\vert}{a}\right )^{3/2} \big\Vert \tilde{X}\big\Vert_6^6.
\label{majoration_terme3deII_III}
\end{multline}
The last probability in (\ref{term_II_III}) is bounded by  (\ref{majoration_var_X}). Hence we obtain the following bound for the term (\ref{II_III}):
\begin{multline}
\mathbb{P}\Bigg[ \vert\widehat{\tau}-\tau \vert \vert\tau \vert\frac{1}{n}\sum_{i=1}^n\tilde{X}_i^4 \geq\frac{B\sqrt{\log(n)}}{36K\sqrt{n}} \Bigg] \leq \frac{1}{\sqrt{n}}\left(1+2c_{\rm{\scriptscriptstyle BE}}
\frac{\big\Vert \tilde{X}\varepsilon\big\Vert_3^3}{
\big\Vert \tilde{X}\varepsilon\big\Vert_2^3}+ \frac{4\big\Vert \tilde{X}\big\Vert_6^6}{\big\Vert \tilde{X}\big\Vert_4^6}
+\frac{2\big\Vert \tilde{X}\big\Vert_4^2}{\sigma_X^2}\right) \\
+ \frac{288 K \vert\tau\vert\big\Vert \tilde{X}\big\Vert_4^4 \sigma_{\varepsilon}}{\big\Vert \tilde{X}\varepsilon\big\Vert_2^2\sigma_X\sqrt{n\log(n)}}+\frac{2}{n}. \label{borne_II_III}
\end{multline}
\underline{Control of the term (\ref{III_I}) of the upper bound of $D_{2,2}$.}
By using the bound on  $\widehat{\tau}$ given in (\ref{majorationtauchapeau}) and Lemma \ref{lemmeproduit} with 
$a=4\big\Vert \tilde{X}^2 \tilde{Y}\big\Vert_1$, $b=\sqrt{2}\sigma_Y$ and $ c=\sqrt{2}/\sigma_X$, 
we get 
\begin{multline}
\mathbb{P}\Bigg[ \left\vert \frac{1}{n}\sum_{i=1}^n\tilde{X}_i^2\vert \widehat{\tau}\vert \vert\overline{\X}-\mu_{X}\vert\vert 2Y_i-\overline{\Y}-\mu_Y\vert\right\vert \geq\frac{B\sqrt{\log(n)}}{36K\sqrt{n}}\Bigg] \\
\leq \mathbb{P}\Bigg[ \vert\overline{\X}-\mu_X\vert\frac{1}{n}\sum_{i=1}^n{\tilde{X}_i^2\vert 2Y_i-\overline{\Y}-\mu_Y\vert} \frac{\sqrt{\frac{1}{n}\sum_{i=1}^n{(Y_i-\overline{\Y})^2}}}{\sqrt{\frac{1}{n}\sum_{i=1}^n{(X_i-\overline{\X})^2}}} \geq\frac{B\sqrt{\log(n)}}{36K\sqrt{n}}\Bigg] \\
\leq \mathbb{P}\Bigg[ abc\vert\overline{\X}-\mu_X\vert \geq\frac{B\sqrt{\log(n)}}{36K\sqrt{n}}\Bigg] 
+ \mathbb{P}\Bigg[ \frac{1}{n}\sum_{i=1}^n{\tilde{X}_i^2\vert 2Y_i-\overline{\Y}-\mu_Y\vert}  \geq a\Bigg] \\
+ \mathbb{P}\Bigg[ \frac{1}{n}\sum_{i=1}^n{(Y_i-\overline{\Y})^2} \geq b^2\Bigg]+ \mathbb{P}\Bigg[ \sqrt{\frac{1}{n}\sum_{i=1}^n{(X_i-\overline{\X})^2}} \leq 1/c\Bigg]. \label{term_III_I}
\end{multline}
The first probability in (\ref{term_III_I}) is bounded by applying Berry-Esseen Inequality (\ref{BE}) and Pollak Inequality (\ref{Pollak}). Since 
$$\displaystyle 0<K\leq \frac{B}{36\,abc\,\sigma_X}=
\frac{\big\Vert \tilde{X}\varepsilon\big\Vert_2^2\sigma_X}{288 \big\Vert \tilde{X}^2\tilde{Y}\big\Vert_1\sigma_Y}
,$$ we get 
\begin{multline*}
\mathbb{P}\Bigg[ \sqrt{n}\frac{\vert\overline{\X}-\mu_X\vert}{\sigma_X}\geq \frac{B\sqrt{\log(n)}}{36Kabc\sigma_X} \Bigg]
\leq \mathbb{P}\Bigg(\vert \mathcal{N}(0,1)\vert \geq \frac{B\sqrt{\log(n)}}{36Kabc\sigma_X}\Bigg)+2c_{\rm{\scriptscriptstyle BE}}\frac{\mathbb{E}\left( \vert \tilde{X}\vert^3\right)}{\sigma_X^3\sqrt{n}}\\
\leq \exp\left(-\frac{1}{2}\frac{B^2\log(n)}{36^2K^2(abc)^2\sigma_X^2}\right)+2c_{\rm{\scriptscriptstyle BE}}\frac{
\big\Vert\tilde{X}\big\Vert_3^3}{\sigma_X^3\sqrt{n}} \leq \frac{1}{\sqrt{n}}\left(1+2c_{\rm{\scriptscriptstyle BE}}\frac{\big\Vert\tilde{X}\big\Vert_3^3}{\sigma_X^3}\right).
\end{multline*}
The second probability in (\ref{term_III_I}) is bounded by applying Markov Inequality: 
\begin{align}
\mathbb{P}\Bigg[ \frac{1}{n}\sum_{i=1}^n\tilde{X}_i^2\vert 2Y_i-\overline{\Y}-\mu_Y \vert \geq a \Bigg] &\leq 
\mathbb{P}\Bigg[ \frac{1}{n}\sum_{i=1}^n\tilde{X}_i^2\vert 2Y_i-\overline{\Y}-\mu_Y\vert -\frac{a}{2}\geq \frac{a}{2} \Bigg] \nonumber\\
&\qquad \leq \frac{2}{a}\left\Vert\frac{1}{n}\sum_{i=1}^n\tilde{X}_i^2\vert 2Y_i-\overline{\Y}-\mu_Y\vert -\frac{a}{2}\right\Vert_1\nonumber\\
&\qquad \quad\leq \frac{2}{a}\left\Vert\frac{1}{n}\sum_{i=1}^n\tilde{X}_i^2 2\vert \tilde{Y}_i\vert -\frac{a}{2}\right\Vert_1
+ \frac{2}{a}\left\Vert\vert \overline{\Y}-\mu_Y\vert\frac{1}{n}\sum_{i=1}^n\tilde{X}_i^2\right\Vert_1 \nonumber\\
&\qquad \qquad \qquad\leq \frac{4}{a\sqrt{n}}\big\Vert\tilde{X}^2\tilde{Y}\big\Vert_2
+ \frac{2\sigma_Y}{a\sqrt{n}}\big\Vert \tilde{X}\big\Vert_4^2
\label{majoration_terme2deIII_I}.
\end{align}
The last two probabilities in (\ref{term_III_I}) are bounded by  (\ref{majoration_var_Y}) and (\ref{majoration_var_X}).
We thus get the following bound for the term (\ref{III_I})
\begin{align}
\mathbb{P}&\Bigg[\vert \overline{\X}-\mu_X\vert \vert\widehat{\tau}\vert\left\vert \frac{1}{n}\sum_{i=1}^n\tilde{X}_i^2 \vert 2Y_i-\overline{\Y}-\mu_Y\vert \right\vert \geq\frac{B\sqrt{\log(n)}}{36K\sqrt{n}}\Bigg]\nonumber\\
&\leq \frac{1}{\sqrt{n}}\left(1+2c_{\rm{\scriptscriptstyle BE}}\frac{\big\Vert \tilde{X}\big\Vert_3^3}{\sigma_X^3}+ \frac{\big\Vert \tilde{X}^2\tilde{Y}\big\Vert_2}{\big\Vert \tilde{X}^2\tilde{Y}\big\Vert_1}
+ \frac{\sigma_Y \big\Vert \tilde{X}\big\Vert_4^2}{2\big\Vert \tilde{X}^2\tilde{Y}\big\Vert_1}
+ \frac{\big\Vert\tilde{Y}\big\Vert_4^2}{\sigma_Y^2}
+ \frac{2\big\Vert\tilde{X}\big\Vert_4^2}{\sigma_X^2}\right)+\frac{2}{n}.\label{borne_III_I}
\end{align}
\underline{Control of the term (\ref{III_II}) of the upper bound of $D_{2,2}$.}
Using the bound of $\widehat{\tau}$ given by  (\ref{majorationtauchapeau}) and Lemma \ref{lemmeproduit} with 
$a=2\mathbb{E}(\vert\tilde{X}\vert^3)$, $b=2\sigma_Y^2$ and $c=2/\sigma_X^2$, we get 
\begin{multline}
\mathbb{P}\Bigg[ \widehat{\tau}^2\vert\overline{\X}-\mu_{X}\vert \frac{1}{n}\sum_{i=1}^n\tilde{X}_i^2\vert X_i-\overline{\X}\vert \geq\frac{B\sqrt{\log(n)}}{36K\sqrt{n}}\Bigg]  \\
\leq \mathbb{P}\Bigg[ \vert\overline{\X}-\mu_X\vert\frac{1}{n}\sum_{i=1}^n{\tilde{X}_i^2\vert X_i-\overline{\X}\vert} \frac{\frac{1}{n}\sum_{i=1}^n{(Y_i-\overline{\Y})^2}}{\frac{1}{n}\sum_{i=1}^n{(X_i-\overline{\X})^2}} \geq\frac{B\sqrt{\log(n)}}{36K\sqrt{n}}\Bigg]  \\
\leq \mathbb{P}\Bigg[ abc\vert\overline{\X}-\mu_X\vert \geq\frac{B\sqrt{\log(n)}}{36K\sqrt{n}}\Bigg] 
+ \mathbb{P}\Bigg[ \frac{1}{n}\sum_{i=1}^n{\tilde{X}_i^2\vert X_i-\overline{\X}\vert}  \geq a\Bigg] \\
+ \mathbb{P}\Bigg[ \frac{1}{n}\sum_{i=1}^n{(Y_i-\overline{\Y})^2} \geq b\Bigg]
+ \mathbb{P}\Bigg[ \frac{1}{n}\sum_{i=1}^n{(X_i-\overline{\X})^2} \leq 1/c\Bigg]. \label{term_III_II}
\end{multline}
The first probability in (\ref{term_III_II}) is bounded by applying Berry-Esseen Inequality (\ref{BE}) and Pollak Inequality (\ref{Pollak}). Since 
$$\displaystyle 0<K\leq \frac{B}{36abc\sigma_X}=
\frac{\big\Vert \tilde{X}\varepsilon\big\Vert_2^2\sigma_X}{288 \big\Vert \tilde{X}\big\Vert_3^3\sigma_Y^2}
,$$ we get 
\begin{multline*}
\mathbb{P}\Bigg[ \sqrt{n}\frac{\vert\overline{\X}-\mu_X\vert}{\sigma_X}\geq \frac{B\sqrt{\log(n)}}{36Kabc\sigma_X} \Bigg]
\leq 
\mathbb{P}\Bigg(\vert \mathcal{N}(0,1)\vert \geq \frac{B\sqrt{\log(n)}}{36Kabc\sigma_X}\Bigg)+2c_{\rm{\scriptscriptstyle BE}}\frac{\mathbb{E}\left( \vert \tilde{X}\vert^3\right)}{\sigma_X^3\sqrt{n}} \\
\leq \exp\left(-\frac{1}{2}\frac{B^2\log(n)}{36^2K^2(abc)^2\sigma_X^2}\right)+2c_{\rm{\scriptscriptstyle BE}}\frac{
\big\Vert\tilde{X}\big\Vert_3^3}{\sigma_X^3\sqrt{n}} 
\leq \frac{1}{\sqrt{n}}\left(1+2c_{\rm{\scriptscriptstyle BE}}\frac{\big\Vert\tilde{X}\big\Vert_3^3}{\sigma_X^3}\right).
\end{multline*}
The last three probabilities in (\ref{term_III_II}) are bounded using (\ref{majorationterme1-2_Markov}), (\ref{majoration_var_Y}) and (\ref{majoration_var_X}).
Hence we obtain the following bound for the term (\ref{III_II})
\begin{multline}
\mathbb{P}\Bigg[ \widehat{\tau}^2\vert\overline{\X}-\mu_{X}\vert \frac{1}{n}\sum_{i=1}^n\tilde{X}_i^2\vert X_i-\overline{\X}\vert \geq\frac{B\sqrt{\log(n)}}{36K\sqrt{n}}\Bigg] \\
\leq \frac{1}{\sqrt{n}}\left(1+2c_{\rm{\scriptscriptstyle BE}}\frac{\big\Vert \tilde{X}\big\Vert_3^3}{\sigma_X^3}+ \frac{\big\Vert \tilde{X}\big\Vert_6^3 +\sigma_X \big\Vert \tilde{X}\big\Vert_4^2}{\big\Vert \tilde{X}\big\Vert_3^3}
+ \frac{\big\Vert\tilde{Y}\big\Vert_4^2}{\sigma_Y^2}
+ \frac{2\big\Vert\tilde{X}\big\Vert_4^2}{\sigma_X^2}\right) +\frac{2}{n}.\label{borne_III_II}
\end{multline}
\underline{Control of the term (\ref{III_III}) of the upper bound of $D_{2,2}$.}
Using the bound (\ref{majorationtauchapeau}) for $\widehat{\tau}$ and Lemma \ref{lemmeproduit} with 
$a=2\vert\tau\vert\mathbb{E}(\vert\tilde{X}\vert^3)$, $b=\sqrt{2}\sigma_Y$ and $c=\sqrt{2}/\sigma_X$, we get 
\begin{multline}
\mathbb{P}\Bigg[ \vert \widehat{\tau}\vert \vert\overline{\X}-\mu_{X}\vert \vert\tau\vert\frac{1}{n}\sum_{i=1}^n\vert\tilde{X}_i\vert^3 \geq\frac{B\sqrt{\log(n)}}{36K\sqrt{n}}\Bigg] \\
\leq \mathbb{P}\Bigg[ \vert\overline{\X}-\mu_X\vert\vert\tau\vert\frac{1}{n}\sum_{i=1}^n{\vert\tilde{X}_i\vert^3} \frac{\sqrt{\frac{1}{n}\sum_{i=1}^n{(Y_i-\overline{\Y})^2}}}{\sqrt{\frac{1}{n}\sum_{i=1}^n{(X_i-\overline{\X})^2}}} \geq\frac{B\sqrt{\log(n)}}{36K\sqrt{n}}\Bigg] \\
\leq \mathbb{P}\Bigg[ abc\vert\overline{\X}-\mu_X\vert \geq\frac{B\sqrt{\log(n)}}{36K\sqrt{n}}\Bigg] 
+ \mathbb{P}\Bigg[ \vert\tau\vert\frac{1}{n}\sum_{i=1}^n{\vert\tilde{X}_i\vert^3}  \geq a\Bigg] \\
+ \mathbb{P}\Bigg[ \frac{1}{n}\sum_{i=1}^n{(Y_i-\overline{\Y})^2} \geq b^2\Bigg]
+ \mathbb{P}\Bigg[ \sqrt{\frac{1}{n}\sum_{i=1}^n{(X_i-\overline{\X})^2}} \leq 1/c\Bigg]. \label{term_III_III}
\end{multline}
The first probability in (\ref{term_III_III}) is bounded by applying Berry-Esseen Inequality (\ref{BE}) and Pollak Inequality (\ref{Pollak}). Since 
$$\displaystyle 0<K\leq \frac{B}{36abc\sigma_X}=
\frac{\big\Vert \tilde{X}\varepsilon\big\Vert_2^2}{144 \vert\tau\vert\big\Vert \tilde{X}\big\Vert_3^3\sigma_Y}
,$$ we get 
\begin{align*}
\mathbb{P}\Bigg[ \sqrt{n}\frac{\vert\overline{\X}-\mu_X\vert}{\sigma_X}\geq \frac{B\sqrt{\log(n)}}{36Kabc\sigma_X} \Bigg]
&\leq 
\mathbb{P}\Bigg(\vert \mathcal{N}(0,1)\vert \geq \frac{B\sqrt{\log(n)}}{36Kabc\sigma_X}\Bigg)+2c_{\rm{\scriptscriptstyle BE}}\frac{\mathbb{E}\left( \vert \tilde{X}\vert^3\right)}{\sigma_X^3\sqrt{n}} \\
&\quad\leq \exp\left(-\frac{1}{2}\frac{B^2\log(n)}{36^2K^2(abc)^2\sigma_X^2}\right)+2c_{\rm{\scriptscriptstyle BE}}\frac{
\big\Vert\tilde{X}\big\Vert_3^3}{\sigma_X^3\sqrt{n}}  \leq \frac{1}{\sqrt{n}}\left(1+2c_{\rm{\scriptscriptstyle BE}}\frac{\big\Vert\tilde{X}\big\Vert_3^3}{\sigma_X^3}\right).
\end{align*}
The last three probabilities in (\ref{term_III_III}) are bounded respectively by  (\ref{majorationterme1-3_Markov}), (\ref{majoration_var_Y}) and (\ref{majoration_var_X}).
We thus get the following bound for the term (\ref{III_III})
\begin{multline}
\mathbb{P}\Bigg[ \vert \widehat{\tau}\vert \vert\overline{\X}-\mu_{X}\vert \vert\tau\vert\frac{1}{n}\sum_{i=1}^n\vert\tilde{X}_i\vert^3 \geq\frac{B\sqrt{\log(n)}}{36K\sqrt{n}}\Bigg] \\
\leq \frac{1}{\sqrt{n}}\left(1+2c_{\rm{\scriptscriptstyle BE}}\frac{\big\Vert \tilde{X}\big\Vert_3^3}{\sigma_X^3}+ \frac{\big\Vert \tilde{X}   \big\Vert_6^3
}{\big\Vert \tilde{X}   \big\Vert_3^3}
+ \frac{\big\Vert\tilde{Y}\big\Vert_4^2}{\sigma_Y^2}
+ \frac{2\big\Vert\tilde{X}\big\Vert_4^2}{\sigma_X^2}\right)+\frac{2}{n}.
\label{borne_III_III}
\end{multline}
\underline{End of the control of the term $D_{2,2}$.}
Since $K$ satisfies (\ref{K_final}), combining the bounds  (\ref{borne_I_I}), (\ref{borne_I_II}), (\ref{borne_I_III}), (\ref{borne_II_I}), (\ref{borne_II_II}), (\ref{borne_II_III}), (\ref{borne_III_I}), (\ref{borne_III_II}) and (\ref{borne_III_III}),
we get the following bound for the term $D_{2,2}$ 
\begin{multline}
D_{2,2} 
\leq
\frac{1}{\sqrt{n}}\left[9+6c_{\rm{\scriptscriptstyle BE}}\left(\frac{\big\Vert \tilde{Y}\big\Vert_3^3}{\sigma_Y^3}+
\frac{
\big\Vert\tilde{X} \varepsilon\big\Vert_3^3}{\big\Vert \tilde{X}\varepsilon\big\Vert_2^3
}+\frac{\big\Vert \tilde{X}\big\Vert_3^3}{\sigma_X^3}\right)+
2\frac{\big\Vert\tilde{X}^2 \tilde{Y}\big\Vert_2}{\big\Vert \tilde{X}^2\tilde{Y}\big\Vert_1}
+ \frac{\sigma_Y\big\Vert \tilde{X}\big\Vert_4^2}{\big\Vert \tilde{X}^2\tilde{Y}\big\Vert_1}+ 2\frac{\sigma_X \big\Vert \tilde{X}\big\Vert_4^2}{\big\Vert \tilde{X}\big\Vert_3^3}
+ 4\frac{\big\Vert \tilde{X}\big\Vert_6^3}{\big\Vert \tilde{X}\big\Vert_3^3} \right. \\
+ \left.5\frac{\big\Vert \tilde{Y}\big\Vert_4^2}{\sigma_Y^2}
+ 14\frac{\big\Vert \tilde{X}\big\Vert_4^2}{\sigma_X^2}+ 4\frac{\big\Vert\tilde{X}^3\tilde{Y}\big\Vert_{3/2}^{3/2}}{\big\Vert \tilde{X}^3 \tilde{Y}\big\Vert_1^{3/2}}+ \frac{\sigma_Y\big\Vert \tilde{X}\big\Vert_6^3}{4\big\Vert \tilde{X}^3 \tilde{Y}\big\Vert_1}
+8\frac{\big\Vert \tilde{X}\big\Vert_6^6}{\big\Vert \tilde{X}\big\Vert_4^6}
+ \frac{\sigma_X\big\Vert \tilde{X}\big\Vert_6^3}{2\big\Vert \tilde{X}\big\Vert_4^4} \right]\\
+ \frac{1}{\sqrt{n\log(n)}}\frac{288K\sigma_{\varepsilon}}{\big\Vert \tilde{X}\varepsilon\big\Vert_2^2
\sigma_X}\left[4\big\Vert \tilde{X}^3 \tilde{Y}\big\Vert_1 +\frac{4 \sigma_Y \big\Vert \tilde{X}\big\Vert_4^4}{\sigma_X}
+ \vert\tau\vert\big\Vert \tilde{X}\big\Vert_4^4
\right]+\frac{14}{n}. 
\label{majorationD22}
\end{multline}
\underline{End of the control of the term $D_{2}$.}
Combining the bounds (\ref{majorationD21}) for $D_{1,2}$ and (\ref{majorationD22}) for $D_{2,2}$, we get 
\begin{multline}
D_{2}
\leq 
\frac{1}{\sqrt{n}}\Bigg[10+2c_{\rm{\scriptscriptstyle BE}}\left(\frac{
\big\Vert\tilde{X}\varepsilon\big\Vert_6^6}{\left[\Var\left(\tilde{X}^2\varepsilon^2\right)\right]^{3/2}}+3\frac{\big\Vert \tilde{Y}\big\Vert_3^3}{\sigma_Y^3}+
3\frac{
\big\Vert\tilde{X} \varepsilon\big\Vert_3^3}{\big\Vert \tilde{X}\varepsilon\big\Vert_2^3
}+3\frac{\big\Vert \tilde{X}\big\Vert_3^3}{\sigma_X^3}\right)+
2\frac{\big\Vert\tilde{X}^2 \tilde{Y}\big\Vert_2}{\big\Vert \tilde{X}^2\tilde{Y}\big\Vert_1}
+ \frac{\sigma_Y\big\Vert \tilde{X}\big\Vert_4^2}{\big \Vert \tilde{X}^2\tilde{Y}\big\Vert_1} \Bigg.\\
\Bigg.+2\frac{\sigma_X \big\Vert \tilde{X}\big\Vert_4^2}{\big\Vert \tilde{X}\big\Vert_3^3}
+ 4\frac{\big\Vert \tilde{X}\big\Vert_6^3}{\big\Vert \tilde{X}\big\Vert_3^3} 
+ 5\frac{\big\Vert \tilde{Y}\big\Vert_4^2}{\sigma_Y^2}
+ 14\frac{\big\Vert \tilde{X}\big\Vert_4^2}{\sigma_X^2}+ 4\frac{\big\Vert\tilde{X}^3\tilde{Y}\big\Vert_{3/2}^{3/2}}{\big\Vert \tilde{X}^3 \tilde{Y}\big\Vert_1^{3/2}}+ \frac{\sigma_Y\big\Vert \tilde{X}\big\Vert_6^3}{4\big\Vert \tilde{X}^3 \tilde{Y}\big\Vert_1}
+8\frac{\big\Vert \tilde{X}\big\Vert_6^6}{\big\Vert \tilde{X}\big\Vert_4^6}
+ \frac{\sigma_X\big\Vert \tilde{X}\big\Vert_6^3}{2\big\Vert \tilde{X}\big\Vert_4^4} \Bigg]\\
+ \frac{1}{\sqrt{n\log(n)}}\frac{288K\sigma_{\varepsilon}}{\big\Vert \tilde{X}\varepsilon\big\Vert_2^2
\sigma_X}\left[4\big\Vert \tilde{X}^3 \tilde{Y}\big\Vert_1 +\frac{4 \sigma_Y \big\Vert \tilde{X}\big\Vert_4^4}{\sigma_X}
+ \vert\tau\vert\big\Vert \tilde{X}\big\Vert_4^4
\right]+\frac{14}{n}. 
\label{majorationD2}
\end{multline}
\subsubsection*{End of the proof of Lemma \ref{lemme_majoration_P1} }
For $K$ given in (\ref{K_final}), we combine \eqref{Majoration_P1}, the bounds (\ref{majorationD1}) for $D_1$ and (\ref{majorationD2}) for $D_2$. The proof of Lemma \ref{lemme_majoration_P1} is complete. 
\subsection{Proof of  Lemma \ref{lemme2}}
\label{section5.8}
Using Notations (\ref{notations_An_Bn_B}), 
we want to bound
$$\mathbb{P}\left(\left\vert B_n-B \right\vert  \geq\frac{B}{2}\right)=\mathbb{P} \Bigg[   \left\vert   \frac{1}{n}\sum_{i=1}^n(X_i-
\overline{\X})^2{\widehat{\varepsilon_i}}^2 -\big\Vert\tilde{X}\varepsilon\big\Vert_2^2
\right\vert \geq\frac{B}{2}    \Bigg].$$
We start by writing
\begin{align*}
\mathbb{P} &\Bigg[   \left\vert   \frac{1}{n}\sum_{i=1}^n(X_i-
\overline{\X})^2{\widehat{\varepsilon_i}}^2 - \big\Vert\tilde{X}\varepsilon\big\Vert_2^2
\right\vert \geq\frac{B}{2}     \Bigg] \leq 
D'_1+D'_2
\end{align*}
where
\begin{align*}
D'_1=\mathbb{P}\Bigg[ \left\vert    \frac{1}{n}\sum_{i=1}^n\Big[(X_i-\overline{\X})^2-\tilde{X}_i^2\Big]{\widehat{\varepsilon_i}}^2
  \right\vert    \geq \frac{B}{4}        \Bigg]
\mbox{ and }
D'_2=\mathbb{P} &\Bigg[   \left\vert   \frac{1}{n}\sum_{i=1}^n\tilde{X}_i^2{\widehat{\varepsilon_i}}^2 -\big\Vert\tilde{X}\varepsilon\big\Vert_2^2\right\vert \geq\frac{B}{4}  \Bigg].\\
\end{align*}
\subsubsection*{Study of the term $D'_1$.} 
We have 
\begin{align*}
D'_1\leq  \mathbb{P}\Bigg[ \left\vert  \overline{\X}-\mu_{X}\right\vert \frac{1}{n}\sum_{i=1}^n   {\widehat{\varepsilon}_i}^2  \left\vert  2X_i-\overline{\X}-\mu_{X}\right\vert \geq\frac{B}{4}  \Bigg].
\end{align*}
Using the  bound (\ref{varepsilonchapo2}) for ${\widehat{\varepsilon_i}}^2$, 
we get  
\begin{align*}
D'_1\leq D'_{1,1}+D'_{1,2}
\end{align*}
with
\begin{align*}
D'_{1,1}&=\mathbb{P} \Bigg[  \left\vert  \overline{\X}-\mu_{X}\right\vert    \frac{1}{n}\sum_{i=1}^n(Y_i-
\overline{\Y})^2\left\vert 2X_i-\overline{\X}-\mu_{X}\right\vert \geq
\frac{B}{16}     \Bigg]
\end{align*}
and
\begin{align*}
D'_{1,2}&=\mathbb{P} \Bigg[  \left\vert  \overline{\X}-\mu_{X}\right\vert  \widehat{\tau}^2  \frac{1}{n}\sum_{i=1}^n(X_i-
\overline{\X})^2\left\vert 2X_i-\overline{\X}-\mu_{X}\right\vert \geq
\frac{B}{16}     \Bigg].
\end{align*}

\noindent
\underline{Study of the term $D'_{1,1}$.}
Applying Lemma \ref{lemmeproduit} with $a=4\big\Vert\tilde{Y}^2 \tilde{X}\big\Vert_1$, we get  
\begin{align}
\label{D'1_1}
D'_{1,1}&\leq \mathbb{P}\Bigg[a\vert \overline{\X}-\mu_{X}\vert \geq \frac{B}{16}\Bigg]+
\mathbb{P}\Bigg[   \frac{1}{n}\sum_{i=1}^n(Y_i-
\overline{\Y})^2\left\vert 2X_i-\overline{\X}-\mu_{X}\right\vert \geq a   \Bigg].
\end{align}
The first probability in  (\ref{D'1_1}) is bounded using Markov Inequality
\begin{align}
\label{D'1_1_terme1}
\mathbb{P}\Bigg[a\vert \overline{\X}-\mu_{X}\vert \geq \frac{B}{16}\Bigg]&\leq 
\frac{16}{B}\big\Vert a( \overline{\X}-\mu_{X})\big\Vert_1
\leq \frac{64\big\Vert\tilde{Y}^2\tilde{X}\big\Vert_1\sigma_X}{\sqrt{n}\big\Vert\tilde{X}\varepsilon\big\Vert_2^2}.
\end{align}
The second probability in  (\ref{D'1_1}) is bounded using (\ref{D11partie2}). 
Then, by combining (\ref{D'1_1_terme1}) and (\ref{D11partie2}), we get 
\begin{multline}
\label{majorationD'11}
D'_{1,1}\leq
\frac{1}{\sqrt{n}}\left(\frac{64
\big\Vert\tilde{Y}^2\tilde{X}\big\Vert_1\sigma_X}{\big\Vert\tilde{X}\varepsilon\big\Vert_2^2}+\frac{\big\Vert
\tilde{Y}^2\tilde{X}\big\Vert_2}{\big\Vert\tilde{Y}^2 \tilde{X}\big\Vert_1}
+\frac{\sigma_X\big\Vert \tilde{Y}\big\Vert_4^2}{2\big\Vert\tilde{Y}^2 \tilde{X}\big\Vert_1}+
\frac{4\sigma_Y
\big\Vert \tilde{X} \tilde{Y}\big\Vert_2}{2\big\Vert\tilde{Y}^2 \tilde{X}\big\Vert_1}\right)\\ 
+\frac{8c_{1,\rm{\scriptscriptstyle R}}^2}{n}\frac{\sigma_X\sigma_Y^2}{\big\Vert\tilde{Y}^2 \tilde{X}\big\Vert_1}+\frac{32c_{2,\rm{\scriptscriptstyle R}}^2}{n^{3/2}}\frac{\sigma_X\big\Vert\tilde{Y}\big\Vert_4^2}{\big\Vert\tilde{Y}^2 \tilde{X}\big\Vert_1}.
\end{multline}
\underline{Study of the term $D'_{1,2}$.} Applying Lemma \ref{lemmeproduit} with  $a=2\sigma_Y^2, b=4 \big\Vert\tilde{X}\big\Vert_3^3$ and $c=2/\sigma_X^2$, we get 
\begin{multline}
\label{D'1_2}
D'_{1,2}\leq \mathbb{P}\Bigg[abc\vert \overline{\X}-\mu_{X}\vert \geq \frac{B}{16}\Bigg]+
\mathbb{P}\Bigg[   \frac{1}{n}\sum_{i=1}^n(Y_i-
\overline{\Y})^2 \geq a   \Bigg] \\
+\mathbb{P}\Bigg[   \frac{1}{n}\sum_{i=1}^n(X_i-
\overline{\X})^2\left\vert 2X_i-\overline{\X}-\mu_{X}\right\vert \geq b   \Bigg]
+\mathbb{P}\Bigg[   \frac{1}{n}\sum_{i=1}^n(X_i-
\overline{\X})^2\leq\frac{1}{c}\Bigg].
\end{multline}
The first probability in (\ref{D'1_2}) is bounded by Markov Inequality:
\begin{eqnarray}
\label{D'1_2_terme1}
\mathbb{P}\Bigg[abc\vert \overline{\X}-\mu_{X}\vert \geq \frac{B}{16}\Bigg]\leq \frac{16}{B}\big\Vert abc(\overline{\X}-\mu_{X})\big\Vert_1 \leq \frac{256\big\Vert\tilde{X}\big\Vert_3^3
\sigma_Y^2}{\sqrt{n}\big\Vert\tilde{X}\varepsilon\big\Vert_2^2
\sigma_X}.
\end{eqnarray}
Then by combining (\ref{D'1_2_terme1}), (\ref{majoration_var_Y}),
(\ref{D12partie2}) and (\ref{majoration_var_X}),  we obtain 
\begin{multline*}
D'_{1,2}\leq \frac{2}{\sqrt{n}}\left(\frac{128\big\Vert\tilde{X}\big\Vert_3^3\sigma_Y^2}{
\big\Vert\tilde{X}\varepsilon\big\Vert_2^2
\sigma_X}+\frac{
\big\Vert\tilde{X}\big\Vert_4^2
}{\sigma_X^2}+\frac{
\big\Vert\tilde{Y}\big\Vert_4^2
}{2\sigma_Y^2}
+\frac{\big\Vert\tilde{X} \big\Vert_6^3 
}{4\big\Vert\tilde{X}\big\Vert_3^3} + \frac{5\sigma_X\big\Vert\tilde{X}\big\Vert_4^2}{8\big\Vert\tilde{X}\big\Vert_3^3} \right)\\
+\frac{2}{n}\left(1+\frac{2c_{1,\rm{\scriptscriptstyle R}}^2\sigma_X^3}{\big\Vert\tilde{X}\big\Vert_3^3}\right)+\frac{16c_{2,\rm{\scriptscriptstyle R}}^2}{n^{3/2}}\frac{\sigma_X\big\Vert\tilde{X}\big\Vert_4^2}{\big\Vert\tilde{X}\big\Vert_3^3 }.
\end{multline*}
From this last inequality combined with inequality (\ref{majorationD'11}), we conclude that
\begin{multline}
D'_{1}\leq 
\frac{1}{\sqrt{n}}\frac{1}{2\big\Vert\tilde{Y}^2 \tilde{X}\big\Vert_1}\left(2\big\Vert\tilde{Y}^2\tilde{X}\big\Vert_2+\sigma_X\big\Vert \tilde{Y}\big\Vert_4^2+
4\sigma_Y
\big\Vert\tilde{X}\tilde{Y}\big\Vert_2\right)  \\ 
+\frac{2}{\sqrt{n}}\left(\frac{32\big\Vert
\tilde{Y}^2 \tilde{X}\big\Vert_1\sigma_X}
{\big\Vert\tilde{X}\varepsilon\big\Vert_2^2}+
\frac{128\big\Vert\tilde{X}\big\Vert_3^3\sigma_Y^2}{
\big\Vert\tilde{X}\varepsilon\big\Vert_2^2
\sigma_X} +\frac{\big\Vert\tilde{X}\big\Vert_4^2}{\sigma_X^2}+\frac{\big\Vert\tilde{Y}\big\Vert_4^2}{2\sigma_Y^2}
+\frac{\big\Vert\tilde{X}\big\Vert_6^3}{4\big\Vert\tilde{X}\big\Vert_3^3} + \frac{5\sigma_X\big\Vert\tilde{X}\big\Vert_4^2}{8\big\Vert\tilde{X}\big\Vert_3^3} \right)\\
+\frac{2}{n}\left(1+\frac{2c_{1,\rm{\scriptscriptstyle R}}^2\sigma_X^3}{\big\Vert\tilde{X}\big\Vert_3^3}+\frac{4c_{1,\rm{\scriptscriptstyle R}}^2\sigma_X\sigma_Y^2}{\big\Vert\tilde{Y}^2 \tilde{X}\big\Vert_1}\right) +\frac{16c_{2,\rm{\scriptscriptstyle R}}^2}{n^{3/2}}\left(\frac{\sigma_X\big\Vert\tilde{X}\big\Vert_4^2}{\big\Vert\tilde{X}\big\Vert_3^3 } +
\frac{2\sigma_X\big\Vert\tilde{Y}\big\Vert_4^2}{\big\Vert\tilde{Y}^2 \tilde{X}\big\Vert_1}\right).
\label{majorationD'1}
\end{multline}
\subsubsection*{Study of the term $D'_2$} 
We have
 \begin{align*}
D'_2&=\mathbb{P} \Bigg[   \left\vert   \frac{1}{n}\sum_{i=1}^n\tilde{X}_i^2{\widehat{\varepsilon_i}}^2 -\big\Vert\tilde{X}\varepsilon\big\Vert_2^2
\right\vert \geq   \frac{B}{4}    \Bigg]\leq D'_{2,1}+D'_{2,2}
\end{align*}
with
\begin{align*}
D'_{2,1}=\mathbb{P} \Bigg[   \left\vert   \frac{1}{n}\sum_{i=1}^n\tilde{X}_i^2{\varepsilon_i}^2 -
\big\Vert\tilde{X}\varepsilon\big\Vert_2^2
\right\vert \geq\frac{B}{8}
    \Bigg]
 \mbox{ and }
D'_{2,2}=\mathbb{P} \Bigg[   \left\vert   \frac{1}{n}\sum_{i=1}^n\tilde{X}_i^2[{\widehat{\varepsilon_i}}^2-{\varepsilon_i}^2] \right\vert \geq
\frac{B}{8}      \Bigg].
\end{align*}
\underline{Study of the term $D'_{2,1}$.} Applying Markov Inequality and von Bahr-Esseen Inequality (\ref{Von_Bahr_Esseen}), we get 
\begin{multline}
D'_{2,1} \leq \left (\frac{8}{\big\Vert\tilde{X}\varepsilon\big\Vert_2^2} \right )^{3/2}  \left\Vert   \frac{1}{n}\sum_{i=1}^n\tilde{X}_i^2{\varepsilon_i}^2 -\big\Vert\tilde{X}\varepsilon\big\Vert_2^2
\right\Vert_{3/2}^{3/2} \leq 
\frac{8^{3/2}}{\big\Vert\tilde{X}\varepsilon\big\Vert_2^3}\frac{\sqrt 2}{\sqrt n} \left\Vert \tilde{X}^2\varepsilon^2-
\big\Vert\tilde{X}\varepsilon\big\Vert_2^2
\right \Vert_{3/2}^{3/2}   
 \leq \frac{64 \sqrt 2\big\Vert\tilde{X}\varepsilon\big\Vert_3^3}{\sqrt n \big\Vert\tilde{X}\varepsilon\big\Vert_2^3} \, .
\label{majorationD'21}
\end{multline}
\underline{Study of the term $D'_{2,2}$.} We bound the term  $D'_{2,2}$ in the same way as $D_{2,2}$ and we obtain that $D'_{2,2}$ is bounded by the sum of the 9 following terms: 
\begin{align}
\mathbb{P}&\Bigg[ \frac{1}{n}\sum_{i=1}^n\tilde{X}_i^2\vert\mu_Y-\overline{\Y}\vert\vert 2Y_i-\overline{\Y}-\mu_Y\vert\geq\frac{B}{72} \Bigg],\label{I_Iprime}\\
 \mathbb{P}&\Bigg[ \frac{1}{n}\sum_{i=1}^n\tilde{X}_i^2\vert \mu_Y-\overline{\Y}\vert \vert\widehat{\tau}\vert\vert X_i-\overline{\X}\vert  \geq\frac{B}{72}\Bigg],
\label{I_IIprime}\\
 \mathbb{P}&\Bigg[ \frac{1}{n}\sum_{i=1}^n\tilde{X}_i^2\vert\mu_Y-\overline{\Y}\vert \vert \tau  \vert   \vert \tilde{X}_i\vert \geq\frac{B}{72}\Bigg]\label{I_IIIprime},\\
 \mathbb{P}&\Bigg[  \frac{1}{n}\sum_{i=1}^n\vert\tilde{X}_i\vert^3\vert \widehat{\tau}-\tau \vert \vert2Y_i-\overline{\Y}-\mu_Y\vert 
\geq\frac{B}{72}\Bigg],\label{II_Iprime}\\
 \mathbb{P}&\Bigg[  \frac{1}{n}\sum_{i=1}^n\vert\tilde{X}_i\vert^3\vert \widehat{\tau}-\tau \vert \vert \widehat{\tau}\vert\vert X_i-\overline{\X}\vert
  \geq\frac{B}{72}\Bigg],\label{II_IIprime} \\
  \mathbb{P}&\Bigg[  \frac{1}{n}\sum_{i=1}^n\vert\tilde{X}_i\vert^3\vert \widehat{\tau}-\tau \vert  \vert \tau   \vert \vert\tilde{X}_i\vert\geq\frac{B}{72}\Bigg],\label{II_IIIprime}\\
 \mathbb{P}&\Bigg[  \frac{1}{n}\sum_{i=1}^n\tilde{X}_i^2\vert \widehat{\tau}\vert \vert\overline{\X}-\mu_{X}\vert\vert2Y_i-\overline{\Y}-\mu_Y\vert\geq\frac{B}{72}\Bigg],\label{III_Iprime}\\
    \mathbb{P}&\Bigg[  \frac{1}{n}\sum_{i=1}^n\tilde{X}_i^2 \widehat{\tau}^2 \vert\overline{\X}-\mu_{X}\vert \vert X_i-\overline{\X}\vert  \geq\frac{B}{72}\Bigg],\label{III_IIprime}\\
    \mathbb{P}&\Bigg[  \frac{1}{n}\sum_{i=1}^n\tilde{X}_i^2\vert \widehat{\tau}\vert \vert\overline{\X}-\mu_{X}\vert   \vert \tau  \vert  \vert\tilde{X}_i\vert
     \geq\frac{B}{72}\Bigg]\label{III_IIIprime}. 
\end{align}
We now study each of these terms. \\

\noindent
\underline{Control of the term (\ref{I_Iprime}) of the upper bound of $D'_{2,2}$.}
Applying Lemma \ref{lemmeproduit} with $a=4\big\Vert\tilde{X}^2\tilde{Y}\big\Vert_1$, we get 
\begin{multline}
\mathbb{P}\Bigg[ \vert\overline{\Y}-\mu_Y\vert\frac{1}{n}\sum_{i=1}^n\tilde{X}_i^2\vert 2Y_i-\overline{\Y}-\mu_Y\vert \geq\frac{B}{72} \Bigg]\leq \mathbb{P}\Bigg[ a\vert\overline{\Y}-\mu_Y\vert \geq\frac{B}{72}\Bigg]+\mathbb{P}\Bigg[ \frac{1}{n}\sum_{i=1}^n\tilde{X}_i^2\vert 2Y_i-\overline{\Y}-\mu_Y\vert \geq a \Bigg]. \label{term_I_Iprime}
\end{multline}
The first probability in (\ref{term_I_Iprime}) is bounded using Markov Inequality
\begin{align*}
\mathbb{P}\Bigg[ a\vert\overline{\Y}-\mu_Y\vert\geq\frac{B}{72} \Bigg]\leq \frac{288\,\sigma_Y\big\Vert\tilde{X}^2 \tilde{Y}\big\Vert_1}{\sqrt{n}
\big\Vert\tilde{X}\varepsilon\big\Vert_2^2}.
\end{align*}
\noindent
The second probability in (\ref{term_I_Iprime}) is bounded by (\ref{majoration_terme2_de_II}). 
Hence we get the following bound for the term (\ref{I_Iprime}): 
\begin{multline}
\mathbb{P}\Bigg[ \vert\overline{\Y}-\mu_Y\vert\frac{1}{n}\sum_{i=1}^n\tilde{X}_i^2\vert 2Y_i-\overline{\Y}-\mu_Y\vert \geq\frac{B}{72} \Bigg] 
\leq 
\frac{1}{\sqrt{n}}\left(\frac{288\,\sigma_Y
\big\Vert\tilde{X}^2 \tilde{Y}\big\Vert_1}{\big\Vert\tilde{X}\varepsilon\big\Vert_2^2}+
\frac{\big\Vert\tilde{X}^2\tilde{Y}\big\Vert_2}{\big\Vert\tilde{X}^2 \tilde{Y}\big\Vert_1}
+ \frac{\sigma_Y\big\Vert \tilde{X}\big\Vert_4^2}{2\big\Vert\tilde{X}^2 \tilde{Y}\big\Vert_1}\right).
\label{borne_I_Iprime}
\end{multline}
\underline{Control of the term (\ref{I_IIprime}) of the upper bound of $D'_{2,2}$.}
We start by bounding the term (\ref{I_IIprime}) using the bound on $ \widehat{\tau}$ given in (\ref{majorationtauchapeau}). Then, applying Lemma \ref{lemmeproduit} with
$a=2\big\Vert\tilde{X}\big\Vert_3^3$, $b=\sqrt{2}\sigma_Y$, $ c=\sqrt{2}/\sigma_X,$ we get
\begin{multline}
\mathbb{P}\Bigg[ \vert\overline{\Y}-\mu_Y\vert\frac{1}{n}\sum_{i=1}^n{\tilde{X}_i^2\vert X_i-\overline{\X}\vert} \frac{\sqrt{\frac{1}{n}\sum_{i=1}^n{(Y_i-\overline{\Y})^2}}}{\sqrt{\frac{1}{n}\sum_{i=1}^n{(X_i-\overline{\X})^2}}} \geq\frac{B}{72}\Bigg] \\
\leq \mathbb{P}\Bigg[ abc\vert\overline{\Y}-\mu_Y\vert \geq\frac{B}{72}\Bigg] 
+ \mathbb{P}\Bigg[ \frac{1}{n}\sum_{i=1}^n{\tilde{X}_i^2\vert X_i-\overline{\X}\vert}  \geq a\Bigg] \\
+ \mathbb{P}\Bigg[ \frac{1}{n}\sum_{i=1}^n{(Y_i-\overline{\Y})^2} \geq b^2\Bigg]
 + \mathbb{P}\Bigg[ \sqrt{\frac{1}{n}\sum_{i=1}^n{(X_i-\overline{\X})^2}} \leq 1/c\Bigg].\label{term_I_IIprime}
\end{multline}
The first probability in (\ref{term_I_IIprime}) is bounded by applying Markov Inequality 
\begin{align*}
\mathbb{P}\Bigg[ abc\vert\overline{\Y}-\mu_Y\vert\geq \frac{B}{72} \Bigg]
\leq \frac{288 \big\Vert\tilde{X}\big\Vert_3^3\sigma_Y^2}{\sqrt{n}\big\Vert\tilde{X}\varepsilon\big\Vert_2^2\sigma_X}.
\end{align*}
The other three probabilities in (\ref{term_I_IIprime}) are bounded by (\ref{majorationterme1-2_Markov}), (\ref{majoration_var_Y}) and (\ref{majoration_var_X}). 
Hence, we obtain the following bound for the term (\ref{I_IIprime}): 
\begin{multline}
\mathbb{P}\Bigg[\vert \overline{\Y}-\mu_Y\vert \vert\widehat{\tau}\vert\frac{1}{n}\sum_{i=1}^n\tilde{X}_i^2 \vert X_i-\overline{\X}\vert  \geq\frac{B}{72}\Bigg]\\
\leq \frac{1}{\sqrt{n}}\left(\frac{288 \big\Vert\tilde{X}\big\Vert_3^3\sigma_Y^2}{\big\Vert\tilde{X}\varepsilon\big\Vert_2^2\sigma_X}+ \frac{\big\Vert \tilde{X}\big\Vert_6^3}{\big\Vert \tilde{X}\big\Vert_3^3}
+ \frac{\sigma_X \big\Vert \tilde{X}\big\Vert_4^2}{\big\Vert \tilde{X}\big\Vert_3^3}
+ \frac{\big\Vert\tilde{Y}\big\Vert_4^2}{\sigma_Y^2}
+ \frac{2\big\Vert\tilde{X}\big\Vert_4^2}{\sigma_X^2}\right)+\frac{2}{n}.\label{borne_I_IIprime}
\end{multline}
\underline{Control of the term (\ref{I_IIIprime}) of the upper bound of $D'_{2,2}$.} We apply Lemma \ref{lemmeproduit} with $a=2\vert\tau\vert \Vert \tilde{X}\Vert_3^3$ to get 
\begin{align}
\mathbb{P}\Bigg[ \vert\overline{\Y}-\mu_Y\vert\vert\tau \vert\frac{1}{n}\sum_{i=1}^n{\vert \tilde{X}_i\vert^3} \geq\frac{B}{72} \Bigg]&\leq \mathbb{P}\Bigg[ a\vert\overline{\Y}-\mu_Y\vert \geq\frac{B}{72}\Bigg]+\mathbb{P}\Bigg[ \vert\tau \vert\frac{1}{n}\sum_{i=1}^n{\vert \tilde{X}_i\vert^3} \geq a \Bigg]. \label{term_I_IIIprime}
\end{align}
The first probability in (\ref{term_I_IIIprime}) is bounded by applying Markov Inequality: 
\begin{align*}
\mathbb{P}\Bigg[ a\vert\overline{\Y}-\mu_Y\vert\geq\frac{B}{72} \Bigg]
\leq \frac{144 \vert\tau \vert\big\Vert\tilde{X}\big\Vert_3^3\sigma_Y}{\sqrt{n}\big\Vert\tilde{X}\varepsilon\big\Vert_2^2}.
\end{align*}
The second probability in (\ref{term_I_IIIprime}) is bounded by applying Inequality
(\ref{majorationterme1-3_Markov}). 
Then, the term (\ref{I_IIIprime}) satisfies
\begin{align}
\mathbb{P}\Bigg[ \vert\overline{\Y}-\mu_Y\vert\vert\tau \vert\frac{1}{n}\sum_{i=1}^n{\vert \tilde{X}_i\vert^3} \geq\frac{B}{72} \Bigg]&\leq 
\frac{1}{\sqrt{n}}\left(
\frac{144 \vert\tau \vert\big\Vert \tilde{X}\big\Vert_3^3\sigma_Y}
{\big\Vert\tilde{X}\varepsilon\big\Vert_2^2}+ \frac{ 
\big\Vert \tilde{X}   \big\Vert_6^3
}{\big\Vert \tilde{X}   \big\Vert_3^3}\right) . \label{borne_I_IIIprime}
\end{align}
\underline{Control of the term (\ref{II_Iprime}) of the upper bound of $D'_{2,2}$.} Applying Lemma \ref{lemmeproduit} with $a=4\Vert \tilde{X}^3 \tilde{Y}\Vert_1$ and $b={2}/{\sigma_X^2}$, we get 
\begin{multline}
\mathbb{P}\Bigg[ \vert\widehat{\tau}-\tau \vert\frac{1}{n}\sum_{i=1}^n\vert \tilde{X}_i\vert^3\vert 2Y_i-\overline{\Y}-\mu_Y\vert \geq\frac{B}{72} \Bigg] \\
\leq \mathbb{P}\Bigg[ ab\left\vert\frac{1}{n}\sum_{i=1}^n{\tilde{X}_i\varepsilon_i}\right\vert  \geq\frac{B}{144}\Bigg]
+\mathbb{P}\Bigg[ ab\vert\overline{\X}-\mu_X\vert\left\vert\frac{1}{n}\sum_{i=1}^n{\varepsilon_i}\right\vert  \geq\frac{B}{144}\Bigg] \\
+\mathbb{P}\Bigg[ \frac{1}{n}\sum_{i=1}^n{\vert \tilde{X}_i\vert^3\vert 2Y_i-\overline{\Y}-\mu_Y\vert} \geq a \Bigg] +\mathbb{P}\Bigg[ \frac{1}{n}\sum_{i=1}^n{(X_i-\overline{\X})^2}\leq 1/b \Bigg]. \label{term_II_Iprime}
\end{multline}
The first two probabilities in (\ref{term_II_Iprime}) are bounded using Markov Inequality: 
\begin{align*}
\mathbb{P}\Bigg[ ab\left\vert 
\frac{1}{n}\sum_{i=1}^n{\tilde{X}_i\varepsilon_i}\right\vert\geq\frac{B}{144} \Bigg]
\leq \frac{1152\big\Vert \tilde{X}^3\tilde{Y}\big\Vert_1}{\sqrt{n}\sigma_X^2\big\Vert\tilde{X}\varepsilon\big\Vert_2}
\end{align*}
and
\begin{align*}
\mathbb{P}\Bigg[ ab\left\vert
\overline{\X}-\mu_X\right\vert\left\vert\frac{1}{n}\sum_{i=1}^n{\varepsilon_i}\right\vert  \geq\frac{B}{144}\Bigg]
&\leq \frac{144\,ab}{B}\big\Vert
\overline{\X}-\mu_X\big\Vert_2\left\Vert\frac{1}{n}\sum_{i=1}^n{\varepsilon_i}\right\Vert_2 \leq \frac{1152\big\Vert\tilde{X}^3\tilde{Y}\big\Vert_1  \sigma_{\varepsilon}}{n\,\sigma_X\big\Vert\tilde{X}\varepsilon\big\Vert_2^2}.
\end{align*}
The last two probabilities in (\ref{term_II_Iprime}) are bounded by (\ref{majoration_Markov_proba3_de_termeII_I}) and (\ref{majoration_var_X}).
Hence we obtain the following bound for the term (\ref{II_Iprime}) 
\begin{multline}
\mathbb{P}\Bigg[ \vert\widehat{\tau}-\tau \vert\frac{1}{n}\sum_{i=1}^n\vert \tilde{X}_i\vert^3\vert 2Y_i-\overline{\Y}-\mu_Y\vert \geq\frac{B}{72} \Bigg] \\
\leq \frac{1}{\sqrt{n}}\left(\frac{1152\big\Vert \tilde{X}^3\tilde{Y} \big\Vert_1}{\sigma_X^2 \big\Vert \tilde{X}\varepsilon \big\Vert_2
}+ 2^{7/2}\frac{\big\Vert\tilde{X}^3\tilde{Y}\big\Vert_{3/2}^{3/2}}{\big\Vert \tilde{X}^3 \tilde{Y}\big\Vert_1^{3/2}}+ \frac{\sigma_Y\big\Vert \tilde{X}\big\Vert_6^3}{2\big\Vert \tilde{X}^3 \tilde{Y}\big\Vert_1}+\frac{2\big\Vert \tilde{X}\big\Vert_4^2}{\sigma_X^2}\right)
+ \frac{1152 \big\Vert\tilde{X}^3 \tilde{Y}\big\Vert_1
\sigma_{\varepsilon}}{
n\sigma_X\big\Vert\tilde{X}\varepsilon\big\Vert_2^2
}+\frac{2}{n}.\label{borne_II_Iprime}
\end{multline}
\underline{Control of the term (\ref{II_IIprime}) of the upper bound of $D'_{2,2}$.}
Applying Lemma \ref{lemmeproduit} with $a=(\sqrt{2}/\sigma_X)^3$, $b=\sqrt{2}\sigma_Y$ and $c=2
\big\Vert\tilde{X}\big\Vert_4^4$, we get 
\begin{multline}
\mathbb{P}\Bigg[ \left\vert
\frac{1}{n}\sum_{i=1}^n{(X_i-
\overline{\X})\varepsilon_i}
\right\vert
\frac{\sqrt{\frac{1}{n}\sum_{i=1}^n{(Y_i-
\overline{\Y})^2}}}{\left[\frac{1}{n}\sum_{i=1}^n{(X_i-
\overline{\X})^2}\right]^{3/2}}
\frac{1}{n}\sum_{i=1}^n\vert \tilde{X}_i\vert^3\vert X_i-\overline{\X}\vert \geq\frac{B}{72} \Bigg]\\
\leq \mathbb{P}\Bigg[ abc\left\vert\frac{1}{n}\sum_{i=1}^n{\tilde{X}_i\varepsilon_i}\right\vert  \geq\frac{B}{144}\Bigg] + \mathbb{P}\Bigg[ abc\vert\overline{\X}-\mu_X\vert\left\vert\frac{1}{n}\sum_{i=1}^n{\varepsilon_i}\right\vert  \geq\frac{B}{144}\Bigg] 
+\mathbb{P}\Bigg[ \frac{1}{n}\sum_{i=1}^n{(X_i-
\overline{\X})^2}\leq \frac{\sigma_X^2}{2} \Bigg] \\
+\mathbb{P}\Bigg[ \frac{1}{n}\sum_{i=1}^n{(Y_i-
\overline{\Y})^2}\geq 2\sigma_Y^2 \Bigg] + 
\mathbb{P}\Bigg[ 
\frac{1}{n}\sum_{i=1}^n\vert \tilde{X}_i\vert^3\vert X_i-\overline{\X}\vert \geq c \Bigg] \label{term_II_IIprime}.
\end{multline}
The first two probabilities in (\ref{term_II_IIprime}) are bounded using Markov Inequality: 
\begin{align*}
\mathbb{P}\Bigg[ abc\left\vert 
\frac{1}{n}\sum_{i=1}^n{\tilde{X}_i\varepsilon_i}\right\vert\geq\frac{B}{144} \Bigg]
\leq \frac{1152 \big\Vert\tilde{X}\big\Vert_4^4\sigma_Y}{\sqrt{n}\,\sigma_X^3\big\Vert\tilde{X}\varepsilon \big\Vert_2}
\end{align*}
and
\begin{align*}
\mathbb{P}\Bigg[ abc\left\vert
\overline{\X}-\mu_X\right\vert\left\vert\frac{1}{n}\sum_{i=1}^n{\varepsilon_i}\right\vert  \geq\frac{B}{144}\Bigg]
\leq \frac{144 abc}{B}\left\Vert
\overline{\X}-\mu_X\right\Vert_2\left\Vert\frac{1}{n}\sum_{i=1}^n{\varepsilon_i}\right\Vert_2 
\leq \frac{1152 \sigma_Y\sigma_{\varepsilon}\big\Vert \tilde{X}\big\Vert_4^4
}{n\sigma_X^2\big\Vert\tilde{X}\varepsilon \big\Vert_2^2}.
\end{align*}
The last three probabilities in (\ref{term_II_IIprime}) are bounded by (\ref{majoration_var_X}), (\ref{majoration_var_Y}) and (\ref{majoration_terme5deII_II}). 
Consequently, the term (\ref{II_IIprime}) is bounded as follows
\begin{multline}
\mathbb{P}\Bigg[ 
\left\vert \widehat{\tau}(\widehat{\tau}-\tau ) \right\vert\frac{1}{n}\sum_{i=1}^n\vert \tilde{X}_i\vert ^3\vert X_i-\overline{\X}\vert 
 \geq\frac{B}{72}\Bigg] 
\\
\leq \frac{1}{\sqrt{n}}\left(\frac{1152\big\Vert\tilde{X}\big\Vert_4^4
\sigma_Y}{\sigma_X^3\big\Vert\tilde{X}\varepsilon\big\Vert_2}
+\frac{\big\Vert \tilde{Y}\big\Vert_4^2}{\sigma_Y^2}+\frac{2\big\Vert \tilde{X}\big\Vert_4^2}{\sigma_X^2} +4\frac{\big\Vert \tilde{X}\big\Vert_6^6}{\big\Vert \tilde{X}\big\Vert_4^6}
+ \frac{\sigma_X\big\Vert \tilde{X}\big\Vert_6^3}{2\big\Vert \tilde{X}\big\Vert_4^4} \right) 
 +\frac{1}{n}\left(\frac{1152 \,\sigma_Y \sigma_{\varepsilon}\big\Vert \tilde{X}\big\Vert_4^4}{\big\Vert\tilde{X}\varepsilon\big\Vert_2^2\sigma_X^2}+2 \right).\label{borne_II_IIprime}
\end{multline}
\underline{Control of the term (\ref{II_IIIprime}) of the upper bound of $D'_{2,2}$.} Applying Lemma \ref{lemmeproduit} with 
$a=2\vert\tau\vert\Vert \tilde{X}\Vert^4_4$ and $b={2}/{\sigma_X^2}$, we get 
\begin{multline}
\mathbb{P}\Bigg[ \vert\widehat{\tau}-\tau\vert \vert\tau\vert \frac{1}{n}\sum_{i=1}^n\tilde{X}_i^4\geq\frac{B}{72}\Bigg] \leq \mathbb{P}\Bigg[ ab\left\vert\frac{1}{n}\sum_{i=1}^n{(X_i-
\overline{\X})\varepsilon_i}\right\vert  \geq\frac{B}{72}\Bigg]
+\mathbb{P}\Bigg[ \vert\tau\vert \frac{1}{n}\sum_{i=1}^n{\vert \tilde{X}_i\vert^4} \geq a \Bigg] \\ 
+\mathbb{P}\Bigg[ \frac{1}{n}\sum_{i=1}^n{(X_i-\overline{\X})^2}\leq 1/b \Bigg] \\
\leq \mathbb{P}\Bigg[ ab\left\vert\frac{1}{n}\sum_{i=1}^n{\tilde{X}_i\varepsilon_i}\right\vert  \geq\frac{B}{144}\Bigg]
+\mathbb{P}\Bigg[ ab\vert\overline{\X}-\mu_X\vert\left\vert\frac{1}{n}\sum_{i=1}^n{\varepsilon_i}\right\vert  \geq\frac{B}{144}\Bigg] \\
+\mathbb{P}\Bigg[ \vert\tau\vert \frac{1}{n}\sum_{i=1}^n{\vert \tilde{X}_i\vert^4} \geq a \Bigg] +\mathbb{P}\Bigg[ \frac{1}{n}\sum_{i=1}^n{(X_i-\overline{\X})^2}\leq 1/b \Bigg]
\label{term_II_IIIprime}.
\end{multline}
The first two probabilities in (\ref{term_II_IIIprime}) are bounded by applying Markov Inequality: 
\begin{align*}
\mathbb{P}\Bigg[ ab\left\vert 
\frac{1}{n}\sum_{i=1}^n{\tilde{X}_i\varepsilon_i}\right\vert\geq\frac{B}{144
} \Bigg]
&\leq \frac{576\vert\tau\vert \big\Vert\tilde{X}\big\Vert_4^4}{\sqrt{n}\,\sigma_X^2\big\Vert\tilde{X}\epsilon\big\Vert_2}
\end{align*}
and
\begin{align*}
\mathbb{P}\Bigg[ ab\left\vert
\overline{\X}-\mu_X\right\vert\left\vert\frac{1}{n}\sum_{i=1}^n{\varepsilon_i}\right\vert  \geq\frac{B}{144}\Bigg]
&\leq \frac{144ab}{B}\left\Vert
\overline{\X}-\mu_X\right\Vert_2\left\Vert\frac{1}{n}\sum_{i=1}^n{\varepsilon_i}\right\Vert_2 \leq \frac{576 \vert\tau\vert \big\Vert \tilde{X}\big\Vert_4^4\sigma_{\varepsilon}}{n\,\sigma_X \big\Vert \tilde{X}\epsilon\big\Vert_2^2}.
\end{align*}
The last two probabilities in (\ref{term_II_IIIprime}) are bounded by
 (\ref{majoration_terme3deII_III}) and 
(\ref{majoration_var_X}).
Hence, the term (\ref{II_IIIprime}) satisfies
\begin{multline}
\mathbb{P}\Bigg[ \vert\widehat{\tau}-\tau \vert \vert\tau \vert\frac{1}{n}\sum_{i=1}^n\tilde{X}_i^4 \geq\frac{B}{72} \Bigg] \leq \frac{1}{\sqrt{n}}\left(\frac{576\vert\tau\vert\big\Vert \tilde{X}\big\Vert_4^4}{\sigma_X^2 \big\Vert \tilde{X}\varepsilon\big\Vert_2}+ 4\frac{\big\Vert \tilde{X}\big\Vert_6^6}{\big\Vert \tilde{X}\big\Vert_4^6}
+2\frac{\big\Vert \tilde{X}\big\Vert_4^2}{\sigma_X^2}\right) + \frac{1}{n}\left(\frac{576  \vert\tau\vert\sigma_{\varepsilon}\big\Vert \tilde{X}\big\Vert_4^4 }{\sigma_X\big\Vert \tilde{X}\varepsilon\big\Vert_2^2}+2\right).\label{borne_II_IIIprime}
\end{multline}
\underline{Control of the term (\ref{III_Iprime}) of the upper bound of $D'_{2,2}$.}
Using the bound on 
 $\widehat{\tau}$ given by (\ref{majorationtauchapeau}) and Lemma \ref{lemmeproduit} with $a=4
 \big\Vert \tilde{X}^2 \tilde{Y} \big\Vert_1$, $b=\sqrt{2}\sigma_Y$ and $ c=\sqrt{2}/\sigma_X$, we get 
\begin{multline}
\mathbb{P}\Bigg[ \left\vert \frac{1}{n}\sum_{i=1}^n\tilde{X}_i^2\vert \widehat{\tau}\vert \vert\overline{\X}-\mu_{X}\vert\vert2Y_i-\overline{\Y}-\mu_Y\vert\right\vert \geq\frac{B}{72}\Bigg] \\
\leq \mathbb{P}\Bigg[ \vert\overline{\X}-\mu_X\vert\frac{1}{n}\sum_{i=1}^n{\tilde{X}_i^2\vert 2Y_i-\overline{\Y}-\mu_Y\vert} \frac{\sqrt{\frac{1}{n}\sum_{i=1}^n{(Y_i-\overline{\Y})^2}}}{\sqrt{\frac{1}{n}\sum_{i=1}^n{(X_i-\overline{\X})^2}}} \geq\frac{B}{72}\Bigg] \\
\leq \mathbb{P}\Bigg[ abc\vert\overline{\X}-\mu_X\vert \geq\frac{B}{72}\Bigg] 
+ \mathbb{P}\Bigg[ \frac{1}{n}\sum_{i=1}^n{\tilde{X}_i^2\vert 2Y_i-\overline{\Y}-\mu_Y\vert}  \geq a\Bigg] \\
+ \mathbb{P}\Bigg[ \frac{1}{n}\sum_{i=1}^n{(Y_i-\overline{\Y})^2} \geq b^2\Bigg]
+ \mathbb{P}\Bigg[ \frac{1}{n}\sum_{i=1}^n{(X_i-\overline{\X})^2} \leq 1/c^2\Bigg].\label{term_III_Iprime}
\end{multline}
The first probability in (\ref{term_III_Iprime}) is bounded by applying Markov Inequality:
\begin{align*}
\mathbb{P}\Bigg[ abc\vert\overline{\X}-\mu_X\vert\geq \frac{B}{72} \Bigg]
&\leq \frac{576 \big\Vert\tilde{X}^2\tilde{Y}\big\Vert_1\sigma_Y}{\sqrt{n}\big\Vert\tilde{X}\varepsilon\big\Vert_2^2\sigma_X}.
\end{align*}
The last three probabilities in (\ref{term_III_Iprime}) are bounded by
(\ref{majoration_terme2deIII_I}), (\ref{majoration_var_Y}) and (\ref{majoration_var_X}).  
Then the term (\ref{III_Iprime}) satisfies
\begin{multline}
\mathbb{P}\Bigg[\vert \overline{\X}-\mu_X\vert \vert\widehat{\tau}\vert\left\vert \frac{1}{n}\sum_{i=1}^n\tilde{X}_i^2 \vert 2Y_i-\overline{\Y}-\mu_Y\vert \right\vert \geq\frac{B}{72}\Bigg]\\
\leq \frac{1}{\sqrt{n}}\left(\frac{576 \big\Vert\tilde{X}^2\tilde{Y}\big\Vert_1\sigma_Y}{\big\Vert\tilde{X}\varepsilon\big\Vert_2^2\sigma_X}+ \frac{\big\Vert \tilde{X}^2\tilde{Y}\big\Vert_2}{\big\Vert \tilde{X}^2\tilde{Y}\big\Vert_1}
+ \frac{\sigma_Y \big\Vert \tilde{X}\big\Vert_4^2}{2\big\Vert \tilde{X}^2\tilde{Y}\big\Vert_1}
+ \frac{\big\Vert\tilde{Y}\big\Vert_4^2}{\sigma_Y^2}
+ \frac{2\big\Vert\tilde{X}\big\Vert_4^2}{\sigma_X^2}\right)+\frac{2}{n}.\label{borne_III_Iprime}
\end{multline}
\underline{Control of the term (\ref{III_IIprime}) of the upper bound of $D'_{2,2}$.}
Using bound (\ref{majorationtauchapeau}) on $\widehat{\tau}$ and Lemma \ref{lemmeproduit} with $a=2\big\Vert \tilde{X} \big\Vert_3^3$, $b=2\sigma_Y^2$ and $c=2/\sigma_X^2$, we get 
\begin{multline}
\mathbb{P}\Bigg[ \widehat{\tau}^2\vert\overline{\X}-\mu_{X}\vert \frac{1}{n}\sum_{i=1}^n\tilde{X}_i^2\vert X_i-\overline{\X}\vert \geq\frac{B}{72}\Bigg] 
\leq \mathbb{P}\Bigg[ \vert\overline{\X}-\mu_X\vert\frac{1}{n}\sum_{i=1}^n{\tilde{X}_i^2\vert X_i-\overline{\X}\vert} \frac{\frac{1}{n}\sum_{i=1}^n{(Y_i-\overline{\Y})^2}}{\frac{1}{n}\sum_{i=1}^n{(X_i-\overline{\X})^2}} \geq\frac{B}{72}\Bigg] \\
\leq \mathbb{P}\Bigg[ abc\vert\overline{\X}-\mu_X\vert \geq\frac{B}{72}\Bigg] 
+ \mathbb{P}\Bigg[ \frac{1}{n}\sum_{i=1}^n{\tilde{X}_i^2\vert X_i-\overline{\X}\vert}  \geq a\Bigg] \\
+ \mathbb{P}\Bigg[ \frac{1}{n}\sum_{i=1}^n{(Y_i-\overline{\Y})^2} \geq b\Bigg]
 + \mathbb{P}\Bigg[ \frac{1}{n}\sum_{i=1}^n{(X_i-\overline{\X})^2} \leq 1/c\Bigg]\label{term_III_IIprime}.
\end{multline}
The first probability in (\ref{term_III_IIprime}) is bounded using Markov Inequality:
\begin{align*}
\mathbb{P}\Bigg[ abc\vert\overline{\X}-\mu_X\vert\geq \frac{B}{72} \Bigg]
\leq \frac{576 \big\Vert \tilde{X}\big\Vert_3^3\sigma_Y^2}{\sqrt{n}\big\Vert \tilde{X}\varepsilon\big\Vert_2^2\sigma_X}.
\end{align*}
The last three probabilities in (\ref{term_III_IIprime}) are bounded by (\ref{majorationterme1-2_Markov}), (\ref{majoration_var_Y}) and (\ref{majoration_var_X}).
Hence, the term (\ref{III_IIprime}) satisfies
\begin{multline}
\mathbb{P}\Bigg[ \widehat{\tau}^2\vert\overline{\X}-\mu_{X}\vert \frac{1}{n}\sum_{i=1}^n\tilde{X}_i^2\vert X_i-\overline{\X}\vert \geq\frac{B\sqrt{\log(n)}}{36K\sqrt{n}}\Bigg] \\
\leq \frac{1}{\sqrt{n}}\left(
\frac{576 \big\Vert \tilde{X}\big\Vert_3^3\sigma_Y^2}{\big\Vert \tilde{X}\varepsilon\big\Vert_2^2\sigma_X}
+ \frac{\big\Vert \tilde{X}\big\Vert_6^3 +\sigma_X \big\Vert \tilde{X}\big\Vert_4^2}{\big\Vert \tilde{X}\big\Vert_3^3}
+ \frac{\big\Vert \tilde{Y}\big\Vert_4^2}{\sigma_Y^2}
+ \frac{2\big\Vert \tilde{X}\big\Vert_4^2}{\sigma_X^2}\right)+\frac{2}{n}.
\label{borne_III_IIprime}
\end{multline}
\underline{Control of the term (\ref{III_IIIprime}) of the upper bound of $D'_{2,2}$.}
Applying the upper bound on $\widehat{\tau}$ given by  (\ref{majorationtauchapeau}) and Lemma \ref{lemmeproduit} with 
$a=2\vert\tau\vert \big\Vert \tilde{X} \big\Vert_3^3$,  $b=\sqrt{2}\sigma_Y$ and $c=\sqrt{2}/\sigma_X$, we get 
\begin{multline}
\mathbb{P}\Bigg[ \vert \widehat{\tau}\vert \vert\overline{\X}-\mu_{X}\vert \vert\tau\vert\frac{1}{n}\sum_{i=1}^n\vert\tilde{X}_i\vert^3 \geq\frac{B}{72}\Bigg] \leq \mathbb{P}\Bigg[ \vert\overline{\X}-\mu_X\vert\vert\tau\vert\frac{1}{n}\sum_{i=1}^n{\vert\tilde{X}_i\vert^3} \frac{\sqrt{\frac{1}{n}\sum_{i=1}^n{(Y_i-\overline{\Y})^2}}}{\sqrt{\frac{1}{n}\sum_{i=1}^n{(X_i-\overline{\X})^2}}} \geq\frac{B}{72}\Bigg] \\
\leq \mathbb{P}\Bigg[ abc\vert\overline{\X}-\mu_X\vert \geq\frac{B}{72}\Bigg] 
+ \mathbb{P}\Bigg[ \vert\tau\vert\frac{1}{n}\sum_{i=1}^n{\vert\tilde{X}_i\vert^3}  \geq a\Bigg] \\
+ \mathbb{P}\Bigg[ \frac{1}{n}\sum_{i=1}^n{(Y_i-\overline{\Y})^2} \geq b^2\Bigg]
+ \mathbb{P}\Bigg[ \sqrt{\frac{1}{n}\sum_{i=1}^n{(X_i-\overline{\X})^2}} \leq 1/c\Bigg]. \label{term_III_IIIprime}
\end{multline}
The first probability in (\ref{term_III_IIIprime}) is bounded by Markov Inequality: 
\begin{align*}
\mathbb{P}\Bigg[ \left\vert\overline{\X}-\mu_X\right\vert\geq \frac{B}{72} \Bigg]
\leq 
\frac{288 \,\vert\tau\vert\big\Vert \tilde{X}\big\Vert_3^3\sigma_Y}{\sqrt{n}\big\Vert \tilde{X}\varepsilon\big\Vert_2^2}.
\end{align*}
The last three probabilities in (\ref{term_III_IIIprime}) are bounded by (\ref{majorationterme1-3_Markov}), (\ref{majoration_var_Y}) and (\ref{majoration_var_X}).
Consequently, the term (\ref{III_IIIprime}) satisfies
\begin{multline}
\mathbb{P}\Bigg[ \vert \widehat{\tau}\vert \vert\overline{\X}-\mu_{X}\vert \vert\tau\vert\frac{1}{n}\sum_{i=1}^n\vert\tilde{X}_i\vert^3 \geq\frac{B}{72}\Bigg] 
\leq \frac{1}{\sqrt{n}}\left(\frac{288 \vert\tau\vert\big\Vert \tilde{X}\big\Vert_3^3\sigma_Y}{\big\Vert \tilde{X}\varepsilon\big\Vert_2^2}+ \frac{\big\Vert \tilde{X}   \big\Vert_6^3
}{\big\Vert \tilde{X}   \big\Vert_3^3}
+ \frac{\big\Vert \tilde{Y}\big\Vert_4^2}{\sigma_Y^2}
+ \frac{2\big\Vert \tilde{X}\big\Vert_4^2}{\sigma_X^2}\right)+\frac{2}{n}. \label{borne_III_IIIprime}
\end{multline}
\underline{End of the control of the term $D'_{2,2}$.}
Combining Bounds  (\ref{borne_I_Iprime}), (\ref{borne_I_IIprime}), (\ref{borne_I_IIIprime}), (\ref{borne_II_Iprime}), (\ref{borne_II_IIprime}), (\ref{borne_II_IIIprime}), (\ref{borne_III_Iprime}), (\ref{borne_III_IIprime}) and (\ref{borne_III_IIIprime}), we get the following bound for 
the term $D'_{2,2}$ 
\begin{multline}
D'_{2,2} 
\leq 
\frac{1}{\sqrt{n}}\left[\frac{288\,\sigma_Y\big\Vert\tilde{X}^2 \tilde{Y}\big\Vert_1}{
\big\Vert\tilde{X}\varepsilon\big\Vert_2^2}+
\frac{864 \big\Vert \tilde{X}\big\Vert_3^3\sigma_Y^2}{\big\Vert \tilde{X}\varepsilon\big\Vert_2^2\sigma_X}+
\frac{432 \,\vert\tau\vert\big\Vert \tilde{X}\big\Vert_3^3\sigma_Y}{\big\Vert \tilde{X}\varepsilon\big\Vert_2^2}+
\frac{1152\big\Vert \tilde{X}^3\tilde{Y}\big\Vert_1}{\sigma_X^2\big\Vert\tilde{X}\varepsilon\big\Vert_2}+
\frac{1152 \big\Vert\tilde{X}\big\Vert_4^4\sigma_Y}{\sigma_X^3\big\Vert\tilde{X}\varepsilon \big\Vert_2} \right.\\
+\frac{576\vert\tau\vert \big\Vert\tilde{X}\big\Vert_4^4}{\sigma_X^2\big\Vert\tilde{X}\epsilon\big\Vert_2}+
\frac{576 \big\Vert\tilde{X}^2\tilde{Y}\big\Vert_1\sigma_Y}{\big\Vert\tilde{X}\varepsilon\big\Vert_2^2\sigma_X}+
\frac{2\big\Vert\tilde{X}^2 \tilde{Y}\big\Vert_2}{\big\Vert \tilde{X}^2\tilde{Y}\big\Vert_1}
+ \frac{\sigma_Y\big\Vert \tilde{X}\big\Vert_4^2}{\big\Vert \tilde{X}^2\tilde{Y}\big\Vert_1} 
+ \frac{2\sigma_X \big\Vert \tilde{X}\big\Vert_4^2}{\big\Vert \tilde{X}\big\Vert_3^3}
+ \frac{4\big\Vert \tilde{X}\big\Vert_6^3}{\big\Vert \tilde{X}\big\Vert_3^3}
 \\
\left. +\frac{5\big\Vert \tilde{Y}\big\Vert_4^2}{\sigma_Y^2}
+ 14\frac{\big\Vert \tilde{X}\big\Vert_4^2}{\sigma_X^2}+ 2^{7/2}\frac{\big\Vert\tilde{X}^3\tilde{Y}\big\Vert_{3/2}^{3/2}}{\big\Vert \tilde{X}^3 \tilde{Y}\big\Vert_1^{3/2}}+ \frac{\sigma_Y\big\Vert \tilde{X}\big\Vert_6^3}{2\big\Vert \tilde{X}^3 \tilde{Y}\big\Vert_1}
+8\frac{\big\Vert \tilde{X}\big\Vert_6^6}{\big\Vert \tilde{X}\big\Vert_4^6}
+ \frac{\sigma_X\big\Vert \tilde{X}\big\Vert_6^3}{2\big\Vert \tilde{X}\big\Vert_4^4}\right]\\
+ \frac{1}{n}\left[\frac{1152 \big\Vert\tilde{X}^3 \tilde{Y}\big\Vert_1
\sigma_{\varepsilon}}{
\sigma_X\big\Vert\tilde{X}\varepsilon\big\Vert_2^2
}+\frac{1152 \,\sigma_Y \sigma_{\varepsilon}\big\Vert \tilde{X}\big\Vert_4^4}{\big\Vert \tilde{X}\varepsilon\big\Vert_2^2\sigma_X^2}
+\frac{576  \vert\tau\vert\sigma_{\varepsilon}\big\Vert \tilde{X}\big\Vert_4^4 }{\sigma_X\big\Vert \tilde{X}\varepsilon\big\Vert_2^2}+14
\right]. 
\label{majorationD'22}
\end{multline}
\underline{End of the control of the term $D'_{2}$.}
Combining bounds (\ref{majorationD'21}) for $D'_{2,1}$ and (\ref{majorationD'22}) for $D'_{2,2}$, 
we get 
\begin{multline}
D'_{2} 
\leq 
\frac{1}{\sqrt{n}}\left[64\sqrt{2}\frac{\big\Vert\tilde{X}\varepsilon\big\Vert_3^3}{\big\Vert\tilde{X}\varepsilon\big\Vert_2^3}+\frac{288\,\sigma_Y\big\Vert\tilde{X}^2 \tilde{Y}\big\Vert_1}{
\big\Vert\tilde{X}\varepsilon\big\Vert_2^2}+
\frac{864 \big\Vert \tilde{X}\big\Vert_3^3\sigma_Y^2}{\big\Vert \tilde{X}\varepsilon\big\Vert_2^2\sigma_X}+
\frac{432 \,\vert\tau\vert\big\Vert \tilde{X}\big\Vert_3^3\sigma_Y}{\big\Vert \tilde{X}\varepsilon\big\Vert_2^2}+
\frac{1152\big\Vert \tilde{X}^3\tilde{Y}\big\Vert_1}{\sigma_X^2\big\Vert\tilde{X}\varepsilon\big\Vert_2}\right.\\
+
\frac{1152 \big\Vert\tilde{X}\big\Vert_4^4\sigma_Y}{\sigma_X^3\big\Vert\tilde{X}\varepsilon \big\Vert_2} 
+\frac{576\vert\tau\vert \big\Vert\tilde{X}\big\Vert_4^4}{\sigma_X^2\big\Vert\tilde{X}\epsilon\big\Vert_2}+
\frac{576 \big\Vert\tilde{X}^2\tilde{Y}\big\Vert_1\sigma_Y}{\big\Vert\tilde{X}\varepsilon\big\Vert_2^2\sigma_X}+
2\frac{\big\Vert\tilde{X}^2 \tilde{Y}\big\Vert_2}{\big\Vert \tilde{X}^2\tilde{Y}\big\Vert_1}
+ \frac{\sigma_Y\big\Vert \tilde{X}\big\Vert_4^2}{\big\Vert \tilde{X}^2\tilde{Y}\big\Vert_1} 
+ \frac{2\sigma_X \big\Vert \tilde{X}\big\Vert_4^2}{\big\Vert \tilde{X}\big\Vert_3^3}\\
\left.+
 \frac{4\big\Vert \tilde{X}\big\Vert_6^3}{\big\Vert \tilde{X}\big\Vert_3^3} 
+ \frac{5\big\Vert \tilde{Y}\big\Vert_4^2}{\sigma_Y^2} 
+ 14\frac{\big\Vert \tilde{X}\big\Vert_4^2}{\sigma_X^2}+ 2^{7/2}\frac{\big\Vert\tilde{X}^3\tilde{Y}\big\Vert_{3/2}^{3/2}}{\big\Vert \tilde{X}^3 \tilde{Y}\big\Vert_1^{3/2}}+ \frac{\sigma_Y\big\Vert \tilde{X}\big\Vert_6^3}{2\big\Vert \tilde{X}^3 \tilde{Y}\big\Vert_1}
+8\frac{\big\Vert \tilde{X}\big\Vert_6^6}{\big\Vert \tilde{X}\big\Vert_4^6}
+ \frac{\sigma_X\big\Vert \tilde{X}\big\Vert_6^3}{2\big\Vert \tilde{X}\big\Vert_4^4}\right]\\
+ \frac{1}{n}\left[\frac{1152 \big\Vert\tilde{X}^3 \tilde{Y}\big\Vert_1
\sigma_{\varepsilon}}{
\sigma_X\big\Vert\tilde{X}\varepsilon\big\Vert_2^2
}+\frac{1152 \,\sigma_Y \sigma_{\varepsilon}\big\Vert \tilde{X}\big\Vert_4^4}{\big\Vert \tilde{X}\varepsilon\big\Vert_2^2\sigma_X^2}
+\frac{576  \vert\tau\vert\sigma_{\varepsilon}\big\Vert \tilde{X}\big\Vert_4^4 }{\sigma_X\big\Vert \tilde{X}\varepsilon\big\Vert_2^2}+14
\right]. 
\label{majorationD'2}
\end{multline}
\subsubsection*{End of the proof of Lemma \ref{lemme2} }
Using the bounds (\ref{majorationD'1}) for $D'_1$ and (\ref{majorationD'2}) for $D'_2$ and the definition of the Lyapunov ratios $G_1$, $G_2$ and $G_3$ given by (\ref{termeG1}), (\ref{termeG2}) and (\ref{termeG3}) in Section \ref{definitionG_ietH_i}, we get 
\begin{align*}
\mathbb{P} \Bigg[   \left\vert   \frac{1}{n}\sum_{i=1}^n(X_i-
\overline{\X})^2{\widehat{\varepsilon_i}}^2 - \big\Vert\tilde{X}\varepsilon\big\Vert_2^2
\right\vert \geq\frac{B}{2}     \Bigg] \leq 
\frac{G_1(X,Y)}{\sqrt{n}}+ \frac{4G_2(X,Y)}{n}+\frac{G_3(X,Y)}{n^{3/2}}
\end{align*}
which completes the proof of Lemma \ref{lemme2}. 

\bigskip




\noindent \textbf{Acknowledgements.}
We would like to thank Marie-Laure Martin for her expertise on the RNA-seq data in Section \ref{dataset}. Her help enabled us to implement our selection procedure on these data.


\begin{thebibliography}{00}{

\bibitem{BE1965}
B. von Bahr and C.G.  Esseen, (1965) Inequalities for the $r$-th absolute moment of a sum of random variables, $1\leq r \leq 2$. {\em Ann. Math. Statist.} {\bf 36}, pp: 299–303. 

\bibitem{BG}
V. Bentkus and F. Götze (1996) The Berry-Esseen bound for Student's statistic. {\em Ann. Probab.} {\bf 24(1)}, pp: 491–503. 

\bibitem{Cai}
T. Cai, Z. Guo and R. Ma (2023) Statistical inference for high-dimensional generalized linear models with binary outcomes. {\em Journal of the American Statistical Association} {\bf 118(542)}, pp: 1319-1332. 

\bibitem{CDMP2022} C.
Cuny, and J. Dedecker and F. Merlevède and M. Peligrad (2022) Berry–Esseen type bounds for the matrix coefficients and the spectral radius of the left random walk on $GL_d({\protect \mathbb{R}})$. {\em C. R. Math. Acad. Sci. Paris}, {\bf 360 (G5)}, pp: 475-482. 


\bibitem{DGT2025} J. Dedecker and O. Guedj and M.L. Taupin, M. L. (2025) Asymptotic confidence interval for $R^2$ in multiple linear regression. {\em Statistics} {\bf 59(1)}, pp: 1–36. 

\bibitem{Esseen1942} C.- G. Esseen (1942), 
On the Liapounoff limit of error in the theory of probability.
{\em Ark. Mat. Astr. Fys. } {\bf 28A (9)}.


\bibitem{Fan2008} J. Fan and J. Lv, (2008). Sure independence screening for ultrahigh dimensional feature space. {\em J. R. Stat. Soc. Ser. B Stat. Methodol.: Series B}, {\bf 70},  pp: 849–911.

\bibitem{ScreeningGLM2010} J. Fan and R. Song, (2010). Sure independence screening in generalized linear models with NP-dimensionality. {\em Ann. Statist.}, {\bf 38 (6)},  pp: 3567-3604.

\bibitem{FSW2009} J. Fan, R. Samworth and Y. Wu, (2009).
Ultrahigh dimensional feature selection: beyond the linear model. 
{\em J. Mach. Learn. Res.} {\bf 10}, 2013–2038.

\bibitem{FN}D. H. Fuk, (1973). Certain probabilistic inequalities for martingales. {\em (Russian) Sibirsk. Mat. Ž. } {\bf 14}, pp: 185–193.

\bibitem{Hou2024} X. Hou and L. Zhang and P. Wang and M. Xie (2024) Repro Samples Method for High-dimensional Logistic Model, 
\textit{{https://arxiv.org/abs/2403.09984}.}

\bibitem{Kaetzel2005} C. S. Kaetzel, (2005).  The polymeric immunoglobulin receptor: bridging innate and adaptive immune responses at mucosal surfaces. {\em Immunol Rev.}, {\bf 206}, pp: 83-99. 


\bibitem{Panetal2019}
W. Pan, X. Wang, W. Xiao, H. Zhu, (2019). A Generic Sure Independence Screening Procedure. {\em J. Am. Stat Assoc.}, {\bf 114(526)}, pp: 928-937. 

\bibitem{Pinelis1994} I. Pinelis (1994)
Optimum bounds for the distributions of martingales in Banach spaces. {\em 
Ann. Probab.} {\textbf 22 no. 4}, pp: 1679–1706.

\bibitem{Pinelis2015} I. Pinelis, (2015) Best possible bounds of the von Bahr–Esseen type. {\em Ann. Funct. Anal.} {\bf 6(4)}, pp: 1–29. 

\bibitem{Pollak1956}
H. O. Pollak, (1956) A remark on ``Elementary inequalities for Mills' ratio'' by Yûsaku Komatu. {\em Rep. Statist. Appl. Res. Un. Japan. Sci. Engrs.} {\bf 4 }. 

\bibitem{Rosenthal1970}
H. P. Rosenthal, (1970) On the subspaces of $L_p$ ($p>2$) spanned by sequences of independent random variables. {\em Israel J. Math.} {\bf 8}, pp: 273–303. 


\bibitem{Shalek2014} A. K. Shalek and  R. Satija and J. Shuga, J. J. Trombetta and D. Gennert and D.  Lu and P.
Chen and  R. S. Gertner, and J. T. Gaublomme and N. Yosef, et al. (2014) Single-cell RNA-seq reveals
dynamic paracrine control of cellular variation. {\em Nature}, {\bf 510(7505)} pp: 363–369.

\bibitem{Simeone2005} M. Simeone and S.A.  Phelan, (2005). Transcripts associated with Prdx6 (peroxiredoxin 6) and related genes in mouse. {\em Mamm Genome}  {\bf 16(2)}, pp: 103-111. 

\bibitem{W}  H. White, (1980).  A heteroskedasticity-consistent covariance matrix estimator and a direct test for heteroskedasticity. {\em Econometrica}, {\bf 48(4)}, pp: 817-838.
}
\bibitem{CoxScreening2012} S. D. Zhao and Y. Li, (2012) Principled sure independence screening for Cox models with ultra-high-dimensional covariates. {\em J. Multivariate Anal.}, {\bf 105},  pp: 397–411.
\end{thebibliography}
\end{document}